\documentclass{amsart}
\usepackage{graphicx}
\usepackage{amssymb}
\usepackage{epstopdf}
\usepackage{nicefrac}
\usepackage{enumerate}
\usepackage{float}
\usepackage{color}
\usepackage{pst-all}
\usepackage{tabularx}
\usepackage{hyperref}
\usepackage{caption}
\usepackage{subcaption}

\newtheorem{thm}{THEOREM}[section]

\newtheorem{cor}[thm]{COROLLARY}
\newtheorem{defn}[thm]{DEFINITION}

\newtheorem{hyp}[thm]{HYPOTHESIS}
\newtheorem{lemma}[thm]{LEMMA}

\newtheorem{prop}[thm]{PROPOSITION}

\newtheorem{remark}[thm]{REMARK}

\newcommand{\ds}{\displaystyle}

\newcommand{\G}{\Gamma}
\newcommand{\e}{{\epsilon}}
\newcommand{\g}{{\gamma}}
\newcommand{\wmW}{\widehat{\mW}_\e}
\newcommand{\bRt}{{\bf R}_{0}}
\newcommand{\bRa}{{\bf R}_{a}}

\newcommand{\bR}{{\bf R}}

\newcommand{\mK}{{\mathbb K}}
\newcommand{\mN}{{\mathbb N}}
\newcommand{\mR}{{\mathbb R}}
\newcommand{\mS}{{\mathbb S}}

\newcommand{\mW}{{\mathbb W}}
\newcommand{\cA}{{\mathcal A}}
\newcommand{\cC}{{\mathcal C}}
\newcommand{\cD}{{\mathcal D}}

\newcommand{\cG}{{\mathcal G}}
\newcommand{\cI}{{\mathcal I}}
\newcommand{\cK}{{\mathcal K}}
\newcommand{\cL}{{\mathcal L}}

\newcommand{\cO}{{\mathcal O}}
\newcommand{\cP}{{\mathcal P}}
\newcommand{\cR}{{\mathcal R}}

\newcommand{\cW}{{\mathcal W}}
\newcommand{\cX}{{\mathcal X}}
\newcommand{\cZ}{{\mathcal Z}}

\newcommand{\whcGe}{\widehat{\cG}_{\epsilon}}

\newcommand{\whPsi}{{\widehat \Psi}}
\newcommand{\whPhi}{{\widehat{\Phi}}}
\newcommand{\whPhie}{\widehat{\Phi^\e}}
\newcommand{\whW}{{\widehat W}}

\newcommand{\fD}{{\mathfrak{D}}}
\newcommand{\fM}{{\mathfrak{M}}}

\newcommand{\ovG}{\overline{\Gamma}}

\newcommand{\oq}{\overline{q}}

\newcommand{\otheta}{\overline{\theta}}
\newcommand{\ovg}{\overline{\gamma}}
\newcommand{\ovl}{\overline{\lambda}}
\newcommand{\ovk}{\overline{\kappa}}
\newcommand{\ovtg}{\overline{\widetilde{\gamma}}}
\newcommand{\ocP}{\overline{\cP}}

 \newcommand{\vp}{{\varphi}}

\newcommand{\psg}{{\rm pseudo}{\star}{\rm group}}

 \newcommand{\wtg}{\widetilde{\gamma}}
 \newcommand{\wtr}{\widetilde{r}}

\begin{document}

\begin{abstract}
We consider the dynamical properties of $C^{\infty}$-variations of the flow on an  aperiodic Kuperberg plug $\mK$. 
Our main result is that there exists a smooth 1-parameter family of plugs $\mK_\e$ for $\e\in (-a,a)$ and $a<1$, such that:
(1) The plug $\mK_0 = \mK$ is a generic Kuperberg plug;
(2) For $\e<0$,  the flow in the plug $\mK_\e$ has two periodic orbits that
  bound an invariant cylinder, all other orbits of the flow are wandering, and the flow has topological entropy  zero;
(3) For $\e>0$, the flow in the plug $\mK_\e$ has positive topological entropy, and an abundance of periodic orbits.
 \end{abstract}

\title{Aperiodicity at the boundary of chaos}


\author{Steven Hurder}
\address{Steven Hurder, Department of Mathematics, University of Illinois at Chicago, 322 SEO (m/c 249), 851 S. Morgan Street, Chicago, IL 60607-7045}
\email{hurder@uic.edu}
\thanks{Preprint date: March 25, 2016; revised September 25}

\author{Ana Rechtman}
\address{Ana Rechtman, Instituto de Matem\'aticas, Universidad Nacional Aut\'onoma de M\'exico, Ciudad Universitaria, 04510 Ciudad de M\'exico, Mexico}
\email{rechtman@im.unam.mx}

\date{}


\maketitle

\section{Introduction} \label{sec-intro}

In this paper, we  analyze the dynamical properties of flows    in a $C^{\infty}$-neighborhood
of the    Kuperberg flows introduced in \cite{Kuperberg1994}, or to be more precise, of  {\it generic
  Kuperberg flows} as introduced in \cite{HR2016}.   The Kuperberg flows are   exceptional for the simplicity of their   explicit construction in   \cite{Kuperberg1994}, 
 and this explicitness makes it a straightforward process to construct 1-parameter families of $C^\infty$-deformations  of a given 
generic Kuperberg flow.
 We show in this work that the Kuperberg flows are furthermore  remarkable, in that there are $C^{\infty}$-nearby    flows with simple   dynamics, and that there are $C^{\infty}$-nearby   flows with positive topological
entropy and an abundance of periodic orbits.

 The construction of a Kuperberg flow is based on the construction of an aperiodic
 plug, which we call a \emph{Kuperberg Plug}, and is noted in this paper by
 $\mK_0$. A plug is a manifold with boundary endowed with a flow, that
 enables  the modification of a given flow inside a flow-box, so that
 after modification, there are orbits that enter the flow-box and never
 exit. Moreover,   Kuperberg's construction  does this modification without introducing additional  periodic orbits.
  Expository treatments of Kuperberg's construction were given by Ghys  \cite{Ghys1995} and Matsumoto \cite{Matsumoto1995}, and in the first chapters of the authors' work \cite{HR2016}.   Our work also introduced new concepts for the study of the dynamical properties of the Kuperberg flows, which allows one to investigate many further remarkable aspects of these flows.

As a consequence of Katok's theorem on $C^2$-flows on $3$-manifolds \cite{Katok1980},
 the topological entropy of a Kuperberg flow is zero.   In
\cite{HR2016} we developed a technique based on the introduction of an ``almost transverse'' 
rectangle inside the plug and a
pseudogroup modeling the return map of the flow to the rectangle, to make an explicit
computation of the topological entropy. This computation revealed
  chaotic behavior for the flow, but such that it evolves    at a very slow rate, and thus it does not
result in positive  topological entropy. We proved that, under some
extra hypotheses, such generic Kuperberg   flows have positive  ``slow entropy'', and that  the rates of chaotic behavior grow at a precise subexponential, but non-polynomial  rate,    as  discussed in the proof of
\cite[Theorem~21.10]{HR2016}. The calculation behind the proof of this result   suggests that   for some
Kuperberg-like flows near to a generic Kuperberg flow, there should be actual chaotic behavior, and  also positive entropy.

In this paper we   prove the following two theorems which make these remarks more precise.

\begin{thm}\label{thm-main1}
There exists a $C^{\infty}$ 1-parameter family of plugs $\mK_\e$ for $\e\in (-1,0]$ such that:
\begin{enumerate}
\item The plug $\mK_0$ is a   Kuperberg plug;
\item For $\e<0$, the flow in the plug $\mK_\e$ has two periodic orbits that
  bound an invariant cylinder, and every other orbit belongs to the wandering set, and thus the flow has topological entropy zero.
\end{enumerate}
\end{thm}
The proof of Theorem~\ref{thm-main1} uses the same  technical tools as developed in the previous works \cite{Kuperberg1994,Kuperbergs1996,Ghys1995,Matsumoto1995,HR2016} for the study of the dynamics of Kuperberg flows, and the result is notable mainly for its contrast with the following result,   that a $C^{\infty}$-neighborhood of a generic Kuperberg flow also contains flows which have exceptionally wild dynamics.  
In particular, we obtain the following result:
\begin{thm}\label{thm-main2}
There exists a $C^{\infty}$ 1-parameter family of plugs $\mK_\e$ for $\e\in [0,a)$, $a > 0$, 
  such that:
\begin{enumerate}
\item The plug $\mK_0$ is a generic Kuperberg plug;
\item For $e>0$, the flow in $\mK_\e$ has positive topological entropy, and an abundance of periodic orbits.
\end{enumerate}
\end{thm}
The proof of Theorem~\ref{thm-main2}  is based on the understanding of the dynamics of standard Kuperberg flows developed in \cite{HR2016}, and in particular uses in a fundamental way  the   technique of relating the dynamics of a Kuperberg-like flow to the dynamics of its return map to an almost transverse rectangle.

The construction of the plugs $\mK_\e$ in the proofs of both Theorems~\ref{thm-main1} and \ref{thm-main2} follows closely the original
construction by K. Kuperberg in \cite{Kuperberg1994}, which begins with a modified version of the original Wilson Plug   \cite{Wilson1966},
 where the modification given in Section~\ref{subsec-wilson} removes the stable periodic orbits for the flow and replaces them with unstable periodic orbits.
The  construction of the 1-parameter
family of plugs is described in the first part of this work,
Section~\ref{subsec-kuperberg}, and closely follows the construction in \cite{HR2016}. 

The   \emph{Radius Inequality} is the main condition in Kuperberg's construction that is used  to prove that the flow is 
aperiodic. The only change in the requirements of the construction of the flows $\Phi_t^\e$ which we study, is that the 
Radius Inequality    gets replaced for $\e \ne 0$ with  the \emph{Parametrized Radius Inequality}, as stated
in Section~\ref{subsec-kuperberg}. For $\e=0$ we recover the
original Radius Inequality for Kuperberg flows.

  Section~\ref{sec-radiuslevel} introduces some of the main tools
for the study of the dynamics of Kuperberg flows, the radius and level functions, and also gives 
some of the immediate consequences for the study of the dynamics that
are independent of the value of $\e$ in the Parametrized Radius Inequality.

As we mentioned above, one of the main techniques in \cite{HR2016} is
to introduce a pseudogroup modeling the dynamics of the Kuperberg
flow. The other main technique is the study of surfaces tangent to the
flow, known as propellers,  that were used to describe the topological properties of the
minimal set of a generic Kuperberg flow. These surfaces are
introduced here in Section~\ref{sec-propeller}.
 
The plugs $\mK_\e$ for $\e<0$ have rather simple dynamical properties, as stated in
Theorem~\ref{thm-main1} and described in Section~\ref{sec-enegative}. The study of
these flows does not requires any extra hypothesis in the
constructions.

On the other hand, the study of the dynamical properties of plugs $\mK_\e$  which satisfy the Parametrized Radius Inequality for   $\e \geq 0$ is extraordinarily complicated. The detailed analysis  in \cite{HR2016} of the standard Kuperberg flows  
 required the introduction of the {\it generic hypothesis}   in that work (see Hypothesis~\ref{hyp-SRI} below), in order to deduce a wide range of properties for the flows. For the case when $\e > 0$, many of the corresponding results do not hold, so we consider in this work only the question of the existence of compact invariant sets for the flows, such that the restricted dynamics of the flow admits a ``horseshoe'' in its transversal model. Even with this more restricted goal, it is still necessary to impose   geometric hypothesis on the
construction as stated in Hypotheses~\ref{hyp-monotone} and \ref{hyp-offset}, in order to obtain 
  the proof of Theorem~\ref{thm-main2}.

  The construction of the pseudogroup
acting on an almost transverse rectangle $\bRt\subset \mK_\e$ is done in
Section~\ref{sec-pseudogroups}. Instead of introducing a pseudogroup analogous to the
one in Chapter 9 of \cite{HR2016}, we introduce the simplest
pseudogroup that allows us to prove
Theorem~\ref{thm-main2}. Section~\ref{sec-horseshoe} is dedicated to the proof that this 
pseudogroup contains ``horseshoe maps'', and the proof introduces surfaces analogous to the 
propellers introduced in \cite{HR2016} which are used to define these   maps. In Section~\ref{sec-entropy} we describe the situation in which these horseshoe maps can be embedded
in the flow $\Phi_t^\e$, so that they generate positive entropy for the flow, and
not just for the pseudogroup.   Then in Section~\ref{subsec-admissible}, we discuss the construction of  examples of $C^{\infty}$-deformations of a generic Kuperberg flow, such that the hypotheses of Theorem~\ref{thm-entropypositive} are satisfied,   completing the proof of Theorem~\ref{thm-main2}.

The results in this paper are  inspired by the work in the paper \cite{HR2016}, and for the sake of brevity, we are forced to refer occasionally to
results in \cite{HR2016}. However, the novel results used in the proof
of Theorems~\ref{thm-main1} and \ref{thm-main2} are explained and proved here.
Further questions and open problems concerning the
 dynamical properties of the many variations of   Kuperberg flows are 
 discussed    in the paper \cite{HKR2016}.

\section{Construction of   parametrized families}\label{sec-plugs}

In this section, we present the construction of the   1-parameter family of plugs $\mK_\e$ for $\e\in (-a,a)$, and introduce some of the techniques used for the study of their dynamical properties. 
This follows closely the outline of the construction and study of the
usual Kuperberg flows in Chapters $2$ and $3$ of \cite{HR2016}.

In Section~\ref{subsec-wilson} we give the construction of the modified Wilson
plug as introduced  by Kuperberg in \cite{Kuperberg1994}.
 The family of plugs is constructed in Section~\ref{subsec-kuperberg}, and the Parametrized Radius Inequality is introduced.

 \subsection{Plugs}\label{subsec-plugs}

A $3$-dimensional plug is a manifold $P$ endowed with a vector field $\cX$ satisfying the following ``plug conditions''. The 3-manifold $P$ is of the form $D \times [-2,2]$, where $D$ is a compact 2-manifold with boundary $\partial D$. Set 
$$\partial_v P = \partial D \times [-2,2] \quad , \quad \partial_h^- P = D \times \{-2\} \quad , \quad \partial_h^+ P = D \times \{2\} \ .$$
Then   the boundary  of $P$ has a decomposition
$$\partial P ~ = ~  \partial_v P \cup \partial_h P ~ = ~  \partial_v P \cup \partial_h^-P \cup \partial_h^+ P \ .$$
Let $\frac{\partial}{\partial z}$ be the \emph{vertical} vector field on $P$, where $z$ is the coordinate of the interval $[-2,2]$.

The vector field $\cX$ must satisfy the conditions:
\begin{itemize}
\item[(P1)] \emph{vertical at the boundary}: $\cX=\frac{\partial}{\partial z}$ in a neighborhood of $\partial P$; thus, $\partial_h^- P$ and $\partial_h^+ P$  are the entry and exit regions of $P$ for the flow of $\cX$, respectively; 
\item[(P2)]   \emph{entry-exit condition}: if a point $(x,-2)$ is in the same trajectory as $(y,2)$, then $x=y$. That is,  an orbit that traverses $P$, exits just in front of its entry point;
\item[(P3)] \emph{trapped orbit}: there is at least one entry point whose entire \emph{forward} orbit is contained in $P$; we will say that its orbit is \emph{trapped} by $P$;
\item[(P4)] \emph{tame}: there is an embedding $i \colon P\to \mR^3$ that preserves the vertical direction.
\end{itemize}

   Note that conditions (P2) and (P3) imply that if the forward orbit of a point $(x,-2)$ is trapped, then the backward orbit of $(x,2)$ is also trapped. 

A {\it semi-plug} is  a manifold $P$ endowed with a vector field $\cX$ as above, satisfying conditions (P1), (P3) and (P4), but not necessarily (P2).  The concatenation of a semi-plug with an inverted copy of it, that is a copy where the direction of the flow is inverted, is then a plug. 

Note that condition   (P4)  implies that given any open ball  $B(\vec{x},\e) \subset \mR^3$ with $\e > 0$, there exists a modified embedding $i' \colon P \to B(\vec{x},\e)$ which preserves the vertical direction again.   Thus, a plug can be used to change a vector field $\cZ$ on any $3$-manifold $M$ inside a flowbox, as follows. Let $\vp \colon U_x \to (-1,1)^3$ be a coordinate chart which maps the vector field $\cZ$ on $M$ to the vertical vector field $\frac{\partial}{\partial z}$. Choose a modified embedding $i' \colon P \to B(\vec{x},\e) \subset (-1,1)^3$, and then  replace the flow $\frac{\partial}{\partial z}$ in the interior of $i'(P)$ with the image of $\cX$. This results in a   flow $\cZ'$ on $M$.

 The   entry-exit condition implies that    a  periodic orbit of $\cZ$ which meets $\partial_h P$ in a non-trapped point,   will remain periodic after this modification. An orbit of $\cZ$ which meets $\partial_h P$ in a  trapped point  never exits the plug $P$, hence  after modification,   limits to a closed invariant   set contained in $P$.  A closed invariant set contains a minimal set for the flow, and thus, a plug  serves as a device    to insert a minimal set  into a flow.

\subsection{The modified Wilson plug $\mW$}\label{subsec-wilson}

  We next introduce the    ``modified Wilson Plug'', which was 
 the starting point of Kuperberg's construction of her plug. We add to the construction one 
 hypothesis  that seems to be implicitly   assumed by   certain
 conclusions stated   in  \cite{Ghys1995,Matsumoto1995}, and that are
 part of the generic hypothesis used in \cite{HR2016}. This
 hypothesis is not needed to construct Kuperberg's aperiodic plug,
 but it seems necessary in order to prove the results in this paper. 
    
\begin{figure}[!htbp]
\centering
 \includegraphics[width=40mm]{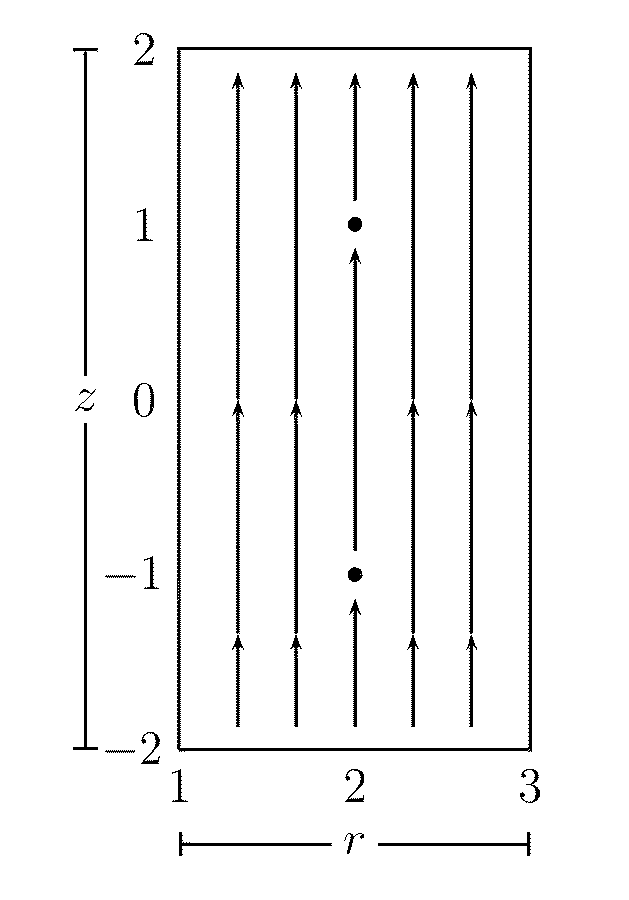}
\caption{Vector field $\cW_v$ \label{fig:wilson1}}
\vspace{-6pt}
\end{figure}

Consider the rectangle, as illustrated in Figure~\ref{fig:wilson1}, 
\begin{equation}\label{eq-rectangle}
{\bR} = [1,3]\times[-2,2] = \{(r,z) \mid 1 \leq r \leq 3 ~ \& -2 \leq z \leq 2\} \ .
\end{equation}
 For a constant $0 < g_0\leq 1$,  choose a $C^\infty$-function $g \colon \bR \to [0,g_0]$   which satisfies the ``vertical'' symmetry condition $g(r,z) = g(r,-z)$. 
Also, require that   $g(2,-1) = g(2,1) = 0$,   that $g(r,z) = g_0$ for $(r,z)$ near the boundary of $\bR$, and that $g(r,z) > 0$ otherwise. We may take $g_0 = 1/10$ for example.
 
We make the following additional hypothesis on the function $g$ in the construction.
\begin{hyp}\label{hyp-genericW}
We require the function $g$ satisfy the above conditions, and in addition:
\begin{equation}\label{eq-generic1}
g(r,z) = g_0 \quad \text{for} \quad  (r-2)^2 + (|z| -1)^2 \geq \e_0^2 
\end{equation}
where $0 < \e_0 < 1/4$ is sufficiently small, as will be specified later in Section~\ref{subsec-kuperberg}. In addition, we also require   that $g(r,z)$ is monotone increasing as a function of the distance $\ds \sqrt{(r-2)^2 + (|z| -1)^2}$ from the special points where it vanishes, and that   $g$  is  non-degenerate.  Non-degenerate means that the matrix of second partial derivatives at each vanishing point is invertible.
\end{hyp}

Define the vector field $\cW_v = g \cdot \frac{\partial}{\partial  z}$ which has two singularities, $(2,\pm 1)$, and   is otherwise everywhere vertical. The flow lines of this vector field are    illustrated in Figure~\ref{fig:wilson1}. 

Next, choose a $C^\infty$-function $f \colon \bR \to [-1,1]$ which satisfies the following conditions:
\begin{enumerate}
\item[(W1)] $f(r,-z) = -f(r, z)$ ~ [\emph{anti-symmetry in z}]
\item[(W2)] $f(r,z) = 0$ for $(r,z)$ near the boundary of $\bR$
\item[(W3)] $f(r,z) \geq 0$ for  $-2 \leq z \leq 0$. 
\item[(W4)] $f(r,z) \leq 0$ for  $0  \leq z \leq 2$. 
\item[(W5)] $f(r,z) =1$ for $5/4 \leq r \leq 11/4$  ~ \text{and} ~  $-7/4 \leq z \leq -1/4$. 
\item[(W6)] $f(r,z) = -1$ for $5/4 \leq r \leq 11/4$ ~ \text{and}~  $1/4 \leq z \leq 7/4$. 
\end{enumerate}
Condition (W1) implies that $f(r,0) =0$ for all $1 \leq r \leq 3$.

  Next, define the manifold with boundary
\begin{equation}\label{eq-wilsoncylinder}
\mW=[1,3] \times \mS^1\times[-2,2] \cong {\mathbf R} \times \mS^1
\end{equation}
with cylindrical coordinates  $x = (r, \theta,z)$.  That is, 
 $\mW$ is a      solid cylinder with an open core removed, obtained by rotating  the rectangle $\bR$,  considered as embedded in $\mR^3$, around the $z$-axis. 
 
\begin{figure}[!htbp]
\centering
\begin{subfigure}[c]{0.3\textwidth}{\includegraphics[height=28mm]{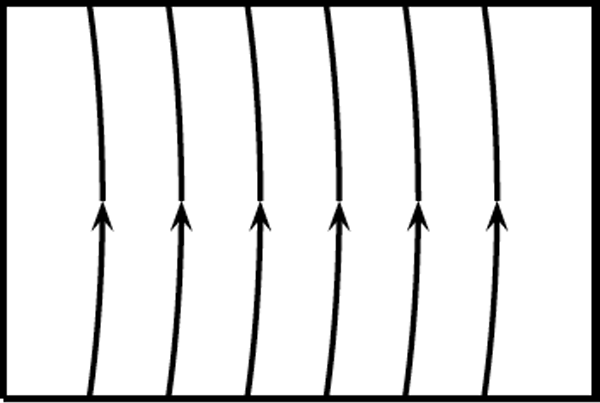}}\caption{$r \approx 1,3$}\end{subfigure}
\begin{subfigure}[c]{0.3\textwidth}{\includegraphics[height=28mm]{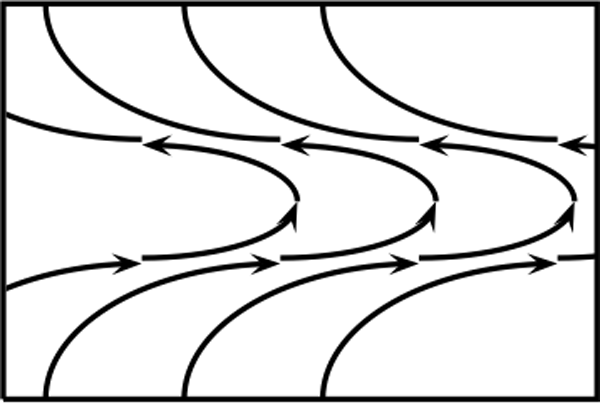}}\caption{$r \approx 0$}\end{subfigure}
\begin{subfigure}[c]{0.3\textwidth}{\includegraphics[height=28mm]{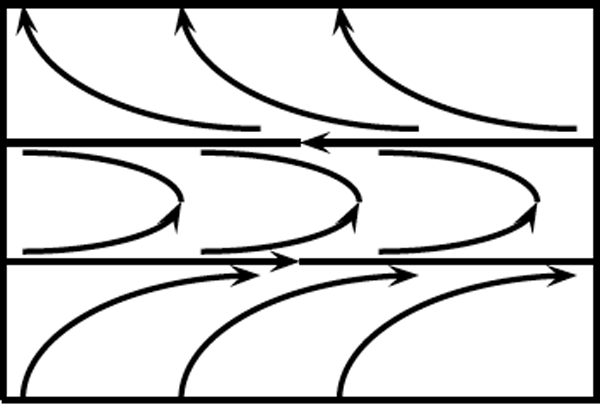}}\caption{$r =  0$}\end{subfigure}
\caption{$\cW$-orbits on the cylinders $\{r=const.\}$ \label{fig:flujocilin}}
\vspace{-6pt}
\end{figure}

\begin{figure}[!htbp]
\centering
\begin{subfigure}[c]{0.4\textwidth}{\includegraphics[height=60mm]{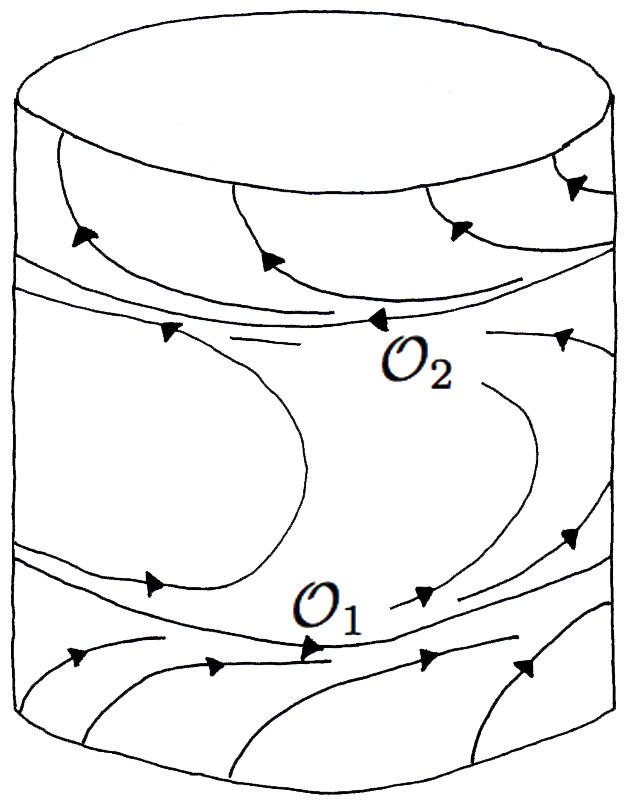}}\end{subfigure}
\begin{subfigure}[c]{0.4\textwidth}{\includegraphics[height=60mm]{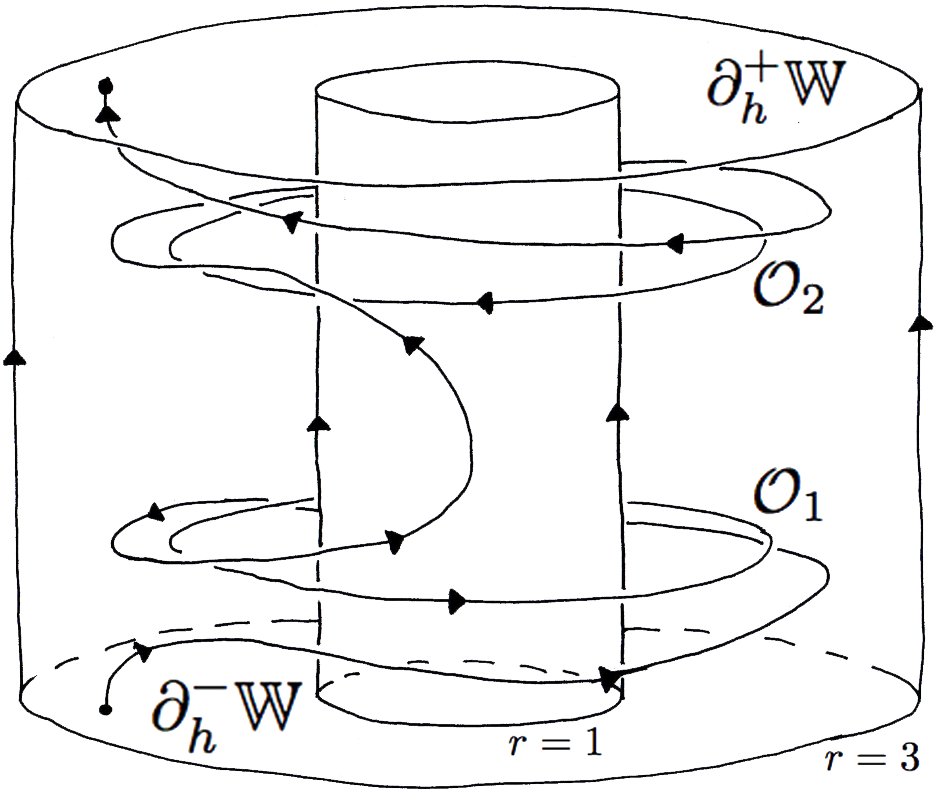}}\end{subfigure}
 \caption[{$\cW$-orbits in the cylinder  $\cC=\{r=2\}$}]{$\cW$-orbits in the cylinder  $\cC=\{r=2\}$ and in $\mW$ \label{fig:Reebcyl}}  
\vspace{-6pt}
\end{figure}

Extend the functions $f$ and $g$ above to $\mW$ by setting $f(r, \theta, z) = f(r, z)$ and $g(r, \theta, z) = g(r, z)$, so that they are   invariant under rotations around the $z$-axes. 
Define the Wilson vector field on $\mW$  by 
\begin{equation}\label{eq-wilsonvector}
\cW =g(r, \theta, z)  \frac{\partial}{\partial  z} + f(r, \theta, z)  \frac{\partial}{\partial  \theta}
\end{equation}
Let $\Psi_t$ denote the flow of $\cW$ on $\mW$.  
Observe that the vector field $\cW$ is vertical near the boundary of $\mW$ and horizontal in the periodic orbits. 
Also, $\cW$ is tangent to the cylinders $\{r=cst\}$. The flow of $\Psi_t$ on the cylinders $\{r=cst\}$ is illustrated (in cylindrical coordinate slices) by  the lines in  Figures~\ref{fig:flujocilin} and \ref{fig:Reebcyl}.

Define the closed subsets:
\begin{eqnarray*}
\cC ~ & \equiv & ~  \{r=2\}   \quad \text{[\emph{The ~ Full ~ Cylinder}]}\\
\cR ~ & \equiv & ~  \{(2,\theta,z) \mid  -1 \leq z \leq 1\} \quad \text{[\emph{The ~ Reeb ~ Cylinder}]}\\
\cA ~ & \equiv & ~  \{z=0\} \quad \text{[\emph{The ~ Center ~ Annulus}]}\\
\cO_i ~ & \equiv & ~  \{(2,\theta,(-1)^i)  \} \quad \text{[\emph{Periodic Orbits, i=1,2}]}
\end{eqnarray*}
Note that  $\cO_1$ is the lower   boundary circle  of the Reeb cylinder $\cR$, and $\cO_2$ is the upper boundary circle.

We give some of the basic properties of the Wilson flow.  
Let $R_{\varphi} \colon \mW \to \mW$ be rotation by the angle $\varphi$. That is, $R_{\varphi}(r,\theta,z) = (r, \theta + \varphi, z)$.  
\begin{prop}\label{prop-wilsonproperties}
Let $\Psi_t$ be the flow on $\mW$ defined above, then:
\begin{enumerate}
\item $R_{\varphi} \circ \Psi_t = \Psi_t \circ R_{\varphi}$ for all $\varphi$ and $t$.
\item The flow $\Psi_t$ preserves the cylinders $\{r=const.\}$   and  in particular 
    preserves the cylinders $\cR$ and $\cC$. 
\item $\cO_i$ for $i=1,2$ are  the periodic orbits for $\Psi_t$.
\item For $x = (2,\theta,-2)$, the forward orbit   $\Psi_t(x)$ for $t > 0$ is trapped.
\item For $x = (2,\theta,2)$, the backward orbit  $\Psi_t(x)$ for $t < 0$ is trapped.
\item For $x = (r,\theta,z)$ with $r \ne 2$, the orbit $\Psi_t(x)$ terminates in the top face $\partial_h^+ \mW$ for some $t \geq 0$, and terminates in $\partial_h^- \mW$ for some $t \leq 0$.
\item The flow $\Psi_t$  satisfies   the entry-exit condition (P2) for plugs.
\end{enumerate}
\end{prop}
\proof The only assertion that needs a comment is the last, which
follows by (W1) and the symmetry condition imposed on the functions
$g$ and $f$.
\endproof

\eject

\subsection{The family of spaces $\mK_\e$}\label{subsec-kuperberg}

The construction of the family of  Kuperberg Plugs $\mK_\e$ begins with    the modified
Wilson Plug   $\mW$ with vector field $\cW$ constructed in Section~\ref{subsec-wilson}, and follows the original
construction of K. Kuperberg, except for the choices of  self-embeddings. 
The parameter $\e$ is a real number, and admits negative and
positive values, though for $\e > 0$ we will later assume that $\e$ is ``sufficiently small''. For $\e=0$ we recover the original Kuperberg Plug.
Moreover, as proved in Sections~\ref{sec-enegative} and 
\ref{sec-entropy},   this is the only plug in the family that has no
periodic orbits.  

\begin{figure}[!htbp]
\centering
{\includegraphics[width=80mm]{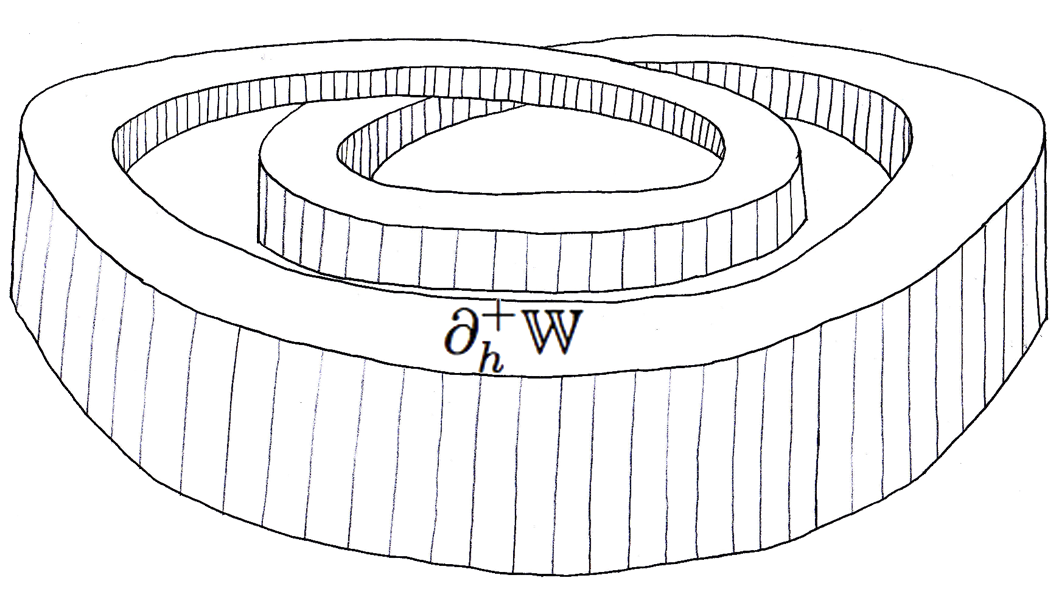}}
\caption{Embedding of Wilson Plug $\mW$ as a {\it folded figure-eight} \label{fig:8doblado} }
\vspace{-6pt}
\end{figure}

The construction follows the steps in Chapter 3 of \cite{HR2016}.   
The first step is to re-embed the manifold $\mW$   in $\mR^3$ as a {\it folded
 figure-eight}, as shown in Figure~\ref{fig:8doblado}, preserving the
vertical direction.

The fundamental   idea of the Kuperberg Plug is to   construct  two
insertions  of $\mW$ in itself. The subtlety of the original
construction arises in the precise requirements on this insertion,
which were chosen so that the periodic orbits of the flow are ``cut open'' by a non-periodic orbit for the inserted plug. 
  Here, we modify this construction so that for values of $\e \ne 0$, the self-insertion again intercepts the periodic orbits, but the resulting variation of the radius inequality results in a family of plugs with much different characteristics than Kuperberg's original construction.
Finally, in Section~\ref{sec-entropy}, we impose further restrictions on the insertion maps for  $\e > 0$, so that the entropies of the resulting Kuperberg flows can be more readily calculated.

\begin{figure}[!htbp]
\centering
{\includegraphics[width=54mm]{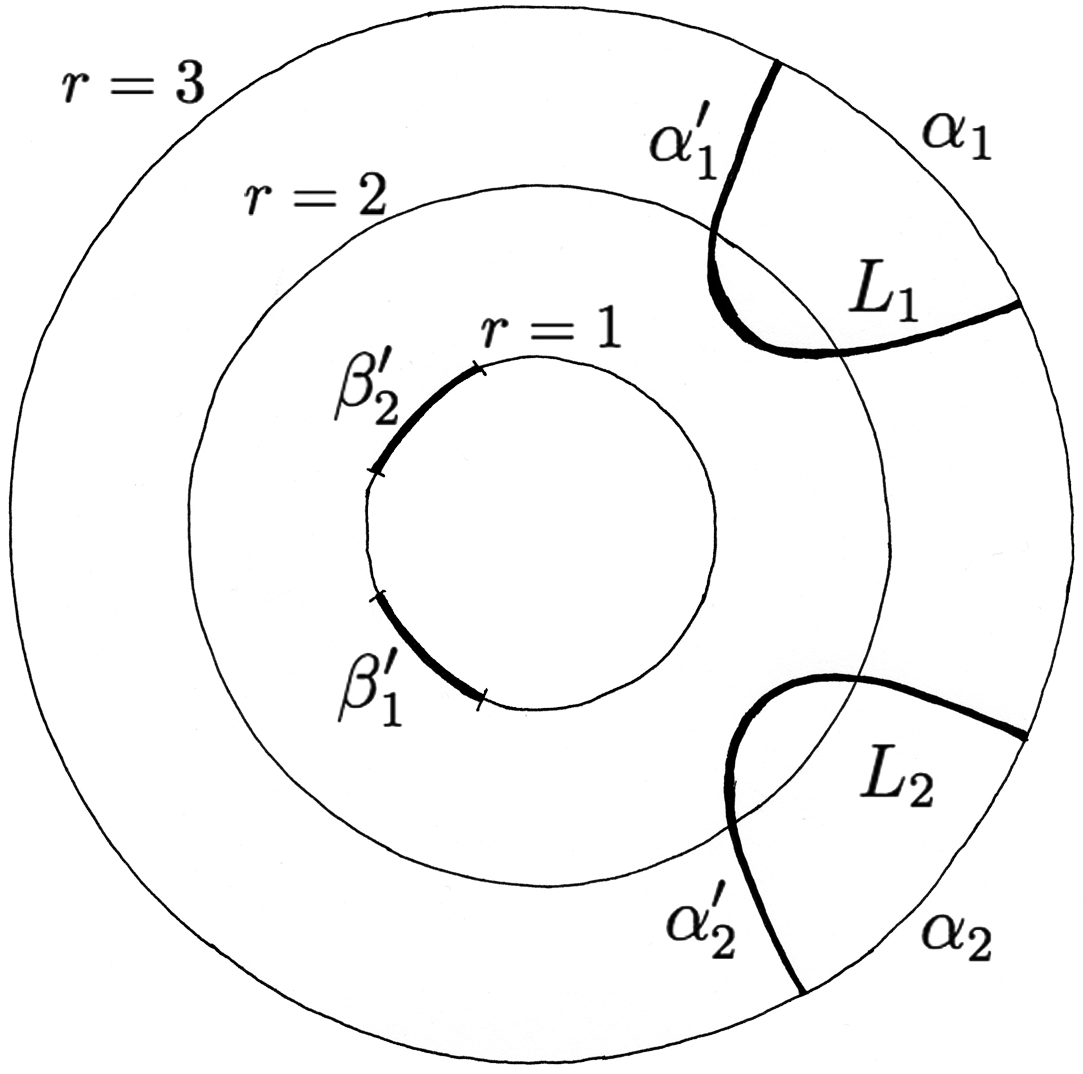}}
\caption{ The disks $L_1$ and $L_2$ \label{fig:insertiondisks}}
\vspace{-6pt}
\end{figure}

 Consider in the annulus $[1,3] \times \mS^1$    two topological closed disks $L_i$, for $i=1,2$,
whose boundaries are composed of two arcs: $\alpha^\prime_i$ in the
interior of $[1,3] \times \mS^1$,  and $\alpha_i$ in the outer boundary circle
$\{r=3\}$, as depicted in Figure~\ref{fig:insertiondisks}. To be precise, let $\zeta_1 = \pi/4$ and $\zeta_2 = - \pi/4$, then   let $\alpha_i$ be the arcs defined by 
$$
\alpha_1   ~ = ~ \{(3, \theta) \mid   ~  |\theta - \zeta_1| \leq 1/10\} \quad , \quad 
\alpha_2 ~ = ~ \{(3, \theta) \mid   ~ |\theta - \zeta_2| \leq 1/10\}
$$
We let $\alpha_i'$ be the curves which in polar coordinates $(r,\theta)$ are parabolas with minimum values $r = 3/2$ and base the line segment $\alpha_i$, as depicted in Figure~\ref{fig:insertiondisks}. We choose an explicit form for the embedded curves, for example,     given by
$\ds \alpha_i' \equiv \{ r =  3/2 + 300/2 \cdot (\theta - \zeta_i)^2  \}$.

Consider the closed sets $D_i \equiv L_i \times[-2,2] \subset \mW$, for $i = 1,2$.  
Note that  each $D_i$ is homeomorphic to a closed $3$-ball, that $D_1 \cap D_2 = \emptyset$,  and each $D_i$ intersects the  cylinder $\cC = \{r=2\}$ in   a rectangle. We label also the top and bottom faces of these regions
\begin{equation}\label{eq-regions}
L_1^{\pm}=  L_1\times \{\pm 2\}  ~ , ~ 
L_2^{\pm} =    L_2\times \{\pm 2\}.
\end{equation}

The next step is to define families of insertion maps   $\sigma_i^\e \colon D_i \to \mW $, for $i=1,2$, in such
a way that for $\e=0$ the periodic orbits $\cO_{1}$ and $\cO_{2}$ for the $\cW$-flow  intersect $\sigma_i^\e(L_i^-)$ in   points
corresponding to   $\cW$-trapped points. 
Consider the two disjoint arcs $\beta_i'$ in the inner boundary circle $\{r=1\}$,
\begin{eqnarray*}
\beta_1' & = & \{(1, \theta) \mid   ~  |\theta - (\zeta_1 + \pi) | \leq 1/10\}\\
\beta_2' & = & \{(1, \theta) \mid   ~ |\theta - (\zeta_2 + \pi) | \leq 1/10\}
\end{eqnarray*}

\begin{figure}[!htbp]
\centering
{\includegraphics[width=70mm]{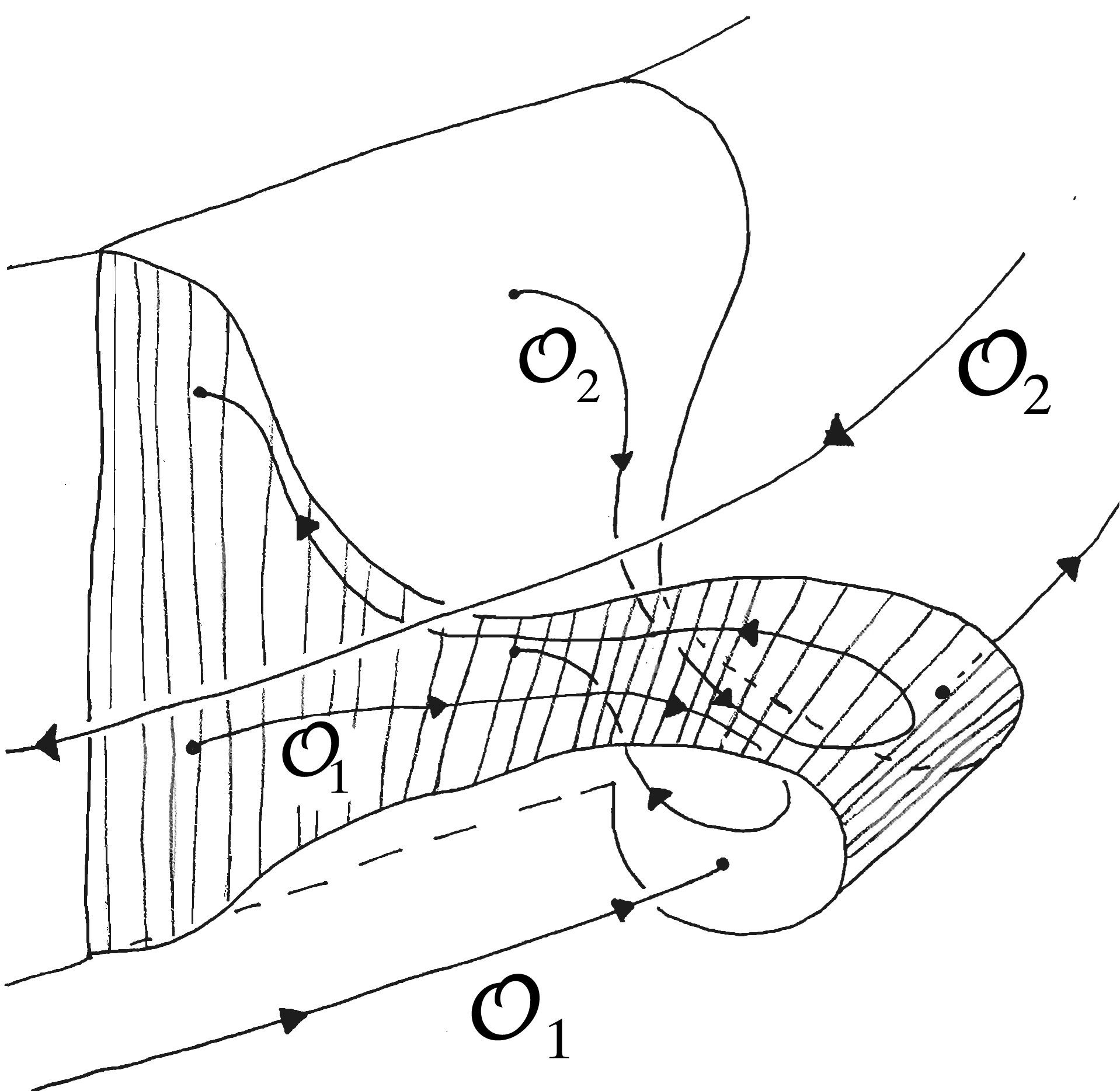}}
\caption{The image of $L_1\times [-2,2]$ in $\mW$ under $\sigma_1$ \label{fig:twisted} }
\vspace{-6pt}
\end{figure}

Now   choose a smooth family of
orientation preserving  diffeomorphisms $\sigma_i^\e \colon \alpha_i'
\to \beta_i'$, $i=1,2$, for $-a \leq \e \leq a$, where $a < \e_0$ is sufficiently small.
Extend these maps to   smooth embeddings
$\sigma_i^\e \colon D_i \to  \mW $, for $i=1,2$, as illustrated in
Figure~\ref{fig:twisted}. We require the following  conditions for all
$\e$ and for $i=1,2$:   
\begin{itemize}
\item[(K1)]   $\sigma_i^\e(\alpha_i'\times z)=\beta_i'\times z$ for all $z\in [-2,2]$,   the interior arc $\alpha_i^\prime$ is mapped to a boundary arc $\beta_i'$. 
\item[(K2)]     $\cD_i^\e = \sigma_i^\e(D_i)$ then $\cD_1^\e \cap \cD_2^\e =\emptyset$;
\item[(K3)]    For every $x \in L_i$, the image   $\cI_{i,x}^\e \equiv \sigma_i ^\e (x \times
  [-2,2])$ is an arc contained in a trajectory of $\cW$;
\item[(K4)] We have $\sigma_1 ^\e (L_1 \times \{-2\}) \subset \{z < 0\}$  and $\sigma_2 ^\e (L_2 \times \{2\}) \subset \{z > 0\}$; 
\item[(K5)]  Each slice $\sigma_i ^\e (L_i\times\{z\})$ is transverse to the vector field $\cW$, for all $-2\leq z \leq 2$. 
\item[(K6)]   $\cD_i ^\e$ intersects the periodic orbit $\cO_i$ and not $\cO_j$, for   $i \ne j$.
\end{itemize}

The ``horizontal faces'' of the embedded regions   $\cD_i^\e\subset \mW$ are labeled by
\begin{equation}\label{eq-tongues}
\cL_1^{\e\pm}= \sigma_1^\e(L_1\times \{\pm 2\})  ~ , ~   
\cL_2^{\e\pm} =   \sigma_2^\e(L_2\times \{\pm 2\}).
\end{equation}

Note that the arcs $\cI_{i,x}^\e$ are  line segments from $\sigma_i^\e(x \times \{-2\})$   to $\sigma_i^\e(x \times \{2\})$ which follow the $\cW$-trajectory, and traverse the insertion from one face to another. 
Since $\cW$ is vertical near the boundary of $\mW$ and horizontal at the two periodic orbits,     (K3-K6) imply that the arcs $\cI_{i,x}^\e$ are vertical near the inserted curve $\sigma_i ^\e (\alpha_i')$ and horizontal at the intersection of the insertion with the periodic orbit $\cO_i$. 
Thus,  the embeddings of the surfaces $\sigma_i ^\e (L_i\times\{z\})$ make a {\it half   turn} upon insertion, for each   $-2 \leq z \leq 2$,  
  as depicted in Figure~\ref{fig:twisted}.  The turning is clockwise for the bottom insertion $i=1$  as in Figure~\ref{fig:twisted}, and counter-clockwise for the upper insertion $i=2$.  
 
The first insertion $\sigma_1^\e(D_1)$ in Figure~\ref{fig:twisted}   
intersects the first periodic orbit of $\cW$ and is disjoint of the
second periodic orbit.  
  The picture of the second insertion
$\sigma_2^\e(D_2)$ is disjoint from the first
insertion and the first periodic orbit, and it intersects
the second periodic orbit.

The embeddings $\sigma_i ^\e$ are also required to satisfy two further conditions:   
\begin{itemize}
\item[(K7)] For $i=1,2$, the disk $L_i$ contains a point $(2,\theta_i)$ such that
  the image under $\sigma_i ^\e$ of the vertical segment
  $(2,\theta_i)\times[-2,2] \subset D_i \subset \mW$ is     contained
  in $\{r=2+\e\} \cap \{\theta_i^- \leq \theta \leq \theta_i^+\}$,
  and for
  $\e=0$ it is contained in $\{r=2\} \cap \{\theta_i^- \leq \theta \leq \theta_i^+\} \cap
  \{z=(-1)^i\}$.\\
\item[(K8)] {\it Parametrized Radius Inequality}: For all $x' = (r', \theta',  -2) \in  L_i^-$, let $x = (r, \theta,z) =
  \sigma_i^\e(r', \theta',  -2) \in \cL_i^{\e-}$,  then $r < r'+\e$ unless  $x' = (2,\theta_i, -2)$ and then $r=2+\e$.
\end{itemize}
 Note that by (K3) we have $r(\sigma_i^\e(r', \theta',  z')) =  r(\sigma_i^\e(r', \theta',  -2))$ for all $-2 \leq z' \leq 2$, so that (K8) holds for all points $x' = (r', \theta',  z') \in  L_i \times [-2,2]$.

\begin{figure}[!htbp]
\centering
\begin{subfigure}[c]{0.3\textwidth}{\includegraphics[height=36mm]{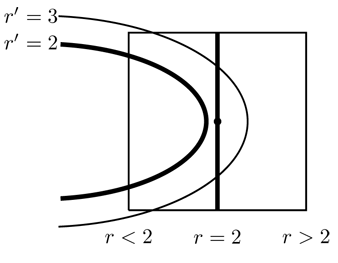}}\caption{$\e < 0$}\end{subfigure}
\begin{subfigure}[c]{0.3\textwidth}{\includegraphics[height=36mm]{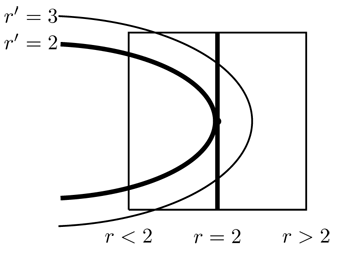}}\caption{$\e = 0$}\end{subfigure}
\begin{subfigure}[c]{0.3\textwidth}{\includegraphics[height=36mm]{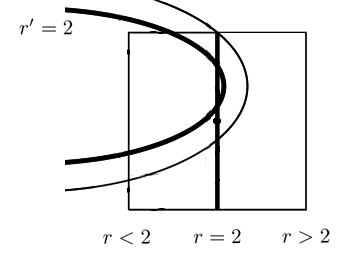}}\caption{$\e > 0$}\end{subfigure}
\caption{\label{fig:modifiedradius}  The modified radius inequality for  the cases $\e<0$, $\e=0$ and  $\e>0$}  
\vspace{-6pt}
\end{figure}

Observe that for $\e=0$, we recover the Radius Inequality of Kuperberg,
 one of the most fundamental  concepts of Kuperberg's
construction. Figure~\ref{fig:modifiedradius} represents the radius
inequality for $\e<0$, $\e =0$,  and $\e>0$. Note that in the third illustration (c) for the case $\e>0$, the insertion as illustrated has a vertical shift upwards. This is not required by conditions (K7) and (K8), but will be used to prove Theorem~\ref{thm-main2} as explained in Section~\ref{subsec-transverse}.  (The vertex points $\{ v_1^\e, v_2^\e\}$ defined in \eqref{eq-boundarypoints} correspond to this vertical offset.)

\begin{remark}\label{rmk-domains}
We add the following hypothesis relating the value of $\e_0$ in
\eqref{eq-generic1} and the insertions: let $\e_0$ as introduced in Hypothesis~\ref{hyp-genericW} be 
sufficiently small so that the $\e_0$-neighborhood of the periodic
orbits $\cO_i$ intersects the insertion regions $\cD_i^\e$ on the
interior of their faces for all $\e$ considered.
\end{remark}

Finally,  define $\mK_\e$ to be the quotient manifold obtained from $\mW$ by
identifying the sets $D_i$ with $\cD_i^\e$. That is, for each point $x
\in D_i$ identify $x$ with $\sigma_i^\e(x) \in \mW$, for $i = 1,2$. 

 The restricted $\cW$-flow on the inserted disk $\cD_i^\e = \sigma_i^\e(D_i)$ is not compatible with the restricted $\cW$-flow on $D_i$.  Thus, to obtain a smooth vector field $\cX_\e$ from this construction,  it is necessary to modify $\cW$ on each insertion $\cD_i^\e$. 
 The idea is to  replace  the vector field $\cW$  on the interior of each  region $\cD_i^\e$ with the image vector field. This requires a minor technical step first. 

Smoothly  reparametrize   the image  of $\cW|_{D_i}$ under $\sigma_i^\e$  on an open neighborhood of the boundary of $\cD_i^\e$ so that it agrees with the restriction of $\cW$ to the same neighborhood. This is possible since the vector field $\cW$ is vertical on a sufficiently small open  neighborhood of    $\partial D_i\cap \partial \mW$,
so is mapped by $\sigma_i^\e$ to an orbit segment of $\cW$ by (K3). We obtain a vector field $\cW_i'$ on $\cD_i^\e$ with the same orbits as the image of $\cW|_{D_i}$.

\begin{figure}[!htbp]
\centering
{\includegraphics[width=120mm]{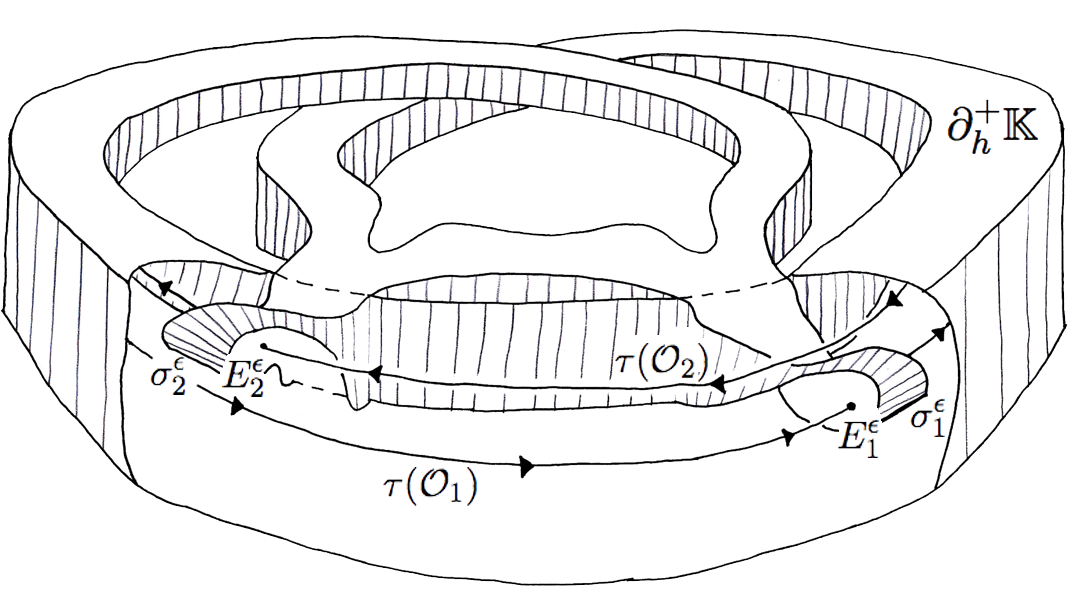}
 \caption{\label{fig:K} The Kuperberg Plug $\mK_\e$}}
\vspace{-6pt}
\end{figure}

 Then modify $\cW$ on each insertion $\cD_i^\e$, replacing it with the modified image   $\cW_i'$.   Let $\cW'$ denote the vector field on  $\mW$ after these modifications, and note that $\cW'$ is smooth.
By the modifications made above, the 
vector field $\cW^\prime$ descends to a smooth vector field on $\mK_\e$ denoted by $\cK_\e$. 
Let $\Phi_t^\e$ denote the flow of the vector field $\cK_\e$ on $\mK_\e$.   
The  resulting space  $\mK_\e \subset \mR^3$ is  illustrated  in
Figure~\ref{fig:K}. Note that  the flow  of $\cK_\e$ on $\mK_\e$
clearly satisfies   the plug conditions    (P1)  and (P4) of
Section~\ref{sec-plugs}, while  Proposition~\ref{prop-ee} below will
show that the condition (P2) is also satisfied.  For the cases $\e < 0$, the trapped
  orbit condition (P3) is shown by  Corollary~\ref{cor-trapped},
  so that the flow of $\cK_\e$ on $\mK_{\e}$ is a plug in the sense of
  Section~\ref{sec-plugs}.

\section{Level and radius functions}\label{sec-radiuslevel}

In this section we recover some of the results on the orbit behaviors of the flows $\mK_\e$ 
  that are independent of $\e$ and that will be useful
for the study of their dynamics in the plugs $\mK_\e$. 
As with the study of the Kuperberg flow $\cK_0$, an orbit in $\mK_\e$ is formed by concatenating
pieces of orbits of the   Wilson plug, where these orbit segments are the result of the construction of $\mK_\e$ via the insertion maps $\sigma_i^\e$.
Thus,  understanding  the way these arcs concatenate is a fundamental aspect of understanding the dynamics of of the flow in $\mK_\e$.

We start by introducing notations that will be used throughout
this work, and also some basic concepts which are fundamental for
relating the dynamics of the two vector fields $\cW$ and $\cK_\e$. These
results are contained in the literature \cite{Ghys1995,HR2016, Kuperberg1994, Kuperbergs1996, Matsumoto1995}, though in a variety
of differing notations and presentations. We adopt the notation of \cite{HR2016}, and parts of the following text is also adapted from this work.

Recall that   $\cD_i^\e = \sigma_i^\e(D_i)$ for $i =1,2$   are   solid $3$-disks embedded in $\mW$.
 Introduce  the sets:
\begin{equation}\label{eq-notchedW}
\mW'_\e ~   \equiv   ~ \mW - \left \{ \cD_1^\e  \cup \cD_2^\e \right\} \quad , \quad 
\wmW~   \equiv   ~ \overline{\mW - \left \{ \cD_1^\e  \cup \cD_2^\e \right\}} ~ .
\end{equation}
The  closure $\wmW$ of $\mW'_\e$  is the  \emph{pi\`{e}ge de Wilson  creus\'{e}}  as defined in \cite[page 292]{Ghys1995}.
The compact space $\wmW \subset \mW$ is the result of  ``drilling out''  the interiors of $\cD_1^\e$ and $ \cD_2^\e$, as the   terminology \emph{creus\'{e}} suggests.

For $x, y \in \mK_\e$, we say that $x \prec_{\cK_\e} y$ if   there exists $t \geq 0$ such that $\Phi_t^\e(x) = y$.
Likewise, for $x',y' \in \mW$, we say that  $x' \prec_{\cW} y'$  if   there exists $t \geq 0$ such that $\Psi_t(x') = y'$.

  Let   $\tau \colon \mW  \to \mK_\e$ denote the quotient map,
   which for $i = 1,2$,    identifies   a point $x \in D_i$ with its
   image $\sigma_i^\e(x)   \in \cD_i^\e$. Even if $\tau$ depends on
   $\e$, we will denote it simply by $\tau$.
Then   the  restriction $\tau' \colon \mW'_\e \to \mK_\e$ is injective and onto. Let   $(\tau')^{-1}  \colon \mK_\e \to \mW'_\e$ denote the inverse map, which followed  by the inclusion $\mW'_\e \subset \mW$, yields  the (discontinuous) map  $\tau^{-1} \colon \mK_\e \to \mW$, where  $i=1,2$, we have:
\begin{equation}\label{eq-radiusdef}
\tau^{-1}(\tau(x)) =x ~ {\rm for} ~ x \in D_i  ~ ,~ {\rm and} ~   \sigma_i^\e(\tau^{-1}(\tau(x))) = x ~ {\rm for} ~ x \in \cD_i^\e~.
\end{equation}
For $x \in \mK_\e$,  let $x=(r,\theta,z)$  be  defined as  the $\mW$-coordinates   of  $\tau^{-1}(x) \subset \mW'_\e$. 
In this way, we obtain   (discontinuous)     coordinates $(r,\theta,z)$ on $\mK_\e$. 

In particular, let $r \colon \mW'_\e \to [1,3]$ be 
 the restriction of the   radius coordinate on $\mW$, then the function is extended to   the \emph{radius function}  of $\mK_\e$, again denoted by $r$, where for $x \in \mK_\e$ set $r(x) = r(\tau^{-1}(x))$.  
    The flow of the vector field $\cW$ on $\mW$ preserves the radius function on $\mW$, so   $x' \prec_{\cW} y'$ implies that      $r(x') = r(y')$.  
However, $x \prec_{\cK_\e} y$ need not imply that $r(x) = r(y)$, and the 
   points of discontinuity are the transition points defined below.

Let $\partial_h^-\mK_\e \equiv  \tau(\partial_h^- \mW\setminus (L_1^-\cup
L_2^-))$ and $\partial_h^+\mK_\e \equiv \tau(\partial_h^+ \mW\setminus
(L_1^+\cup L_2^+))$  denote  the   bottom and top horizontal faces of
$\mK_\e$, respectively.   Note that the vertical boundary component  $\partial_v\mK_\e \equiv \tau(\wmW \cap \partial_v \mW)$   is  tangent to the flow.

Points $x' \in \partial_h^- \mW$ and $y' \in \partial_h^+ \mW$ are said to be \emph{facing}, and we  write $x' \equiv y'$,  
if      $x' = (r, \theta, -2)$ and $y' = (r, \theta, 2)$ for some $r$ and $\theta$.  
There is also a notion of facing points for $x, y \in \mK_\e$,   if either of   two cases are satisfied: 
\begin{itemize}
\item For  $x = \tau(x') \in \partial_h^-\mK_\e$ and $y = \tau(y') \in \partial_h^+\mK_\e$, if  $x' \equiv y'$ then $x \equiv y$.
\item For    $i=1,2$, for $x', y' \in \partial^{\pm} \mW$ and $x = \sigma_i^\e(x')$ and $y = \sigma_i^\e(y')$,     if  $x' \equiv y'$ then $x \equiv y$.
\end{itemize} 
 The context in which the notation $x \equiv y$ is used dictates which usage applies.

Consider the embedded disks  $\cL_i^{\e\pm} \subset \mW$ defined by \eqref{eq-tongues}, which appear  as the   faces of the
insertions  in $\mW$.  Their images in the quotient manifold $\mK_\e$ are denoted by:
\begin{equation}\label{eq-sections}
E^{\e}_1= \tau(\cL_1^{\e-}) ~ , ~ S^{\e}_1= \tau(\cL_1^{\e+}) ~ , ~   
E^{\e}_2= \tau(\cL_2^{\e-}) ~ , ~ S^{\e}_2= \tau(\cL_2^{\e+})   ~ .
\end{equation}
Note that $\tau^{-1}(E^{\e}_i) = L_i^-$, while $\tau^{-1}(S^{\e}_i) = L_i^+$.  The
{\it transition   points}  of an orbit of $\cK_\e$ are those points which meet
the $E^{\e}_i$, $S^{\e}_i$, $\partial_h^-\mK_\e$ or $\partial_h^+\mK_\e$, for $i=1,2$.
 They    are then either \emph{primary} or \emph{secondary} transition
 points, where $x\in \mK_\e$ is: 
\begin{itemize}
\item a \emph{primary entry  point} if $x\in \partial_h^-\mK_\e$;
\item  a \emph{primary exit  point} if $x \in \partial_h^+\mK_\e$;
\item a \emph{secondary entry  point} if $x \in E^{\e}_1 \cup E^{\e}_2$;
\item a \emph{secondary exit  point}  $x \in S^{\e}_1 \cup S^{\e}_2$.
\end{itemize}
If a $\cK_\e$-orbit contains no transition points, then it lifts to a $\cW$-orbit in $\mW$ which flows   from    $\partial_h^- \mW$ to     $\partial_h^+ \mW$.

\begin{figure}[!htbp]
\centering
{\includegraphics[width=80mm]{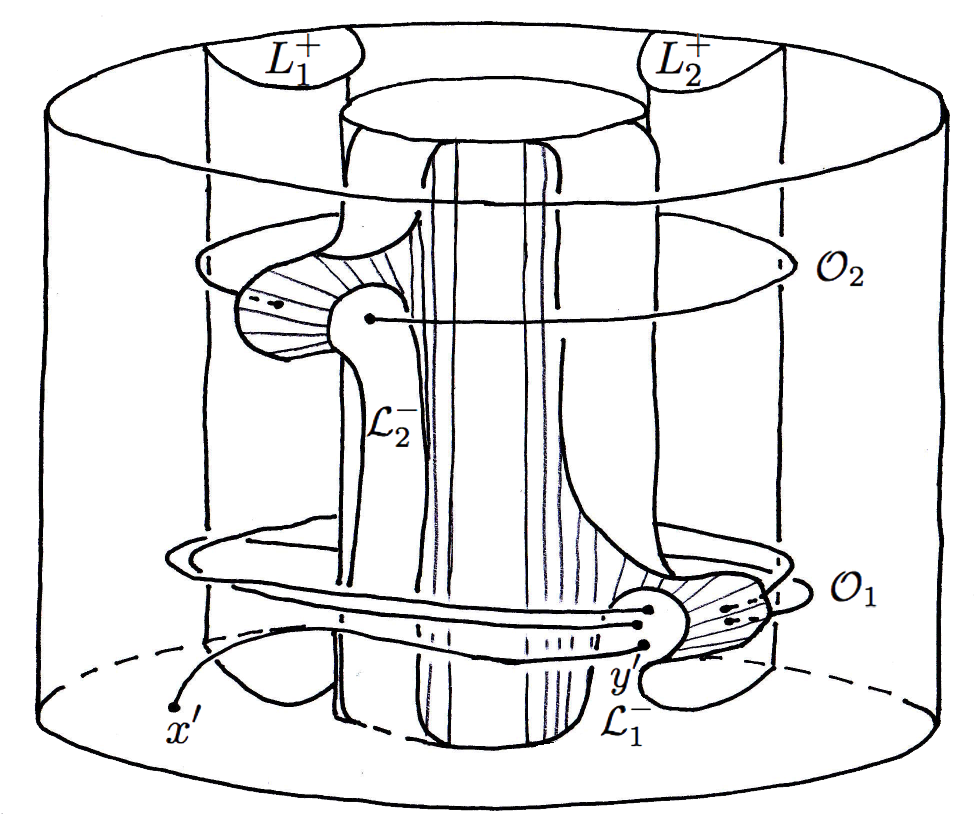}}
\caption{$\cW$-arcs lifted to $\mW'_\e$ \label{fig:cWarcs}}
\vspace{-6pt}
\end{figure}

 A \emph{$\cW$-arc} is a closed segment $[x,y]_{\cK_\e} \subset \mK_\e$ of
 the flow of $\cK_\e$ whose endpoints $\{x,y\}$ are the only transition
 points in $[x,y]_{\cK_\e}$. The open interval $(x,y)_{\cK_\e}$  is then the image under $\tau$ of a unique $\cW$-orbit segment in $\mW'_\e$, denoted by $(x',y')_{\cW}$ (see Figure~\ref{fig:cWarcs}.) 
 Let $[x', y']_{\cW}$ denote the closure of $(x',y')_{\cW}$   in $\wmW$, then  we say that $[x', y']_{\cW}$ is the \emph{lift} of $[x,y]_{\cK_\e}$.  Note that   the radius function $r$ is constant along   $[x', y']_{\cW}$.

 The properties of the Wilson flow $\cW$ on $\wmW$ determine   the
 endpoints of lifts   $[x',y']_{\cW}$.  We state the six cases which
 arise explicitly, as they will be cited in  later arguments.  For a
 proof we refer to Lemma~4.1 of \cite{HR2016}.  Figure~\ref{fig:cWarcs} helps in visualizing these cases.  
 \begin{lemma}\label{lem-cases1}
Let  $[x,y]_{\cK_\e} \subset \mK_\e$ be a $\cW$-arc, and let $[x',y']_{\cW} \subset \wmW$ denote its lift.
 \begin{enumerate}
 
\item {\bf (p-entry/entry)}   If $x$ is a primary entry point, 
then $\ds x' \in    \partial_h^-  \mW \setminus (L_1^-   \cup L_2^- )$, 
and  if $y$ an entry point, then $y$ is   a secondary   entry point and  $y' \in \cL_i^{\e-}$ for  $i =1$ or $2$.  
  
\item   {\bf (p-entry/exit)} If $x$ is a primary entry point, then $x' \in   \partial_h^-
  \mW \setminus ( L_1^-   \cup L_2^- )$, and  if  $y$ is an exit point, 
  then   $\ds y' \in \partial_h^+\mW\setminus (L_1^+\cup L_2^+)$ is a primary exit point, and by  the entry/exit condition on $\mW$ we have   $x\equiv y$.     

\item  {\bf   (s-entry/entry)} If $x$ is a secondary entry point, then $x' \in    L_i^-$,  for  $i =1$ or $2$, and  if $y$ is   an entry point,  then we have   $y'  \in    \cL_j^{\e-} $ where $j =1,2$   is not necessarily    equal to $i$.

\item   {\bf  (s-entry/exit)} If $x$ is a secondary entry point, then $x' \in    L_i^-$,  for  $i =1$ or $2$, and  if  $y$ is an   exit point, then $y$ is a secondary exit point,  $y'
  \in    L_i^+$ and $x \equiv y$ by the entry/exit condition  of $\mW$.  

\item  {\bf  (s-exit/entry)} If $x$ is a secondary exit point, then $x' \in   \cL_i^{\e+}$,  for  $i =1$ or $2$, and  if $y$ is   an entry point, so that $y' \in   \cL_j^{\e-}$ then $j=2$ if  $i=2$, and $j =1,2$  if $i=1$.

\item {\bf  (s-exit/exit)} If $x$ is a secondary exit point, then $x' \in   \cL_i^{\e+}$,  for  $i =1$ or $2$, and  if   $y$ is a  primary
  exit point, $y' \in  \{ \partial_h^+ \mW \setminus ( L_1^+ \cup L_2^+) \}$. If $y$ is a secondary exit point, then $y'  \in    L_j^+$, where  $j =1$ or $2$ is not necessarily    equal to $i$. 
\end{enumerate}
\end{lemma}

 Figure~\ref{fig:cWarcs} illustrates some of the notions discussed in this section. 
 The disks $L_1^-$ and $L_2^-$ contained in  $\partial_h^-\mW$ are drawn in the bottom face, though they are partially obscured by the   inner cylindrical boundary   $\{r=1\}$.
The  image of $L_1^-$ under $\sigma_1^\e$ is  the entry face of the insertion region in the lower half of the cylinder, 
while the  image of $L_2^-$ under $\sigma_2^\e$ is  the entry face of the insertion region in the upper half of the core cylinder. 
Analogously, the   disks $L_1^+$ and $L_2^+$ in $\partial_h^+\mW$  are mapped to the exit region of the insertion regions. 
The intersection of the $\cW$-periodic orbits $\cO_i$ with the
insertions in $\mW$ are illustrated,
as well as two  $\cW$-arcs in $\mW'_\e$ that belong to the same orbit. One $\cW$-arc  goes from
$\partial_h^-\mW$ to $\cL_1^{\e-}$, hence from a principal entry point to
a secondary entry point (as in Lemma \ref{lem-cases1}(1)).  The second $\cW$-arc   goes from $\cL_1^{\e+}$ to
$\cL_1^{\e-}$, thus from a secondary exit point to a secondary entry point (as in Lemma \ref{lem-cases1}(5)).

    Introduce the radius coordinate function along $\cK_\e$-orbits, where 
for $x \in \mK_\e$,  set     $\rho_x^\e(t) \equiv r(\Phi_t^\e(x))$.
Note that if $\Phi_{t}^\e(x)$ is not a transition point then the    function $\rho_x^\e(t)$  is locally constant at $t$, and thus if a $\cK_\e$-arc $\{\Phi_t^\e(x) \mid t_0 \leq t \leq t_{1}\}$  contains no transition point, then $\rho_x^\e(t)=\rho_x^\e(t_0)$ for all $t_0 \leq  t \leq  t_1$.

The  \emph{level function}  along an orbit  indexes the discontinuities of the  radius function. 
Given $x \in \mK_\e$, set      $n_x(0) = 0$,   and for $t > 0$, define 
\begin{equation}\label{def-level+}
n_x(t) = \#   \left\{ (E^{\e}_1 \cup E^{\e}_2) \cap \Phi_s^\e(x) \mid  0 < s \leq t \right\} - \#   \left\{ (S^{\e}_1 \cup S^{\e}_2) \cap \Phi_s^\e(x) \mid  0 < s \leq t  \right\} .
\end{equation}
 That is,  $n_x(t)$ is the total number of secondary entry points, minus
 the total number of secondary exit points,  traversed by the flow of $x$ over the interval 
 $0 < s \leq  t$.  

The function can be extended to negative time by setting, for $t < 0$,
 \begin{equation}\label{def-level-}
n_x(t) =  \#   \left\{ (S^{\e}_1 \cup S^{\e}_2) \cap \Phi_s^\e(x) \mid  t < s \leq 0  \right\}    -   \#   \left\{ (E^{\e}_1 \cup E^{\e}_2) \cap \Phi_s^\e(x) \mid  t < s \leq 0 \right\}.
\end{equation}

We use throughout this work a Riemannian metric on   the tangent bundle to $\mK_\e$. The Wilson plug $\mW$ has a natural product Riemannian metric, where the rectangle $\bR$ in \eqref{eq-rectangle} has the product euclidean metric, and the circle factor $\mS^1$ has length $2\pi$. Then modify this metric along the insertions $\sigma_i^\e$ to smooth them out, and so obtain a Riemannian metric on $\mK_\e$.

 For $x' \prec_{\cW} y'$ in $\mW$,  let   $d_{\mW}(x',y')$ denote the
 path length of the $\cW$-orbit segment $[x' ,y']_{\cW}$ between
 them. Similarly, for $x \prec_{\cK_\e} y$ in $\mK_\e$,  let
 $d_{\mK_\e}(x,y)$ denote the path length of the $\cK_\e$-orbit segment $[x,y]_{\cK_\e}$. Note that if $[x,y]_{\cK_\e}$ is a $\cW$-arc with lift $[x',y']_{\cW}$ then we have $d_{\mK_\e}(x,y) = d_{\mW}(x',y')$ by the choice of the metric on $\mW$.

In Chapter 4 of \cite{HR2016} (Lemmas~4.3 and 4.4 and
Corollary~4.5) the following basic length
estimates are established for the Kuperberg Plug $\mK_0$, and the proofs for this case carry over directly to the flows on the plugs $\mK_\e$:
\begin{enumerate}
\item Let $0 < \delta < 1$. There exists $L(\delta) > 0$ such that for any $\xi \in \mW$ with $|r(\xi) -2| \geq \delta$,    the total $\cW$-orbit segment  $[x',y']_{\cW}$ through $\xi$ has   length bounded above by $L(\delta)$.
\item There exists $0 < d_{min} <
  d_{max}$ such that  if  $[x',y']_{\cW} \subset \wmW$ is the lift of
  a $\cW$-arc $[x,y]_{\cK_\e}$, then we have the uniform estimate
\begin{equation}\label{eq-segmentlengths}
d_{min} \leq  d_{\mW}(x',y') \leq d_{max} ~ .
\end{equation}
\end{enumerate}
Thus for  $[x,y]_{\cK_\e} \subset \mK_\e$  a   \emph{$\cW$-arc}, there is a uniform length estimate
 \begin{equation}\label{eq-segmentlengths2}
d_{min} \leq  d_{\mK_\e}(x,y) \leq d_{max} ~ . 
\end{equation}
 
We end this section by recalling   some technical results, which are  key to analyzing the orbits of the flows $\Phi_t^\e$  on $\mK_\e$. For a proof
of the first of these, we refer to     \cite[Proposition 5.5]{HR2016}, whose proof carries over directly to the case of the flow $\Phi_t^{\e}$ on $\mK_\e$.

\begin{prop}  \label{prop-shortcut4}
Let $x   \in \mK_\e$.  For $n \geq 3$, assume that we  are given
successive transition points   $x_{\ell} = \Phi_{t_{\ell}}^{\e}(x)$
for $0 = t_0 < t_1 < \cdots < t_{n-1} < t_n$. Let
$[x_\ell',y_{\ell+1}']_\mW$ be te lift of the $\cW$-arc
$[x_\ell,x_{\ell+1}]_\mK$ for all $0\leq \ell \leq n-1$.
  Suppose that $n_{x_0}(t) \geq 0$ for all $0 \leq  t < t_n$, and that $n_{x_0}(t_{n-1}) = 0$. 
Then    $x_0'  \prec_{\cW} y_n'$, and hence $r(x_0') = r(y_n')$.
Moreover, if $x_0$ is an    entry point and $x_n$ is an exit point, then 
  $x_0\equiv x_n$.
\end{prop}

For a proof of the next result, we refer to the proof of \cite[Proposition 5.6]{HR2016}.
\begin{prop}\label{prop-ee} 
Let   $x \in \partial_h^-\mK_\e$   be a primary entry   point, and label the successive transition points by   $x_{\ell} = \Phi_{t_{\ell}}^\e(x)$ for $0 = t_0 < t_1 < \cdots < t_{n-1} < t_n$. 
   If $x_n$ is a primary exit point,  then $x \prec_{\cW} x_n$ and hence $x
\equiv x_n$.  Moreover, $n_x(t) \geq 0$ for $0 \leq t < t_{n}$.
\end{prop}

Note that this implies the flow  of $\cK_\e$ on $\mK_\e$  satisfies
condition  (P2)  of Section~\ref{sec-plugs}.  The following result
proves that $\mK_\e$ satisifies condition (P3) for all $\e$, implying
that $\mK_\e$ is a plug.
 
\begin{cor}\label{cor-trapped}
Let $x\in \partial_h^-\mK_\e$ be a primary entry point with
$r(x)=2$. Then the forward $\cK_\e$-orbit of $x$ is trapped.
\end{cor}

\proof
Assume that the forward $\cK_\e$-orbit of $x$ is not trapped, and let
$x_{\ell} = \Phi_{t_{\ell}}^\e(x)$ for $0 = t_0 < t_1 < \cdots <
t_{n-1} < t_n$ be the transition points. Then $x_n$ is a primary exit
point and by Proposition~\ref{prop-ee}, $x\prec_\cW x_n$, which is a
impossible since $r(x)=2$.
\endproof

\section{Propellers and double propellers in $\mW$}\label{sec-propeller}

Propellers, simple and double, as defined in Chapters 11, 12 and 13 of
\cite{HR2016}, are obtained  from the Wilson flow of  selected arcs
in $\partial_h^-\mW$. Double propellers are formed by the union of two
(simple) propellers, as described below. These geometric structures associated to the Kuperberg flows $\Phi_{t}^{\e}$ are an essential tool for analyzing their dynamical properties. 

Consider a   curve in the entry region, $\g \subset \partial_h^- \mW$, with a parametrization  $\g(s) = (r(s), \theta(s), -2)$ for $0 \leq s \leq 1$, 
 and assume that the map $\g \colon [0,1] \to \partial_h^- \mW$ is a homeomorphism onto its image. We use the notation $\g_s = \g(s)$ when convenient, so that 
  $\g_0 = \g(0)$ denotes the initial point and $\g_1 = \g(1)$  denotes  the terminal point  of $\g$. 
For    $\delta > 0$,  assume    that $r(\g_0) = 3$,   $r(\g_1)  = 2+\delta$, and $2 +  \delta < r(\g_s) < 3$ for $0 < s < 1$.

The $\cW$-orbits of the
points in $\g$ traverse $\mW$ from $ \partial_h^- \mW$ to $ \partial_h^+ \mW$, and hence the flow of $\g$ generates a compact invariant surface $P_{\g} \subset \mW$.
The surface $P_{\g}$ is  parametrized by $(s,t) \mapsto \Psi_t(\g(s))$ 
for $0 \leq s \leq 1$ and $0 \leq t \leq T_s$, where $T_s$ is the exit time for the $\cW$-flow of $\g(s)$. Observe that as  $s \to 1$ and $\delta \to 0$, the exit time $T_s    \to \infty$.

\begin{figure}[!htbp]
\centering
\begin{subfigure}[c]{0.4\textwidth}{\includegraphics[height=54mm]{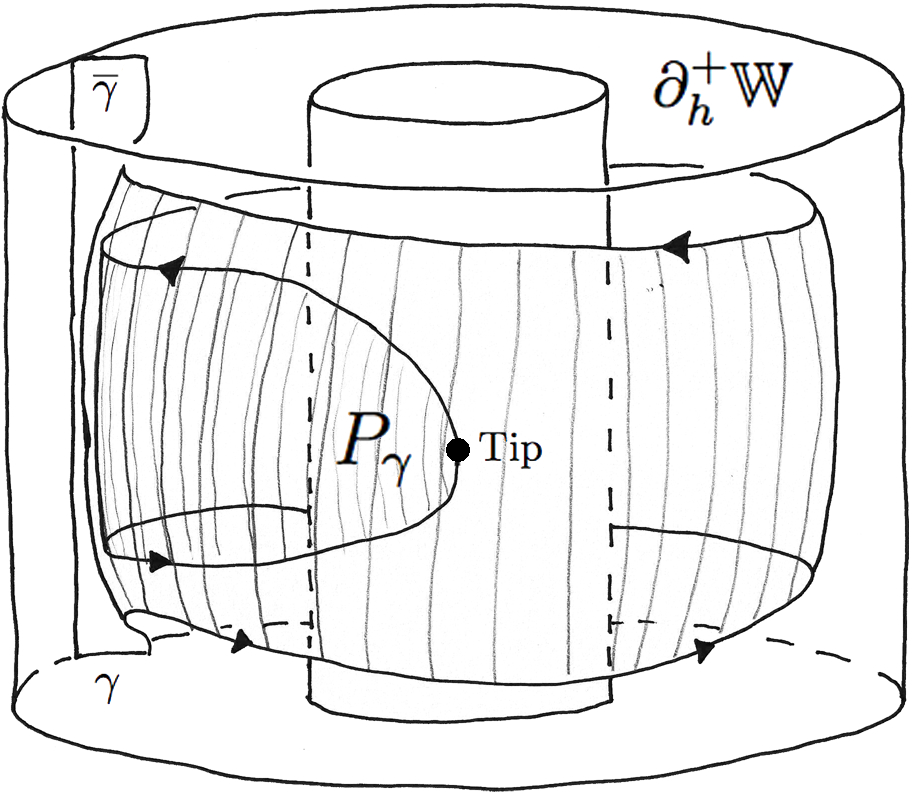}}\end{subfigure}
\begin{subfigure}[c]{0.5\textwidth}{\includegraphics[height=48mm]{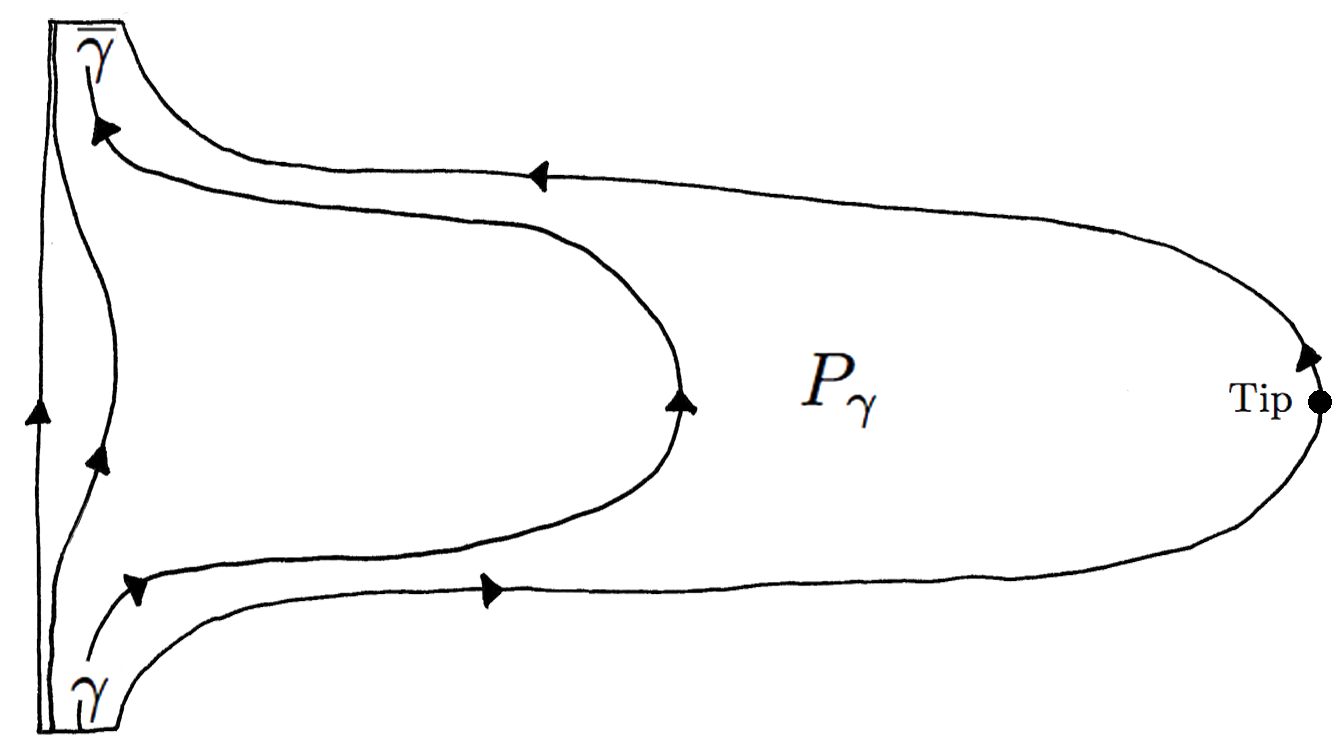}}\end{subfigure}
\caption{\label{fig:propeller}   Embedded and flattened finite propeller}
\vspace{-6pt}
\end{figure}

The surface  $P_{\g}$ is called a {\it propeller}, due to the nature of its shape in $\mR^3$. It takes the form of a ``tongue'' wrapping around the core cylinder $\cC(2 + \delta)$ which contains the orbit of $\g_1$. To visualize the shape of this surface, consider the   case where $\g$
 is topologically transverse to the cylinders $\cC(r_0) = \{r=r_0\}$ for $2+ \delta \leq r_0 \leq 3$. 
The transversality assumption implies that the radius  $r(\g_s)$ is \emph{monotone decreasing} as $s$ increases.
   Figure~\ref{fig:propeller}  illustrates the surface $P_{\g}$ as  a ``flattened'' propeller on the right, and   its embedding in $\mW$ on the left. As $\delta \to 0$ the surface  approaches  the   cylinder $\cC= \{r=2\}$ in an infinite spiraling manner. 

We comment on the details in Figure~\ref{fig:propeller}.
The horizontal boundary   $\partial_h  P_{\g}$ is composed of the
initial curve   $\g \subset \partial_h^- \mW$, and its      mirror
image   $\ovg \subset \partial_h^+ \mW$ via the entry/exit condition
on the Wilson Plug. The vertical boundary   $\partial_v P_{\g}$ is
composed of the vertical segment $\g_0 \times [-2,2]$ in $\partial_v \mW$, and the
orbit  $\{\Psi_t(\g_1) \mid 0 \leq t \leq T_1\}$ which is the inner (or long) edge in the interior of $\mW$.
One way to visualize the surface,   is to consider the product surface $\g \times [-2,2]$, and  then start deforming it by an isotopy which follows the flow lines of $\cW$, as illustrated in Figure~\ref{fig:propeller}.
In the right hand side of the figure,  some of the orbits in the propeller
are illustrated, while in the left hand side,  just the boundary orbit is
illustrated. 

 Consider the orbit $\{\Psi_t(\g_1) \mid 0 \leq t \leq T_1\}$  of the endpoint $\g_1$ with $r(\g_1) = 2+\delta$.
  The path $t \mapsto \Psi_t(\g_1)$    makes a certain
number of turns in the positive $\mS^1$-direction before reaching
the core annulus $\cA$ at $z = 0$. The Wilson vector field $\cW$ is vertical on the plane $\cA$,     so the flow of $\g_1$ then crosses   $\cA$, after which the orbit $\Psi_t(\g_1)$ starts turning in the negative direction and ascending until it reaches $\partial_h^+ \mW$. 
The point where the flow $\Psi_t(\g_1)$ intersects $\cA$ is called the \emph{tip} of the propeller $P_{\g}$.

The anti-symmetry of the vector field $\cW$ implies that the number of turns in one direction (considered as a
real number)   equals   the number of turns in the opposite direction. To be precise, for $\g_1 = (r_1, \theta_1, -2)$ in coordinates,  let $\Psi_t(\g_1) = (r_1(t), \theta_1(t), z_1(t))$ in coordinates. The function $z_1(t)$ is monotone increasing, and by the symmetry, we have $z_1(T_1 /2) = 0$. Thus, the tip is the point 
$\Psi_{T_1 /2}(\g_1)$.

 Next,  for fixed $0 \leq a < 2\pi$, consider the intersection of $P_{\g}$ with a slice 
 \begin{equation}\label{eq-slices}
\bRa \equiv \{\xi = (r, a, z) \mid ~ 1 \leq r \leq 3 ~,   ~ -2 \leq z \leq 2\} ~  .
\end{equation}
Each rectangle $\bRa$ is tangent to the Wilson flow along the annulus $\cA$, and also near the boundaries of $\mW$, but is transverse to the flow at all other points. The case when  $a = \theta_1(T_1/2)$ is special,   as  the tip of the propeller is tangent to $\bRa$.

\begin{figure}[!htbp]
\centering
\begin{subfigure}[c]{0.4\textwidth}{\includegraphics[width=50mm]{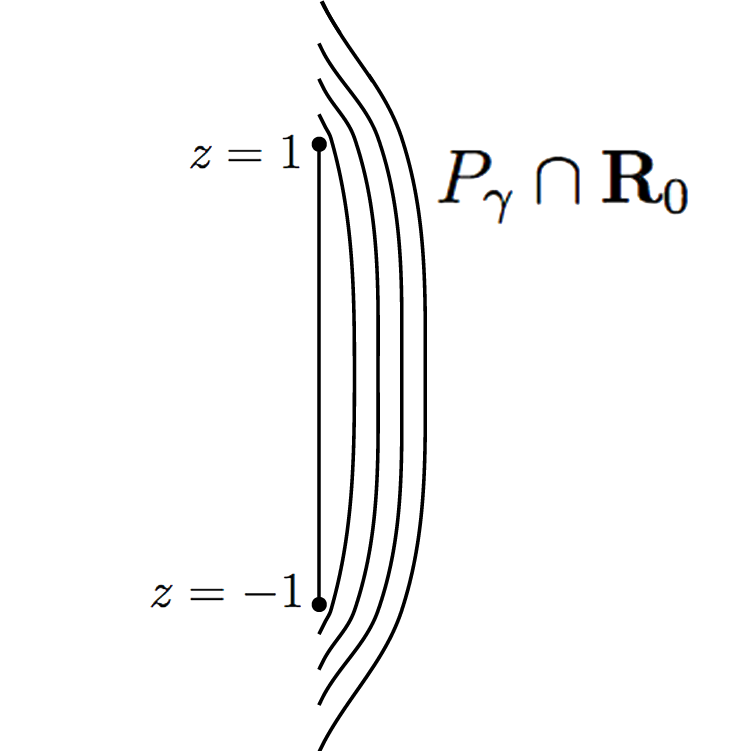}}\caption{Infinite propeller $P_{\g}$}\end{subfigure}
\begin{subfigure}[c]{0.4\textwidth}{\includegraphics[width=50mm]{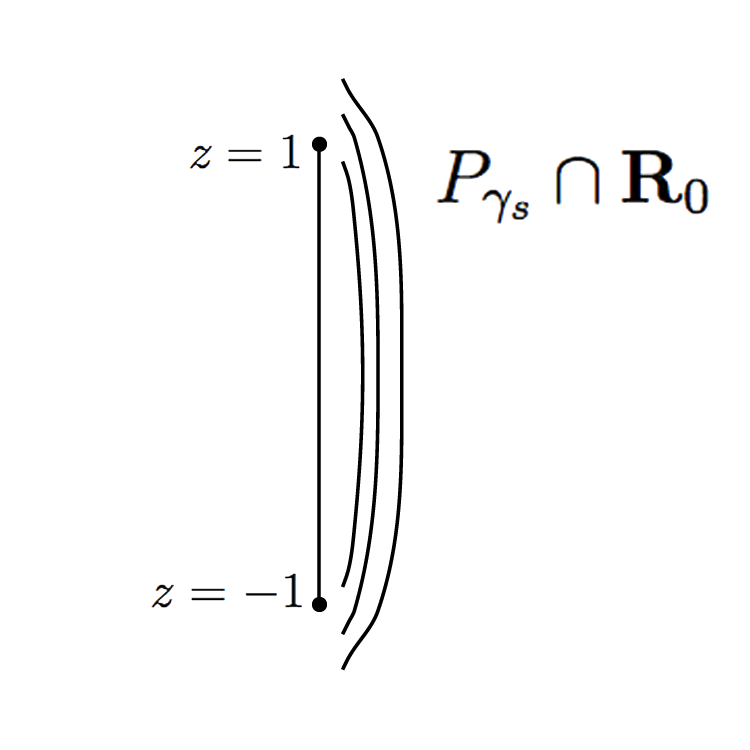}}\caption{Finite propeller $P_{\g}$}\end{subfigure}
\caption{\label{fig:arcspropeller} Trace of   propellers   in $\bRt$}
\vspace{-6pt}
 \end{figure}

Assume that $a \ne  \theta_1(T_1/2)$, then the orbit 
$\Psi_t(\g_1)$ intersects $\bRa$ in a series of points
on the line $\cC(2+\delta) \cap \bRa$ that are paired, as illustrated
in the right hand side of Figure~\ref{fig:arcspropeller}. 
Moreover, the intersection $P_{\g} \cap \bRa$ consists of a finite sequence of arcs between the symmetrically paired points of $\Psi_t(\g_1) \cap \bRa$.
 The number of such arcs is equal to, plus or minus one,  the number of times the curve
 $\Psi_t(\g_1)$ makes a complete turn around the cylinder $\cC(2+\delta)$. 
  
 We comment on the details of Figure~\ref{fig:arcspropeller}, which illustrates the case of $\bRt$ where $a=\theta_0$.
 The vertical line between the points $(2,\theta_0,-1)$ and $(2,\theta_0,1)$ (marked
 in the figure simply by $z=-1$ and $z=1$, respectively) is the trace
of the Reeb cylinder in $\bRt$. The trace of a propeller in $\bRt$ is a collection of arcs that
have their endpoints in the vertical line $\{r=2 + \delta\}$. In the left
hand figure, $r(\g_1)=2$ and the propeller in consideration is
infinite, as defined below. The curves form an infinite family, here just four arcs are shown,
accumulating on the vertical line. The right hand figure illustrates the case  $r(\g_1)>2$,  and the propeller is
finite.

Finally,  consider the case where $\delta \to 0$, so that    the endpoint $\g_1$ of   $\g$
lies in  the cylinder $\cC$.  Then for $0 \leq s < 1$,  we have $r(\g(s)) > 2$, so the $\Psi_t$-flow of 
$\g_s \in    \partial_h^- \mW$   escapes from $\mW$. Define the curve $\ovg$ in $\partial_h^+ \mW$ to be the trace of these facing endpoints in $\partial_h^+ \mW$, parametrized by $\ovg(s)$ for $0 \leq s < 1$, where $\g(s) \equiv \ovg(s)$. 
Define $\ds \ovg_1 = \lim_{s\to 1} ~ \ovg(s)$ so that $\ovg_1 \equiv \g_1$ also. 

Note that  the forward $\Psi_t$-orbit of $\g_1$ is   asymptotic to the periodic orbit
$\cO_1$,  while the backward   $\Psi_t$-orbit of $\ovg_1$ is   asymptotic to the periodic orbit
$\cO_2$. Introduce their ``pseudo-orbit'', 
\begin{equation}\label{eq-Zset}
\cZ_{\g} = \cZ_{\g}^- ~ \cup ~ \cZ_{\g}^+ ~ , ~ \cZ_{\g}^- =
\{\Psi_t(\g_1) \mid t \geq 0 \}  ~ \text{and} ~ \cZ_{\g}^+ =  \{\Psi_t(\ovg_1) \mid t \leq 0 \}
\end{equation}
Each curve $\cZ_{\g}^{\pm}$ traces out a semi-infinite ray in $\cC$ which spirals from the bottom or top face to a periodic orbit, and thus 
$\cZ_{\g}$ traces out two semi-infinite curves in $\cC$ spiraling to the periodic orbits $\cO_1 \cup \cO_2$. 
 
For $0 < \delta \leq 1$,  denote by $\g^{\delta}$ the curve with image $\g([0, \delta])$, parametrized by 
\begin{equation}\label{eq-sparem}
\g^{\delta}(s) = \g(\delta \cdot s) . 
\end{equation}

\begin{defn}\label{def-infpropeller}
Let $\g$ be a curve parametrized by $\g \colon [0,1]\to \mW$   as above, with $r(\g_0) = 3$ and $r(\g_1)=2$.   Introduce the \emph{infinite propeller} and its closure in $\mW$:
\begin{equation}\label{eq-infpropeller}
P_{\g} ~ \equiv ~  \cZ_{\g} ~ \cup ~ \bigcup_{\delta > 0} ~ P_{\g^\delta}   \quad , \quad \overline{P}_{\g} ~ \equiv ~ \overline{\bigcup_{\delta > 0} ~ P_{\g^\delta}}
\end{equation}
\end{defn}

 As observed in \cite{HR2016}, the closure $\overline{P}_{\g}$ of an
 infinite propeller contains the Reeb cylinder $\cR$,    with $\overline{P}_{\g} ~ = ~ P_{\g} ~ \cup ~ \cR$.

\begin{figure}[!htbp]
\centering
{\includegraphics[width=54mm]{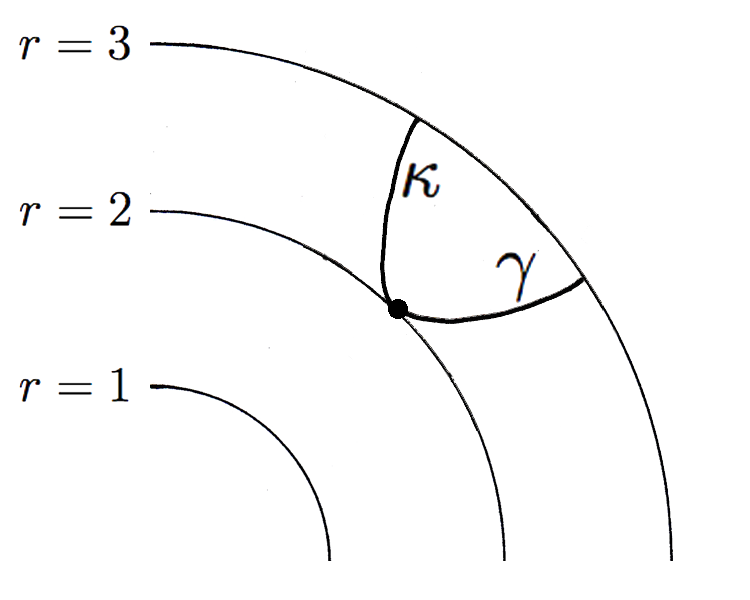}
\caption{\label{fig:Gamma1}   A curve $\G=\g\cup \kappa$ in  $\partial_h^-\mW$}}
\vspace{-6pt}
 \end{figure}

We now introduce the notion of double propellers in $\mW$.
  Consider a smooth curve $\G \subset \partial_h^-\mW$ parametrized by $\G \colon [0,2] \to \partial_h^- \mW$, with the notation $\G_s = \G(s)$, such that:  
  \begin{enumerate}
\item $r(\G_s) \geq 2$ for all $0 \leq s \leq 2$; 
\item  $r(\G_0) = r(\G_2) = 3$, so that both endpoints lie in the boundary $\partial_h^- \mW \cap \partial_v \mW$; 
 \item $\G$ is  topologically transverse to the   cylinders $\cC(r)$ for $2 \leq r \leq 3$,  except   at the midpoint $\G_1$.
\end{enumerate}
It then follows that   $r(\G_s) \geq  r(\G_1) = 2+\delta$ for all $0 \leq s \leq 2$, and some $\delta\geq 0$.
 See Figure~\ref{fig:Gamma1} for an illustration in the case when $\delta = 0$.

 Assume that $\delta > 0$, so that $r(\G_s) > 2$  for all $0 \leq s \leq 2$,   then    the $\cW$-orbit of each $\G_s$  traverses $\mW$. 
 The $\Psi_t$-flow of the points in $\G$   form a compact   surface embedded in $\mW$,   whose boundary is contained in the boundary of $\mW$, and thus the surface separates $\mW$ into two connected components. This surface is denoted   $P_{\G}$ and called  the \emph{double propeller} defined by the $\Psi_t$-flow of $\G$.

Consider the   curves $\g, \kappa \subset   \partial_h^- \mW$ obtained by dividing 
  the  curve $\G$ into two  segments at the midpoint $s=1$.  Parametrize these curves as follows: 
\begin{eqnarray*}
\g = \G \ | \ [0,1] ~ & , &  ~   \G(s) ~ \text{for} ~ 0 \leq s \leq 1 \\
\kappa= \G \  | \ [1,2] ~ & , &  ~  \G(2-s) ~ \text{for} ~ 0 \leq s \leq 1 
\end{eqnarray*}
The orbit  $\ds \{\Psi_t(\G_1) \mid 0 \leq t \leq T_1\}$ forms the \emph{long boundary} of the propellers $P_{\g}$ and $P_{\kappa}$ generated by the $\cW$-flow of these curves. Then $P_{\G}$ is viewed as the gluing of   $P_{\g}$ and $P_{\kappa}$ along the long boundary, hence  the notation ``double propeller'' for $P_{\G}$. 

If $\delta =0$, define   two infinite propellers   $P_{\g}$ and $P_{\kappa}$ 
as in Definition~\ref{def-infpropeller}, and then define 
    $\ds P_{\G} = P_{\g} \cup P_{\kappa}$, 
 where   the  $\Psi_t$-orbit $\cZ_{\g}$  of the midpoint $\G_1$,    defined as in \eqref{eq-Zset},  is again common to both $P_{\g}$ and $P_{\kappa}$, and  $P_{\G}$ is viewed as the gluing of the two propellers along an ``infinite zipper''.

\begin{defn}\label{def-doublepropeller}
Let $\G$ be as above with $r(\G_1)=2$. Let $\g^{\delta}$  and $\kappa^{\delta}$
for $0 <\delta\leq 1$ be the curves as defined in \eqref{eq-sparem}. The  \emph{infinite double propeller} is the union:
\begin{equation}\label{eq-doublepropeller}
P_{\G} ~ \equiv ~ \cZ_{\g} ~ \cup ~ \bigcup_{\delta > 0} ~ \left\{P_{\g^{\delta}} \cup P_{\kappa^{\delta}}\right\}  
\end{equation}
\end{defn}

\section{Global dynamics for $\e<0$}\label{sec-enegative}

We now analyze   the global dynamics of the flows $\Phi_t^{\e}$ on $\mK_\e$   for the case when $\e < 0$.
Recall that the assumption $\e < 0$ means that the self-insertion maps $\sigma_i^\e$ for $i=1,2$ of the Wilson Plug $\mW$ do not penetrate far enough into $\mW$  to break the two periodic orbits of the Wilson flow $\Psi_t$ by trapping them with the attracting orbits to these periodic orbits, as is the case when $\e=0$.  As a consequence, we show that the flow $\Phi_t^{\e}$ has simple dynamics:  there is an invariant set that is a cylinder whose boundary components are the two
periodic orbits, and these are the only periodic orbits of the flow in the plug $\mK_\e$.

On the technical level, by the Parametrized Radius Inequality condition (K8), we have the strict inequality $r' < r$,   as illustrated in Figure~\ref{fig:modifiedradius}(A).
That is,  the radius coordinate is strictly increasing at a secondary entry
point and strictly decreasing at a secondary exit point. Thus, there
exists some $\Delta>0$, depending on $\e$, such that whenever an orbit
hits an entry point, the radius increases by at least $\Delta$. The main
consequence of this fact is that, unless an orbit hits the cylinder
$\tau(\cC)$ where is $\cC$ is the cylinder of radius 2 in $\mW$, its behavior is the same as in the Wilson plug. In spite of the simplicity of this portrait of the dynamics of the flow $\Phi_t^{\e}$, the proofs of these claims   uses key aspects on the analysis of the flow  $\Phi_t^{0}$ that were developed in \cite{HR2016}.

Note that as $\e < 0$ tends to $0$, the constant $\Delta$ also tends to zero, and as a consequence the length of a given non-trapped orbit in $\mK_{\e}$ grows increasingly long, as the number of possible reinsertions of the orbit through the faces of the insertions increases to infinity. The geometry of the parametrized flow $\Phi_t^{\e}$ as a function of $\e < 0$ is reminiscent of the Moving Leaf Lemma in \cite{EMS1977}  which is the key to analyzing the dynamics of   flows in counter-examples to the Periodic Orbit Conjecture, as constructed in \cite{EpsteinVogt1978,Sullivan1976}.

\subsection{Orbit lifting property}\label{subsec-neglemma}

The first step in the analysis the dynamics of the flow $\Phi_t^{\e}$ is  to establish the following key technical result, which   should be compared  with Propositions 6.5 and 6.7 of \cite{HR2016}.
Recall that we defined  $\rho_x^\e(t) \equiv r(\Phi_t^\e(x))$.
\begin{prop}\label{prop-negativemain}
Let $x$ be a primary or secondary entry point of $\mK_\e$ and $y$ an exit point with $x\equiv y$. Assume that either:
\begin{enumerate}
\item $r(x)>2$, 
\item $r(x)<2$ and $\rho_x^\e(t)\neq 2
$ for all $t\geq 0$,
\end{enumerate} 
then $x\prec_{\cK_\e} y$. Moreover, the collection of lifts of the $\cW$-arcs in $[x,y]_{\cK_\e}$
contains all the $\cW$-arcs of the $\cW$-orbit of $x'$ that are in
$\wmW$, where $\tau(x')=x$.
\end{prop}

\proof
Let $x=x_0$ and 
   $0 = t_0 < t_1 < \cdots < t_n < \cdots $ with $x_{\ell} =
   \Phi_{t_{\ell}}^\e(x)$ be  the transition points for the positive
   $\cK_\e$-orbit  of $x=x_0$. Note that if  the forward orbit of $x_0$ is finite, meaning that there exists $t>0$ such that $\Phi_t^{\e}(x_0)$  is a primary exit point, then there are a finite number of transitions points.
For $\ell \geq 0$, let $[x_{\ell}', y_{\ell +1}']_{\cW} \subset \wmW$ be the lift of the
$\cK_\e$-arc $[x_{\ell}, x_{\ell + 1}]_{\cK_\e} \subset \mK_\e$.

Assume that $y$ is not in the $\cK_\e$-orbit of $x_0$. We prove that this implies that there exists $s>0$ such that for every $t\geq s$ we have $n_{x_0}(t)>0$ and deduce a contradiction from this.

If $x_1$ is an exit point, then $[x_0',y_1']_\cW$ is a complete $\cW$-orbit traveling from $\partial_h^-\mW$ to $\partial_h^+\mW$. Thus, by definition we have $x_0'\equiv y_1'$, or equivalently that $x_0\equiv x_1=y$ and we obtain a contradiction to our assumption. Thus, we can assume that $x_1$ is a secondary entry point and $n_{x_0}(t_1)=1$.

Next, we prove that $n_{x_0}(t)\geq 0$ for all $t >0$. Suppose not, then  
there exists $t>0$ such that $n_{x_0}(t)<0$.  Let $\ell>1$ be the
first index such that $n_{x_0}(t_\ell)=0$ and
$n_{x_0}(t_{\ell+1})=-1$. Then $x_\ell$ and $x_{\ell+1}$ are both exit points,  so
by Proposition~\ref{prop-shortcut4} we obtain that $x_0'\prec_\cW
y_\ell'$ and $x_0\equiv x_{\ell+1}$, implying that $x_{\ell+1}=y$ which contradicts the assumption that  $y$ is not in the $\cK_\e$-orbit of $x_0$.
Thus $n_{x_0}(t) \geq 0$ for all $t >0$.
 
Next, suppose there exists $\ell>1$ such that $n_{x_0}(t_\ell)=0$, then let $\ell > 1$ be the least such index. If $x_{\ell}$ is a primary exit point, then by Proposition~\ref{prop-ee}, the claim of the Proposition follows. Otherwise, $x_{\ell +1}$ is defined, and 
both  $x_{\ell-1}$ and $x_{\ell}$ are   exit points,  so
by Proposition~\ref{prop-shortcut4} we obtain that $x_0'\prec_\cW y_{\ell +1}'$.
   Then by Proposition~\ref{prop-shortcut4} $x_1\equiv x_{\ell-1}$, and  it follows that  the $\cW$-arc $[x_\ell',y_{\ell+1}']_\cW$  is the $\cW$-arc following $[x_0',y_{1}']_\cW$  in the intersection of the $\cW$-orbit of $x_0'$ with $\widehat{\mW}_\e$.

We then apply this argument repeatedly to every   $\cK_\e$-arc   $[x_{\ell}, x_{\ell+1}]_{\cK_\e}$ for which $n_{x_0}(s) = 0$   for $t_{\ell} \leq s < t_{\ell+1}$, 
to conclude that $[x_{\ell}, x_{\ell+1}]_{\cK_\e}$ lifts to an arc $[x_{\ell}', y_{\ell+1}']_{\cW} \subset \wmW$. Note that for each such  $\ell$, 
the  $\cW$-arc  $[x_{\ell}', y_{\ell+1}']_{\cW}$ is contained in the
$\cW$-orbit of $x_0'$, and for distinct $\ell > 0$ these $\cW$-arcs
are disjoint. On the other hand,   the $\cW$-orbit of $x_0'$ is finite
as $r(x_0')\neq 2$ by assumption. We conclude that there is at most a finite number of indices $\ell$ such that $n_{x_0}(t_\ell)=0$. 
Thus, there exists $s_0>0$ such that $n_{x_0}(t)>0$ for all $t\geq s_0$.

There are now two possible cases to analyze:   either the $\cK_\e$-orbit of $x_0$ is finite and thus exits $\mK_\e$ at a primary exit point, or the $\cK_\e$-orbit of $x_0$ is infinite. 

Assume first that the   $\cK_\e$-orbit of $x_0$  is finite. Then there
exists $\ell>1$ such that $x_\ell$ is a primary exit point. Let
$n=n_{x_0}(t_{\ell-1})\geq 0$. Since $n_{x_0}(t)>0$ for $t\geq s_0$ we
can assume that $n>0$. Then 
  $x_{\ell-1}$ must be a secondary exit point. Otherwise,  the $\cW$-arc  $[x_{\ell-1}',y_\ell']_\cW$ goes from an entry to an exit point in $\mW$, which must then be facing. 
  This would imply that  $x_{\ell-1}$ is a primary entry point, contrary to assumptions. 
  
  Thus, $x_{\ell-1}$ is a secondary exit point and therefore
  $n_{x_0}(t_{\ell-2})=n+1$. Then there exists $k>0$, chosen to be
  the smallest index such that $n_{x_0}(t_{k})=n$ and $n_{x_0}(t)\geq
  n$ for every $t_k\leq t<t_\ell$, implying that
  $n_{x_0}(t_{k-1})=n-1$. Then $x_k$ is an entry point and
  Proposition~\ref{prop-shortcut4} applied to  the $\cK_\e$-arc $[x_k,
  x_{\ell}]_{\cK_\e}$  implies that $x_k\equiv x_\ell$. But then $x_k$
  is a primary entry point and thus $x_k=x_0$, which is again a contradiction as we assumed that $n_{x_0}(t_{k})=n>0$.

Finally, consider the case where the $\cK_\e$-orbit of $x_0$  is  infinite. Recall that we may assume that  $n_{x_0}(t)>0$ for all $t\geq s_0$.
We then have to consider   the two following situations: either the level function $n_{x_0}(t)$ admits an upper bound for $t \geq 0$, or the function $n_{x_0}(t)$ is unbounded.

Assume that $n_{x_0}(t)$ grows without bound. We show that this leads to a contradiction. Let $\ell_i$ be the first
index such that $n_{x_0}(t_{\ell_i})=i$, thus $\ell_0=0$,
$\ell_1=1$ and all the points $x_{\ell_i}$ are secondary entry points. We claim that
\begin{equation}\label{eq-radiusunbounded}
\rho_{x_0}^\e(t_{\ell_i})\geq r(x)+ n_{x_0}(t_{\ell_i}) \cdot \Delta  =r(x)+i \cdot \Delta  \ .
\end{equation}
Let us prove the claim by recurrence on $i$. For $i=1$ condition (K8)
in the construction of $\mK_\e$ implies that
$\rho_{x_0}^\e(t_1)\geq r(x)+\Delta$. Assume that for $i=k-1$ the inequality is
satisfied. If $\ell_k=\ell_{k-1}+1$ the claim follows from (K8). If
not, observe that $n_{x_0}(t_{\ell_k-1})=k-1$, thus by Proposition~\ref{prop-shortcut4} we have that $x_{\ell_{k-1}}'\prec_\cW y_{\ell_k-1}'$
and $\rho_{x_0}^\e(t_{\ell_{k-1}})=\rho_{x_0}^\e(t_{\ell_k-1})$. Since
$x_{\ell_k}$ is a secondary entry point we have that
$\rho_{x_0}^\e(t_{\ell_k})\geq r(x)+k \cdot \Delta$ proving the claim.
Since $\Delta>0$ and the function $\rho_{x_0}^\e(t)$ is bounded above by $3$, this is impossible, and thus the level function $n_{x_0}(t)$ must be bounded.

It remains to consider the case where there exists some $N > 0$ such that 
$n_{x_0}(t)\leq N$ for all $t\geq 0$.   Then there exists      $0 \leq n_0 \leq N$ which is  the least integer such that there exists    
     $0 < \ell_0 < \ell_1 < \cdots < \ell_k < \cdots$ such that  $n_{x_0}(t_{\ell_j}) = n_0$. That is, 
     $\ds n_0 = \liminf_{\ell > 0} ~ n_{x_0}(t_{\ell}) \leq N$. 
          Since $n_0$ is the least such integer, there exists $k \geq 0$ such that  $n_{x_0}(t) \geq n_0$ for all $t \geq t_{\ell_k}$. 
      Then the segment $[x_{\ell_k}, x_{\ell_{k+1} + 1}]_{\cK_\e}$ 
     satisfies the hypotheses of Proposition~\ref{prop-shortcut4}, so
     we have   $x_{\ell_k}' \prec_{\cW} y_{\ell_{k+1} + 1}'$. Thus the
     lifted $\cW$-arcs $[x_{\ell_m}', y_{\ell_m+1}']_\cW$ for $m\geq
     k$ belong to the $\cW$-orbit of $x_{\ell_k}'$.
Then by the estimates \eqref{eq-segmentlengths} and \eqref{eq-segmentlengths2} above, considering the lengths of the lifted  $\cW$-orbit segments $[x_{\ell_m}', y_{\ell_{m} + 1}']_{\cW}$  yields the   estimate  
$\ds  d_{min}\cdot (m-k)\leq L(2-r(x_{\ell_k}'))$ for all $m\geq k$. 
   However,   we can choose $m$ arbitrarily large, and so 
    also $(m-k)$, which yields a contradiction. 

We proved so far that $x\prec_{\cK_\e} y$. Observe that $x_1$ is a
secondary entry point and assume $y=x_n$, for some $n>0$. Let $\ell_i$
be all the indices $1\leq \ell_i\leq n$ such that
$n_{x_0}(t_{\ell_i})=0$. Then $n_{x_0}(t_{\ell_i-1})=1$ and
Proposition~\ref{prop-shortcut4} implies $x_1\equiv
x_{\ell_1-1}$. Thus the $\cW$-arc $[x_{\ell_1}',y_{\ell_1+1}']_\cW$ is
the second arc in the $\cW$-orbit of $x'$ that belongs to $\wmW$. The
same argument proves that $x_{\ell_1+1}\equiv x_{\ell_2-1}$ and $[x_{\ell_2}',y_{\ell_2+1}']_\cW$ is
the third arc in the $\cW$-orbit of $x'$ that belongs to $\wmW$. We
conclude that the arcs $[x_{\ell_i}',y_{\ell_i+1}']_\cW$, for all $i$,
are all the arcs in the $\cW$-orbit of $x'$ that belong to $\wmW$.
\endproof

We point out two important corollaries of Proposition~\ref{prop-negativemain} and its proof.

\begin{cor}\label{cor-negfinite}
Let $x\in \mK_\e$ be such that $\rho_x^\e(t)\neq 2$ for every $t\geq 0$. If $x$ is not a transition point, then the forward $
\cK_\e$-orbit of $x$ exits $\mK_\e$.
\end{cor}

\proof
Let $x\in \mK_\e$ be such that  $\rho_x^\e(t)\neq 2$ for every $t\geq 0$, and assume that the forward orbit never exits $\mK_\e$. Since for every $t\geq 0$ the $\cW$-orbit of $\tau^{-1}(\Phi_t^\e(x))$ is finite, the assumption implies that the forward $\cK_\e$-orbit of $x$ has infinitely many transition points. Then either the level function $n_{x_0}(t)$ is upper and lower bounded, or is not.  The proof of Proposition~\ref{prop-negativemain} gives a contradiction in both cases.
\endproof

\begin{cor}\label{cor-negperiodic}
The $\cK_\e$-orbit $\cO_i^\e$ containing $\tau(\cO_i)$, for $i=1,2$, is periodic.
\end{cor}

\proof
The arc $\tau(\cO_i)$, for $i=1,2$, intersects the entry region $E^{\e}_i$ at a
point $p_i^- \in E^{\e}_i$ whose radius coordinate satisfies $r(p_i^-) >  2$, as illustrated in Figure~\ref{fig:modifiedradius}(A). Then by the proof of Proposition~\ref{prop-negativemain}, the $\cK_{\e}$-orbit of $x_i$ contains the facing point   $y_i \in S^{\e}_i$, which is by construction the other endpoint of the arc $\tau(\cO_i)$.
\endproof

We end this section with a brief discussion of the trapped set of
$\mK_\e$ for $\e<0$. Corollary~\ref{cor-trapped} states that for $x\in \partial_h^-\mK_\e$ with
$r(x)=2$ the $\cK_\e$-orbit of $x$ is trapped. The same argument
can be applied to secondary entry points with radius 2. Thus if
$x$ is a primary entry point with $r(x)\neq 2$ and such that there
exists $t>0$ with $\rho_x^\e(t)=2$, the orbit of $x$ is trapped. Thus
the trapped set changes with $\e$: it is composed of the curve 
$$\{\tau(x')\in \partial_h^-\mK_\e \mid x'=(2,\theta,-1)\},$$
and the curve of primary entry points $x\in \partial_h^-\mK_\e$ such
that the $\cK_\e$-orbit of $x$ hits $E_i^\e$ in a
point with radius 2.
Observe that these conclusions hold for negative orbits also.

\subsection{Existence of exactly two periodic orbits}\label{subsec-periodic}

Corollary~\ref{cor-negperiodic} implies that the plugs $\mK_\e$ have
at least two periodic orbits, corresponding to the periodic orbits of the Wilson plug $\mW$. We next show that  these are the only
periodic orbits of the flow $\Phi_t^\e$.   Moreover, these orbits   form the boundary of an invariant cylinder
$\fM_\e \subset \mK_\e$, whose limit as $\e\to 0$ is the invariant open set $\fM_0\subset \mK_0$, introduced in \cite{HR2016} and discussed in further detail below.

The proof of the following result is based in spirit, and also in many details, on the proof of the aperiodicity of the plug $\mK_0$ as given in  \cite[Theorem 8.1]{HR2016}.
\begin{thm}\label{thm-2periodicnegative}
For $\e<0$, the flow $\Phi_t^\e$ has exactly two periodic orbits, $\cO_1^\e$ and $\cO_2^\e$.
\end{thm}

\proof
 
 Suppose there exist a periodic orbit inside $\mK_\e$ and let  $x$ be a point on  it.  
Let    $0 \leq t_0 < t_1 < \cdots < t_n < \cdots $ with $x_{\ell} = \Phi_{t_{\ell}}^\e(x)$ be  the transition points for the $\cK_\e$-orbit    $\{\Phi_t ^\e (x) \mid t \geq 0\}$, and suppose that  $t_n > 0$ is the first subsequent transition point   with $x_n = \Phi_{t_n}^\e (x_0) = x_0$, so that $\Phi_{t + t_n}^\e (x) = \Phi_{t}^\e (x)$ for all $t$.  

Let $[x_{\ell}', y_{\ell +1}']_{\cW} \subset \wmW$ be the lift of the
$\cK_\e$-arc $[x_{\ell}, x_{\ell + 1}]_{\cK_\e} \subset \mK_\e$, for
$0 \leq \ell \leq n$ and let $r_\ell$ be the radius coordinate of the
arc $[x_{\ell}', y_{\ell +1}']_{\cW}$. Periodicity of the orbit implies that $[x_{0}', y_{1}']_{\cW} = [x_{n}', y_{n+1}']_{\cW}$. Moreover, the   $r$-coordinate is constant on each $\cW$-arc 
$[x_{\ell}', y_{\ell +1}']_{\cW}$, so it has a minimal value
$r_0$. Without loss,   assume   this minimum occurs for  $[x_{0}',
y_{1}']_{\cW}$. 

Observe that since the radius is minimum on $[x_{0}',
y_{1}']_{\cW}$ and the radius strictly increases at secondary
entry points, $x_0$ is a secondary exit point and $x_1$ is a secondary entry point. Assume that
$n_{x_0}(t_{n-1})=k$ for some integer $k$, we next analyze the
different situations depending on the value of $k$.

If $k<0$, let $\ell>1$ be the first index such that
$n_{x_0}(t_\ell)=-1$, then $n_{x_0}(t_{\ell-1})=0$ and Proposition~\ref{prop-shortcut4} implies that $x_0'\prec_{\cW}y_\ell'$. Thus
$r_{\ell-1}=r_0$ and $r_\ell<r_0$, contradicting the fact that $r_0$
is minimum. Thus $n_{x_0}(t)\geq 0$ for all
$t\geq0$.

If $k=0$, Proposition~\ref{prop-shortcut4} implies that $x_0'\prec_{\cW}y_n'$, and thus
$r_{n-1}=r_0=r_n$, that is a contradiction because the orbit
cannot pass a transition point with the radius coordinate staying
constant.

If $k>0$, let us consider the toy case where the level is strictly
increasing, or in other words, the case where all the transition
points $x_\ell$ for $1\leq \ell \leq n-1$ are secondary entry
points. Since $x_0=x_n$ is a secondary exit point,
Proposition~\ref{prop-shortcut4} implies that
$x_{n-1}\equiv x_n$ and $x_{n-2}'\prec_{\cW} y_{n+1}'=y_1'$. Thus
$r_{n-3}<r_{n-2}=r_n=r_0$, contradicting the minimality of $r_0$.

If the level is not strictly increasing, consider first the case when $k>1$. We use a technique introduced by Ghys in \cite[page 299]{Ghys1995}, which   defines a monotone increasing function $i(a)$ derived from the level function $n_{x_0}(t)$.   For
$0\leq a\leq k-1$, set $i(a)$ such that $n_{x_0}(t_{\ell_{i(a)}})=a$
and $n_{x_0}(t)\geq a$ for all $t\geq t_{\ell_{i(a)}}$. Then
$0=i(0)<i(1)<\cdots <i(k-1)<n-1$ and $x_{\ell_{i(a)}}$ is a secondary
entry point for all $a$. Moreover, $x_{\ell_{i(k-1)}}\prec_\cW x_n$, by
Proposition~\ref{prop-shortcut4}, and thus
$r_{\ell_{i(k-1)}}=r_n=r_0$. Then
$r_0=r_{\ell_{i(0)}}<r_{\ell_{i(1)}}<\cdots<r_{\ell_{i(k-1)}}=r_0$, a contradiction.

 \begin{figure}[!htbp]
\centering
{\includegraphics[width=80mm]{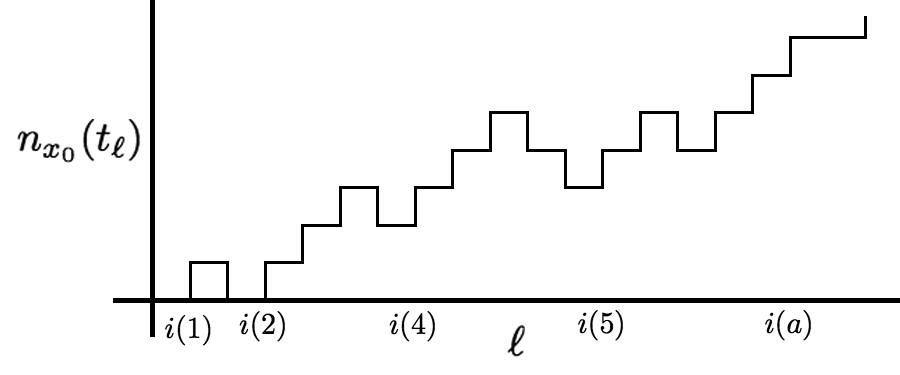}}
\caption{The functions $n_{x_0}$ and $i(a)$  \label{fig:steps} }
\vspace{-6pt}
\end{figure}

We are left with the case $k=1$ and $n_{x_0}(t_n)=0$. Thus
$1=n_{x_0}(t_1)=n_{x_0}(t_{n-1})$ and $x_1'\prec_\cW y_n'$. Moreover,
since $x_1$ is a secondary entry point and $x_n$ is a secondary exit
point, then $x_1\equiv x_n=x_0$. This implies that the in $\mW$, $x_0'\prec_\cW x_1'\prec_\cW x_0'$ and hence the arc $[x_0',y_1']_\cW$ is one of
the arcs $\cO_i\cap \mW'_\e$ for $i=1,2$. 

Thus the only  periodic orbits are the $\cK_\e$-orbits of $\tau(\cO_i)$ for
$i=1,2$. 
\endproof

The discussion of the trapped set at the end of Section~\ref{subsec-neglemma} and  Theorem~\ref{thm-2periodicnegative}, yields the following assertion in Theorem~\ref{thm-main1}.
\begin{cor}\label{cor-wandering}
All orbits of the flow $\Phi_t^\e$  are wandering, except for the periodic orbits $\tau(\cO_i)$ for $i=1,2$.
\end{cor}

\subsection{Invariant sets}\label{subsec-invariant}

The flow $\Psi_t$ on the modified Wilson Plug $\mW$ preserves the   Reeb cylinder $\cR\subset \mW$ introduced in Section~\ref{subsec-wilson}. The restriction of the flow to this invariant set consists of the two periodic orbits $\cO_i$, $i=1,2$, for  the flow, along with orbits asymptotic to these orbits, as illustrated by the restriction of the flow to the central band in  Figure~\ref{fig:flujocilin}(C).

\begin{figure}[!htbp]
\centering
{\includegraphics[width=60mm]{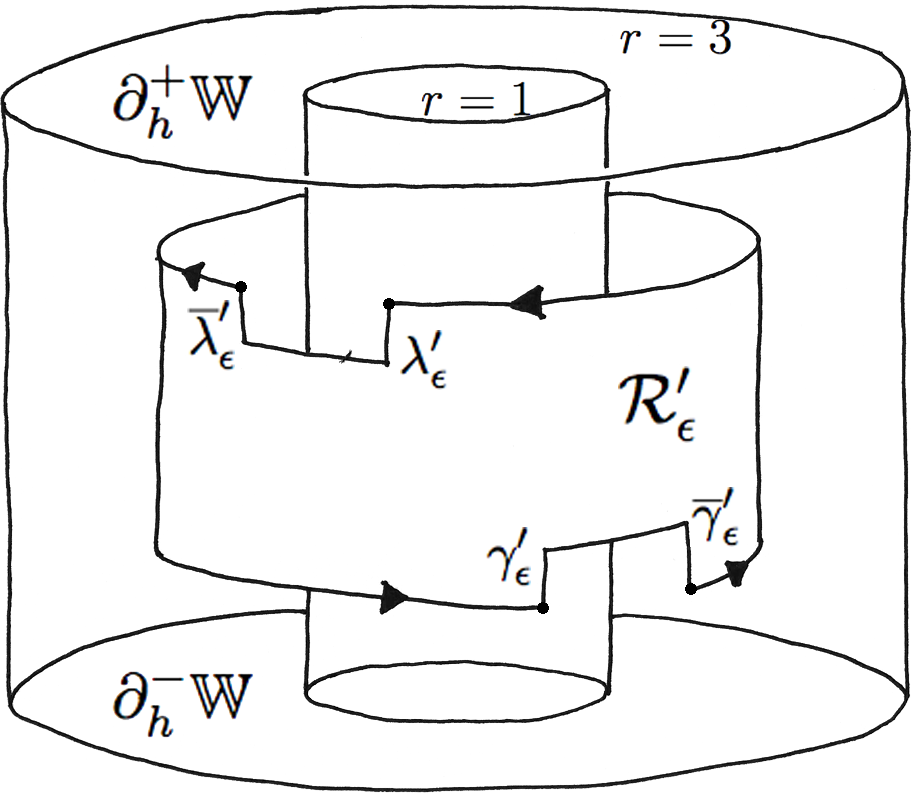}}
\caption{\label{fig:notches} The notched cylinder $\cR'$ embedded in $\mW$}
\vspace{-6pt}
\end{figure}

 Consider the intersection $\cR_\e' = \cR \cap \wmW$, where $\wmW$ is the closure of $\mW'_\e$ as defined in \eqref{eq-notchedW}. 
Then $\cR_\e'$ is a compact submanifold of $\cR$ with boundary,  which is the union of the two arcs   $\cO_i\cap \mW'_\e$ which are contained in the $\cW$-periodic orbits, and the boundary of the ``notches'', that is the intersections of $\cR$ with $\cD_i^\e$ for $i=1,2$. 
This set is  illustrated in Figure~\ref{fig:notches}. 
 We consider in this section the image of $\tau(\cR_\e')$ under the flow   $\Phi_t^\e$.

First, recall that for $\e =0$,   the images $\tau(\cO_i \cap \wmW) \subset \mK_0$ of these orbits are not periodic. 
 One of the main results in \cite{HR2016} is that under some
generic assumptions, the unique minimal set for the flow $\Phi_t^0$ equals the closure of the $\cK_0$-flow of the image of the Reeb cylinder
$\tau(\cR_0')$, and has the structure of  a \emph{zippered lamination}.  For a precise definition of a zippered lamination
we refer to   \cite[Chapter~19]{HR2016}.

Observe that the set $\tau(\cR_\e')  \subset \mK_\e$  is not   invariant under the flow $\Phi_t^\e$    since the Reeb cylinder $\cR$ intersects the inserted regions $\cD_i^\e$ for $i=1,2$. 
Instead, we must consider the submanifold $\fM_\e$ obtained by  applying the flow $\Phi_t^\e$  to the image  $\tau(\cR_\e')$: 
\begin{equation}\label{eq-invariantMe}
\fM_\e ~ \equiv ~  \{\Phi_t^\e(\tau(\cR_\e')) \mid -\infty < t < \infty \} \subset \mK_\e ~.
\end{equation}
For $\e< 0$,  Corollary~\ref{cor-negperiodic} implies that $\cK_\e$-orbits of the boundary segments  $\tau(\cO_i\cap \mW'_\e)$ are periodic orbits. In fact, a much stronger statement is true.

\begin{prop}
For $\e<0$, the set $\fM_\e$ is homeomorphic to a cylinder.
\end{prop}

\proof

Consider the intersections of $\cR$ with $\cD_i^\e$ for $i=1,2$. The
notches formed by deleting the intersection $\cD_1^\e$  contains 3
boundary curves: label the curve in $\cL_1^{\e-}$ by  $\g_\e'$. The facing curve transverse to the flow is labeled by $\ovg_\e'$. Analogously, the notches formed by deleting the intersection $\cD_2^\e$  contains 3 boundary curves,
where    $\lambda_\e'$  labels the curve transverse to the flow in $\cL_2^{\e-}$, and $\ovl_\e'$ is the facing curve.  These curves are 
   illustrated in Figure~\ref{fig:notches}.

The projections of the curves $\g_\e'$ and $\lambda_\e'$ to $\mK_\e$ are denoted by $\g_\e = \tau(\g_\e')$ and $\lambda_\e = \tau(\lambda_\e')$. These are curves in the entry regions $E^{\e}_1$ and $E^{\e}_2$, respectively. 
In the same way consider the  curves $\ovg_\e=\tau(\ovg_\e')$ and $\ovl_\e=\tau(\ovl_\e')$ which are contained in the exit regions $S^{\e}_1$ and $S^{\e}_2$.

The Parametrized Radius Inequality for $\e<0$ implies that each point on these four curves lies in the
region $\{r> 2\}$. Thus, by Proposition~\ref{prop-negativemain} the $\cK_\e$-orbit of each point
$x\in \g_\e$  passes through the facing point $y\in \overline{\g_\e}$.
Thus, the $\Phi_t^\e$-flow of the curve $\g_\e$ generates  a compact surface denoted by $\Sigma_{\g_\e} \subset \mK_\e$ that ``fills the
notch'' created by deleting  the rectangular region $\cR \cap \cD_1^\e$.  Analogously, the $\Phi_t^\e$-flow of the curve $\lambda_\e$ generates  a compact surface denoted by $\Sigma_{\lambda_\e} \subset \mK_\e$ that ``fills the
notch''   created by deleting  the rectangular region $\cR \cap \cD_2^\e$.  Both of these surfaces with boundary are diffeomorphic to a rectangular region, but the embeddings of this region in $\mK_\e$ becomes increasingly ``twisted'' as $\e \to 0$. In any case, the surface $\fM_\e$ obtained by attaching these surfaces  $\Sigma_{\g_\e}$  and $\Sigma_{\lambda_\e}$ along their boundary to the boundary of the notches in the image $\tau(\cR_\e')$ is homeomorphic to a cylinder whose boundary components are the periodic
orbits $\cO_i^\e$ for $i=1,2$.
\endproof

Let us describe the surfaces $\Sigma_{\g_\e}$ and  $\Sigma_{\lambda_\e}$. The description below proceeds analogously to the surface $\fM_0$ obtained from the $\Phi_t^0$-flow of the Reeb cylinder   $\cR_0' \subset \mK_0$ which is described in complete detail in Chapter~12 of \cite{HR2016}. We describe the surface $\Sigma_{\g_\e}$ below, and the description of the surface $\Sigma_{\lambda_\e}$ is exactly analogous.

Consider the curve $\g_\e \in E^{\e}_1$ and the corresponding curve $\g_\e'' = \tau^{-1}(\g_\e) \subset  L_1^-$. 
Note that the curves $\g_\e'$ and $\g_\e''$ satisfy $\tau(\g_\e') = \tau(\g_\e'') = \g_\e$, but the curve $\g_\e'$ lies in the cylinder of points $r=2$ in $\mW$, while the curve 
$\g_\e''$ is contained in the bottom face of $\mW$. In particular,   the radius coordinate along $\g_\e''$ is strictly
greater than 2.  Let $r_1\geq 2+\Delta$ be the minimum radius attained along the
curve. Then the $\Psi_t$-flow of $\g_\e''$ generates a finite propeller $P_{\g_\e}$ as
described in Section~\ref{sec-propeller}. Consider the notched propeller
$P_{\g_\e}'=P_{\g_\e}\cap \mW'_\e$. Thus $\tau(P_{\g_\e}')\subset \fM_\e$
and is attached to $\tau(\cR_\e')$ via the insertion map $\sigma_1^\e$.

Since this propeller is finite, its intersection with $E^{\e}_1$
and $E^{\e}_2$ is along a finite number of curves, whose radius
coordinate is bounded below by $r_1+\Delta$. Each of these curves generates a
finite propeller, simple or double, that is attached to $\tau(P_{\g_\e}')$ via the map
$\sigma_i^\e$, for $i=1$ or 2. Iterating this process a
finite number of steps, we construct the compact surface $\Sigma_{\g_\e}$. 
Analogously, we obtain the surface $\Sigma_{\lambda_\e}$. We conclude that $\fM_\e$ is formed from the union of 
$\tau(\cR')$  with a finite number of finite propellers. An important
observation is that the points in $\fM_\e$ with radius coordinate
equal to 2 are exactly those in $\tau(\cR_\e')$.

\begin{remark}
The propellers in the construction might have ``bubbles'' formed by
double propellers, as described in
Chapters~15 and 18 of \cite{HR2016}. These are compact surfaces that are
attached to internal notches, which are  notches in a propeller that do
not intersect the boundary of the propeller.
\end{remark}

Observe that as $\e\to 0$, the first propeller attached to $\tau(\cR_\e')$ becomes arbitrarily
long, and when $\e=0$ we obtain  an infinite propeller. The set
$\fM_0$ is no longer a cylinder, and has a very complicated
structure. We can thus see $\fM_0$ as the limit of the embedded cylinders $\fM_\e$.
This phenomenon is analogous to the behavior of the   leaves for a compact foliation as they approach its bad set, as described in \cite{EMS1977}.

To complete the proof of   Theorem~\ref{thm-main1}, choose a $C^{\infty}$-family of embeddings $\sigma_i^\e$ which satisfy the Parametrized Radius inequality, for $-1 \leq \e \leq 0$. Let $\Phi_t^\e$ be the resulting flows on the plug $\mK_{\e}$. Then the results of Section~\ref{sec-enegative} show that these flows satisfy the assertions of Theorem~\ref{thm-main1}.

\section{Geometric hypotheses for $\e \geq 0$}\label{sec-generic}

  The  dynamical properties  of the Kuperberg flows $\Phi_t^\e$ when $\e \geq 0$ are far more subtle than when   $\e < 0$, and to obtain our results on the global dynamics of $\Phi_t^\e$ when $\e \geq 0$ 
  requires that we impose a variety of additional  hypotheses  on the construction of $\Phi_t^\e$. 
In this section, we formulate some basic additional geometric assumptions for  the insertion maps $\sigma_i^\e$ beyond what we specified in Section~\ref{subsec-kuperberg}. The basic point is to  require  that they have   the geometric shape which is intuitively 
implicit in  Figures~\ref{fig:insertiondisks}, \ref{fig:twisted} and \ref{fig:modifiedradius}(C).    
 
First, we note   a straightforward
consequence of the Parametrized Radius Inequality (K8) in Section~\ref{subsec-kuperberg}. 
Recall that $\theta_i$ is the radian angle coordinate specified in (K8) such that for $x' = (2, \theta_i,  -2) \in  L_i$ we have 
$r(\sigma_i^\e(2, \theta_i,  -2)) = 2+\e$.

\begin{lemma}\label{lem-r0}
For $\e>0$ there exists $2+\e <r_\e<3$ such that $r(\sigma_i^\e(r_\e,\theta_i, -2))= r_\e$.
\end{lemma}

\proof
Since  $r(\sigma_i^\e(2,\theta_i, -2))= 2+\e$ and $r(\sigma_i^\e(3,\theta_i, -2))<3$, by the continuity in $r$ of the function  
$r(\sigma_i^\e(r,\theta_i, -2))$  we conclude that there
exists $2+\e <r_\e<3$ such that $r(\sigma_i^\e(r_\e,\theta_i, -2))= r_\e$.
\endproof

 We then add an additional assumption on the   insertion maps $\sigma_i^\e$ for $i=1,2$ that the radius  is \emph{decreasing} under the insertion map, for $r \geq r_\e$.
 
 \begin{hyp} \label{hyp-monotone}
If $r_\e$ is the smallest  $2+\e <r_\e<3$  such that
$r(\sigma_i^\e(r_\e,\theta_i, -2))= r_\e$, then assume that $r(\sigma_i^\e(r,\theta_i, -2)) < r$ for $r > r_\e$.  
 \end{hyp}

Next, we introduce an  hypothesis on the insertion maps   $\sigma_i^\e$ for $i=1,2$ which, in essence, guarantees that the images of the level curves for $r'=c$ are ``quadratic'' for $c$ near the value $r=2$, as pictured in Figure~\ref{fig:modifiedradius}. 
Recall our notational conventions. For    $i=1,2$, let  $x' = (r', \theta',  -2) \in  L_i^-$   denote a point in the domain of $\sigma_i^\e$   and denote its image by  
 $(r, \theta,z)  = \sigma_i^\e(x') \in \cL_i^{\e-} \subset \mW$.

  Let  $\pi_z \colon \mW \to \partial_h^- \mW$   denote the projection   along the $z$-coordinate, so $\pi_z(r, \theta, z) = (r, \theta,-2)$.

We first assume that $\sigma_i^\e$ restricted to the bottom face, $\sigma_i^\e \colon L_i^- \to \cL_i^{\e-} \subset \mW$,  has image  transverse to the vertical fibers of $\pi_z$.  
 Then  $\pi_z \circ \sigma_i^\e \colon L_i^- \to \mW$ is a diffeomorphism into the face $\partial_h^- \mW$. Denote the image set of this map by    $\fD_i \subset \partial_h^- \mW$. Then we can define 
 the inverse map 
 \begin{equation}\label{eq-inverse}
\Upsilon_i^\e = (\pi_z \circ \sigma_i^\e)^{-1} \colon \fD_i \to L_i^- \ .   
\end{equation}
In particular,  express the inverse map $x' = \Upsilon_i^\e(x)$ in polar coordinates as:
\begin{equation}\label{eq-coordinatesVT}
x' = (r',\theta', -2) = \Upsilon_i^\e(r, \theta,-2) = (r(\Upsilon_i^\e(r,\theta,-2)), \theta(\Upsilon_i^\e (r,\theta,-2)), -2)  = (R_{i,r}^\e(\theta)  , \Theta_{i,r}^\e(\theta), -2)  ~.
\end{equation}
 
In the following hypothesis, we    impose uniform conditions on the derivatives of the maps $\Upsilon_i^\e$. Recall that $0 < \e_0 < 1/4$ was specified in Hypothesis~\ref{hyp-genericW}, and we assume that $0 < \e < \e_0$.

 \begin{hyp} [{\it Strong Radius Inequality}]  \label{hyp-SRI}
  For   $i=1,2$, assume that: 
  \begin{enumerate}
\item \label{item-SRI-2} $\sigma_i^\e \colon L_i^- \to \mW$ is transverse to the fibers of $\pi_z$; 
\item \label{item-SRI-1} $r = r(\sigma_i^\e(r',\theta',-2))< r+\e$, except for  $x' = (2,\theta_i, -2)$ and  then  $r=2+\e$; 
  \item \label{item-SRI-3} $\Theta_{i,r}^\e(\theta)$ is an increasing function of $\theta$ for each fixed $r$; 
\item  \label{item-SRI-4}  For $2-\e_0 \leq r \leq 2+\e_0$ and  $i = 1,2$, assume that   $R_{i,r}^\e(\theta)$ has non-vanishing derivative, except when    $\theta = \overline{\theta_i}$ as 
  defined by  $\Upsilon_i^\e(2+\e,\otheta_i,-2)= (2,\theta_i,-2)$;   
 \item    For $2-\e_0 \leq r \leq 2+\e_0$      and $\theta_i -\e_0 \leq \theta \leq \theta_i + \e_0$ for  $i = 1,2$, assume that   
\begin{equation}\label{eq-quadratic}
\frac{d}{d\theta} \Theta_{i,s}(\theta) > 0 \quad, \quad  \frac{d^2}{d\theta^2} R_{i,s}(\theta) > 0.
\end{equation}
Thus for $2-\e_0 \leq r \leq 2+ \e_0$,  the graph of $R_{i,r}(\theta)$ is   parabolic with vertex       $\theta = \otheta_i$. 
\end{enumerate}
Consequently, each  surface $\cL_i^{\e-}$   is transverse to the coordinate vector fields $\partial/\partial \theta$ and $\partial/\partial z$ on $\mW$.
  \end{hyp}

 \begin{remark}
 Hypotheses~\ref{hyp-monotone} and \ref{hyp-SRI} combined, imply  that $r_\e$ is the unique value of $2+\e <r_\e<3$ for which 
 $r(\sigma_i^\e(r_\e,\theta_i, -2))= r_\e$. It follows that  the
 radius function  $\rho_x^\e$ at a secondary entry point $x$ with $r(x)>r_\e$ is strictly increasing.
 \end{remark}

\section{A pseudogroup model}\label{sec-pseudogroups}

The analysis  of the dynamical properties  of the standard Kuperberg flow $\Phi_t$,  as made in   \cite{HR2016}, was based on the introduction of an ``almost transverse'' rectangle $\bRt \subset \mK$, and  a detailed  study of the dynamics of the induced pseudogroup for the return map    $\whPhi$ of the flow to   $\bRt$. Our analysis of  the dynamics of the flows $\Phi_t^{\e}$  for the case when $\e > 0$ follows a similar approach.  We utilize the same rectangle $\bRt$ as   defined in \eqref{eq-goodsection} below, which is the same as in   \cite{HR2016},  but we use a simplified model for the pseudogroup, which incorporates   the induced return maps to $\bRt$ for both of the  flows $\Psi_t$ and $\Phi_t^\e$. The dynamics for the induced map $\whPsi$ from the Wilson flow is fairly straightforward to analyze, while that of the  induced map $\whPhie$ from the Kuperberg flow is extraordinarily complicated, so we restrict to analyzing the dynamics of selected maps defined by $\whPhie$, and   the task  is then more manageable.

The goal is to show that for $\e > 0$, the dynamics of the return map
$\whPhie$ contains disjoint families of ``horseshoes'', as will be
proved in Section~\ref{sec-horseshoe}, which then yields the
conclusions of  Theorem~\ref{thm-main2}.  The first step is to
introduce the pseudogroup $\whcGe$ and the $\psg$ $\whcGe^*$ in this
section, as Definitions~\ref{def-pseudogroup1} and
\ref{def-pseudogroup1a} below. In the next
Section~\ref{sec-horseshoe}, we  show that the action of
$\whcGe^*\subset \whcGe$  contains invariant  horseshoe dynamical subsystems.

 \subsection{A good rectangle}\label{subsec-crosssection}

We   introduce the ``almost transversal'' rectangle   to the flows $\Psi_t$ and  $\Phi_t^{\e}$ which is used for the study of the return dynamics of their flows, and for the construction of   associated pseudogroups.

\begin{figure}[!htbp]
\centering
{\includegraphics[width=100mm]{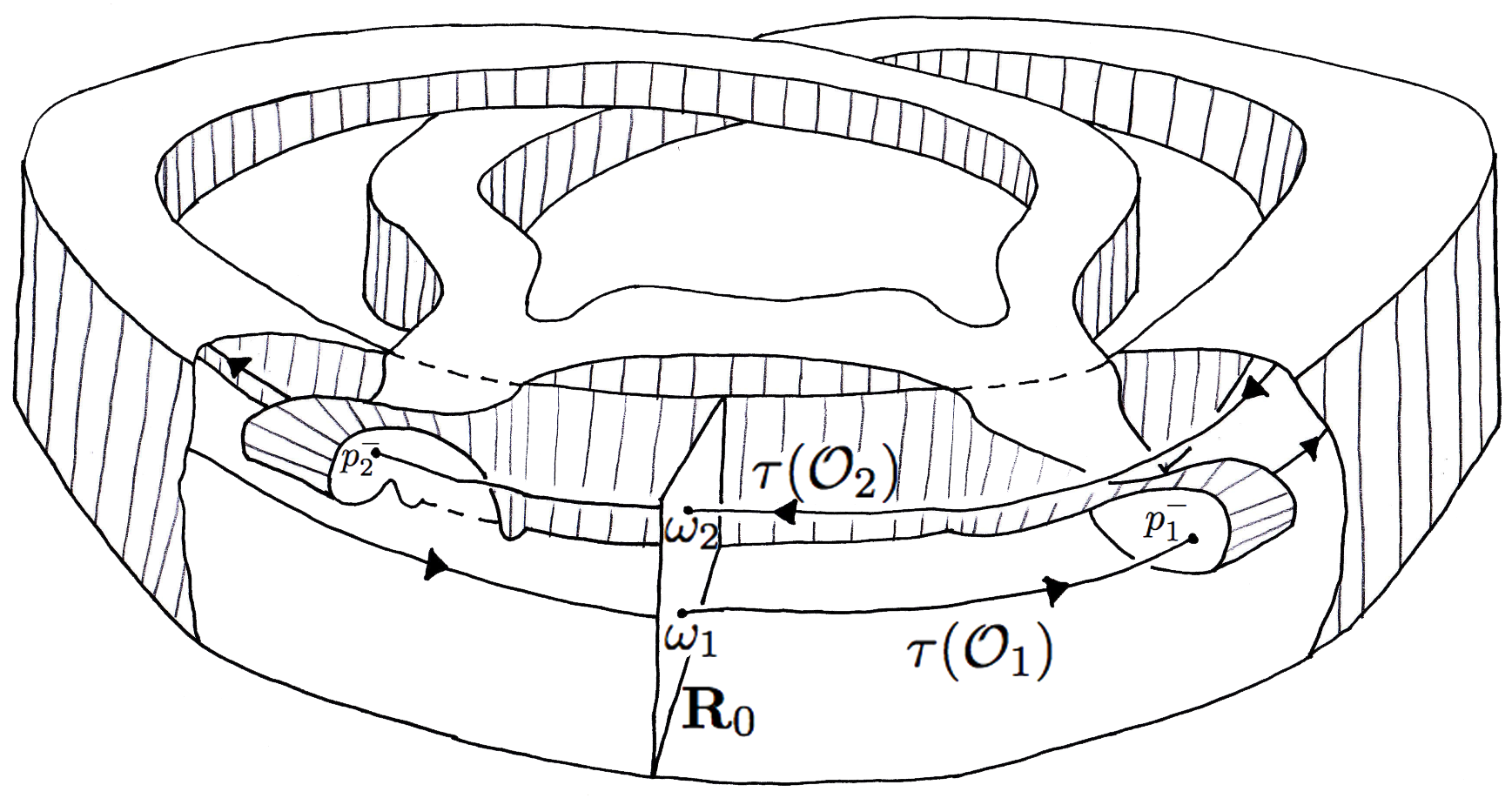}}
 \caption{\label{fig:KR} The rectangle $\bRt$ in the Kuperberg plug
   $\mK_\e$}
 \vspace{-6pt}
\end{figure}

Choose a value of $\theta_0$ such that the rectangle $\bRt$ as defined in cylindrical coordinates, 
\begin{equation}\label{eq-goodsection}
\bRt \equiv \{\xi = (r, \theta_0, z) \mid ~ 1 \leq r \leq 3 ~,   ~ -2
\leq z \leq 2\}  \,\subset \mW'_\e ~ , 
\end{equation}
is disjoint from both the regions $D_i$ and their insertions $\cD_i^\e$
for $i=1,2$, as defined in Section~\ref{subsec-kuperberg}. 
 For example, for the curves  $\alpha_i$ and $\beta_i'$  defined in
Section~\ref{subsec-kuperberg}, we take $\theta_0 = \pi$ so that
$\bRt$ is between the embedded regions  $\cD_i^\e$ for $i=1,2$ as illustrated in Figure~\ref{fig:KR}.    

As $\bRt \subset \mW'_\e$,  the quotient map $\tau \colon \mW \to \mK_\e$ is injective on $\bRt$. 
We use a slight abuse of notation, and also denote  the image    $\tau(\bRt) \subset \mK_\e$      by $\bRt$ with coordinates   $r = r(\xi)$ and $z = z(\xi)$  for  $\xi \in \bRt$. 
 
The periodic points in $\bRt$ for the Wilson flow are denoted by
 \begin{equation}\label{eq-omegas}
 \omega_1  = \cO_1 \cap \bRt = (2,\theta_0, -1)  \quad , \quad \omega_2 = \cO_2 \cap \bRt = (2,\theta_0, 1)   ~ .
\end{equation}
For $i=1,2$,   the first transition point for the forward orbit
of    $\omega_i$  is denoted by $p_i^{-} = \tau(\cL_i^{\e-} \cap \cO_i)$, and for the backward orbit    the first transition point is the special exit point 
 $p_i^{+} = \tau(\cL_i^{\e+} \cap \cO_i)$.

 Define a metric on   $\bRt$ by    $\ds d_{\bRt}(\xi, \xi') = \sqrt{(r' - r)^2 + (z'-z)^2}$,  for $\xi = (r, \theta_0, z)$ and 
$\xi' = (r', \theta_0, z')$.

We next introduce the first return map   $\whPsi$  on $\bRt$ for the Wilson flow $\Psi_t$. 
 The map $\whPsi$  is defined at $\xi \in \bRt$ if there is a $\cW$-orbit segment $[\xi, \eta]_{\cW}$ with $\eta \in \bRt$ and its interior $(\xi, \eta)_{\cW}$ is disjoint from $\bRt$. We then set $\whPsi(\xi) = \eta$. 
Thus, the domain of $\whPsi$ is the set:
\begin{equation}\label{eq-domainwhPsi}
 Dom(\whPsi) \equiv \left\{ \xi \in \bRt \mid \exists ~ t > 0 ~
    \text{ such ~ that} ~ \Psi_t(\xi) \in \bRt ~ \text{and}  ~
  \Psi_s(\xi)\notin \bRt ~ \text{for}~   0<s<t  \right\}.
\end{equation}
The radius function is constant  along the orbits of the Wilson flow,  so that $r(\whPsi(\xi)) = r(\xi)$ for all $\xi \in Dom(\whPsi)$. 
Also, note that the points $\omega_i$ for $i=1,2$ defined in \eqref{eq-omegas} are   fixed-points for $\whPsi$.
 For all  other points $\xi \in \bRt$ with $\xi \ne \omega_i$,  it was assumed in Section~\ref{subsec-wilson} that the function $g(r,\theta, z) > 0$, so the   $\cW$-orbit of $\xi$   has a ``vertical drift'' arising from the term 
$g(r, \theta, z)  \frac{\partial}{\partial  z}$ in the formula \eqref{eq-wilsonvector} for $\cW$. 

The precise description of the domain $ Dom(\whPsi)$ is discussed in detail in Chapter~9 of \cite{HR2016}, to which we refer the reader for further details. For our applications here, it suffices to note that the domain $ Dom(\whPsi)$ contains an open neighborhood of the vertical line segment $\cR \cap \bRt$. The dynamical properties   of $\Psi_t$ on $ Dom(\whPsi)$ are  described in Proposition~\ref{prop-wilsonproperties}, and illustrated in  Figures~\ref{fig:flujocilin} and \ref{fig:Reebcyl}.

Next, let   $\whPhie$ denote the first return map on $\bRt$ for the Kuperberg flow $\Phi_t^\e$. 
The domain of $\whPhie$ is the set:
\begin{equation}\label{eq-domainwhPhi}
Dom(\whPhie) \equiv \left\{ \xi \in \bRt \mid \exists ~ t > 0 ~ \text{ such ~ that} ~ \Phi_t^\e(\xi) \in \bRt ~ \text{and} ~  \Phi_s^\e(\xi)\notin \bRt 
  ~ \text{for} ~   0<s<t   \right\} .
\end{equation}
 The precise description of the domain $Dom(\whPhie)$ is very complicated, due to the nature of the orbits of $\Phi_t^\e$ as the union of   $\cW$-arcs for the Wilson flow. Moreover, 
the map $\whPhie \colon Dom(\whPhie) \to \bRt$ has many points of discontinuity, which arise when an orbit is tangent to section $\bRt$  along the line $\bRt\cap \cA$. 
 An extensive discussion of the properties of the return map $\whPhie$ for  the case $\e=0$  is discussed in detail in  Chapter~9 of \cite{HR2016}.  We  adapt these results as required, for the case   $\e > 0$.

 \subsection{The pseudogroup $\whcGe$}\label{subsec-pseudogroup}

A key idea, introduced in the work \cite{HR2016}, is to associate a pseudogroup $\cG_K$ to the return map $\whPhi$ and study the dynamics of this pseudogroup. This approach allows a more careful analysis of the interaction of the return maps $\whPsi$ and $\whPhi$ in determining the dynamics of $\Phi_t^\e$. In this paper, we work with   a   pseudogroup $\whcGe$ generated by   the return maps \emph{for both flows}, and then show that under the proper hypotheses, the orbits of $\whcGe$ in an invariant ``horseshoe'' subset of $\bRt$ agrees with the orbits of $\whPhie$. 

First, we recall the formal definition of a pseudogroup   modeled on a space $X$. 
\begin{defn}\label{def-pseudogroup}
A pseudogroup $\cG$ modeled on a topological space $X$ is a collection of  homeomorphisms between open subsets of $X$ satisfying the following properties:
\begin{enumerate}
\item For every open set $U \subset X$, the identity $Id_U \colon U \to U$ is in $\cG$.
\item For every $\varphi \in \cG$  with  $\varphi \colon U_{\varphi} \to V_{\varphi}$ where $U_{\varphi}, V_{\varphi} \subset X$ are open subsets of $X$, then    also  $\varphi^{-1} \colon V_{\varphi} \to U_{\varphi}$ is in $\cG$.
\item  For every $\varphi \in \cG$  with  $\varphi \colon U_{\varphi} \to V_{\varphi}$  and each open subset $U' \subset U_{\varphi}$, then the restriction $\varphi \mid U'$ is in $\cG$.
\item  For every $\varphi \in \cG$  with  $\varphi \colon U_{\varphi} \to V_{\varphi}$ and   every $\varphi' \in \cG$  with  $\varphi' \colon U_{\varphi'} \to V_{\varphi'}$, 
if $V_{\varphi} \subset U_{\varphi'}$ then the composition $\varphi' \circ \varphi$ is in $\cG$.
\item If $U \subset X$ is an open set,   $\{U_{\alpha} \subset X \mid \alpha \in \cA\}$ are open sets whose union is $U$, $\varphi \colon U \to V$ is a homeomorphism to an open set $V \subset X$ and for each $\alpha \in \cA$ we have $\varphi_{\alpha} = \varphi \mid U_{\alpha} \colon U_{\alpha} \to V_{\alpha}$ is in $\cG$, then $\varphi$ is in $\cG$.
\end{enumerate}
\end{defn}

 We first introduce a map which encodes a part of the insertion dynamics of the return map $\whPhie$.   The  reader interested in more more complete development of these ideas can consult   Chapter~9 of \cite{HR2016}.
 
 Let $U_{\phi_1^+} \subset Dom(\whPhie)$ be the
  subset of $\bRt$ consisting of points 
 $\xi \in Dom(\whPhie)$ with $\eta = \whPhie(\xi)$,  such that the $\cK_\e$-arc $[\xi, \eta]_{\cK_\e}$
  contains a single transition point $x$, with  $x \in E^{\e}_1$. 
  
Note that for such $\xi$, we see from Figures~\ref{fig:cWarcs} and ~\ref{fig:KR}, that its $\cK_\e$-orbit exits the surface $E^{\e}_1$ as the $\cW$-orbit of a point  $x' \in L_1^-$ with $\tau(x') = x$, flowing upwards from $\partial_h^- \mW$ until it   intersects $\bRt$ again. If the $\cK_\e$-orbit of $\xi$ enters   $E^{\e}_1$  but   exits   through $S^{\e}_1$ before crossing $\bRt$, then it is not considered to be in the domain $U_{\phi_1^+}$ as it contains more than one transition point between $\xi$ and $\eta$.

Let  $\phi_1^+ \colon U_{\phi_1^+} \to V_{\phi_1^+}$ denote the map defined by the restriction of $\whPhie$. 
As the $\cK_\e$-arcs $[\xi, \eta]_{\cK_\e}$ defining $\phi_1^+$ do not intersect $\cA$, the restricted map $\phi_1^+$ is continuous. 
Observe that the action of the map    $\phi_1^+$   corresponds to a flow through a transition point which increases the level function $n_x(t)$ by $+1$.

  We can now define the pseudogroup $\whcGe$ acting on $\bRt$ associated to the return maps $\whPsi$ and $\whPhie$.
 
   \begin{defn}\label{def-pseudogroup1}
  Let   $\whcGe$ denote the  pseudogroup generated by  the collection of all maps   formed by compositions of the maps 
 \begin{equation}\label{eq-generators}
\{Id, \phi_1^+,   \whPsi|U \mid U \subset  Dom(\whPsi) ~ {\rm is ~ open ~ and}~ \whPsi|U ~{\rm is~ continuous}\}
\end{equation}
   and their restrictions to open subsets in their domains.
\end{defn}

The notion of a  \emph{$\psg$} was introduced by Matsumoto in \cite{Matsumoto2010},  which is a subset of the maps in a pseudogroup that satisfy conditions (1) to (4) of  Definition~\ref{def-pseudogroup}, but 
need not satisfy the condition (5) on unions of maps.

     \begin{defn}\label{def-pseudogroup1a}
  Let   $\whcGe^*$ denote the    collection of all maps   formed by compositions of the maps 
 in \eqref{eq-generators} above,   and their restrictions to open subsets in their domains. 
\end{defn}

  Note that  $\whcGe^*$  is a $\psg$ contained in  $\whcGe$ but is not a pseudogroup itself. A key point in the proof of 
 Proposition~\ref{prop-iteration0} later in this work, is that for appropriately chosen $k \geq \ell(\e)$, there is a well-defined non-trivial element    $\vp_k=\whPsi^k\circ\phi_1^+ \in \whcGe^*$ which is a concatenation of a power of the Wilson return map $\whPsi$ with the first return map $\phi_1^+$ of the flow  $\Phi_t^\e$.

We conclude with an observation and a fundamental technical result which relates  the orbits of the map $\whPsi$   with those of the map $\whPhie$, that encodes  the existence of ``shortcuts'' as given in Proposition~\ref{prop-shortcut4}. 
We   require the following result, which is analogous  to Proposition~\ref{prop-negativemain}.
Recall that the constant $r_\e> 2$ was introduced in Lemma~\ref{lem-r0}, and that we assume that Hypothesis~\ref{hyp-monotone} is satisfied. 
\begin{prop}\label{prop-positivemain}
Let $x$ be a primary or secondary entry point of $\mK_\e$ and $y$ the exit point with
$x\equiv y$. Assume that $r(x)\geq r_\e$,
then $x \prec_{\cK_\e} y$, and the collection of lifts of the $\cW$-arcs in $[x,y]_{\cK_\e}$
contains all the $\cW$-arcs of the $\cW$-orbit of $x'$ that are in
$\wmW$, where $\tau(x')=x$.
\end{prop}
\proof
The proof follows in  the same way as that of Proposition~\ref{prop-negativemain}, where we note  that for  $r > r_\e$ the radius strictly grows along any orbit of $\Phi_t^{\e}$ when entering a face  of one of the insertions.
\endproof

 Introduce  the local coordinate function  $\wtr(\xi)$ which is defined
for points in $\bRt$ whose forward orbit hits $E_i^\e$, for
$i=1,2$, before reaching any other transition point or returning to $\bRt$. Let
$x\in E_i^\e$ be the first secondary entry point in the forward orbit
$\xi$, set $\wtr(\xi)=r(x)$. Then, for $\xi$ in the domain
$U_{\phi_1^+} \subset \bRt$  of the map $\phi_1^+ \in \whcGe$, we have 
$\wtr(\xi) =
r(\phi_1^+(\xi))$. 
Recall that  $\xi \in U_{\phi_1^+}$   if the forward $\cK_\e$-orbit of
$\xi$ hits $E^{\e}_1$ before returning to $\bRt$, so in particular
$U_{\phi_1^+}$ contains an open neighborhood of    the set
$\{\wtr=2\}\subset \bRt$.
 Next,  define  the  subset 
  \begin{equation}\label{eq-domainre}
U_{r_\e}  ~     =  ~ \{ \xi \in D(\whPhie)  \mid \wtr(\xi)  > r_\e\} \subset \bRt \ .
\end{equation}

 \begin{cor}\label{cor-positivemain}
 Let  $\xi  \in U_{r_\e}$ be such that   $\eta$ is
 contained in the forward $\cW$-orbit of $\xi$.  Then  there exists some $\ell > 0$ such that $\eta = (\whPhie)^{\ell}(\xi)$.
 \end{cor}
  \proof
   If the $\cW$-arc $[\xi, \eta]_{\cW}$ does not intersect an entry region $\cL_i^{\e-}$ then it is also a $\cK$-arc, and so the result follows. Otherwise, let $x \in \cL_i^{\e-}$ be the first transition point along $[\xi, \eta]_{\cW}$ and let $y \in \cL_j^{\e+}$ be the last exit point. Then $x\equiv y$, and $\xi \in U_{r_\e}$ implies that 
 $r(x)>  r_\e$, so by  Proposition~\ref{prop-positivemain} we have that
 $y$ and $\eta$ are in the forward $\cK_\e$-orbit of $\xi$. Thus,  there exists some $\ell > 0$ such that $\eta = (\whPhie)^{\ell}(\xi)$. 
 \endproof

\section{Horseshoes for the pseudogroup dynamics}\label{sec-horseshoe}

  In this section, we show that for $0<\e<\e_0$ sufficiently small, the  $\psg$   $\whcGe^*$ contains a map $\varphi$ with ``horseshoe dynamics''.   In fact, the precise statement, Theorem~\ref{thm-horseshoe} below, gives an even stronger conclusion, showing that there exists a countable collection of such maps, each with an invariant Cantor set for which the action is coded by a shift map, and each disjoint from the others. 
The proof of Theorem~\ref{thm-horseshoe} will be via a construction, which uses   the relations between the action of the maps $\whPsi$ and $\phi_1^+$ on $\bRt$, and the geometry of the traces of propellers and double propellers in $\bRt$.

\subsection{Traces of propellers}\label{subsec-props}

We describe the surfaces in the Wilson plug $\mW$ that are generated by the Wilson flow of curves in the bottom face $\partial_h^-\mW$, and especially those curves which cross the circle $\{r=2\}\subset \partial_h^-\mW$.  These surfaces are similar to the double propellers introduced in Definition~\ref{def-doublepropeller}.

\begin{figure}[!htbp]
\centering
{\includegraphics[width=50mm]{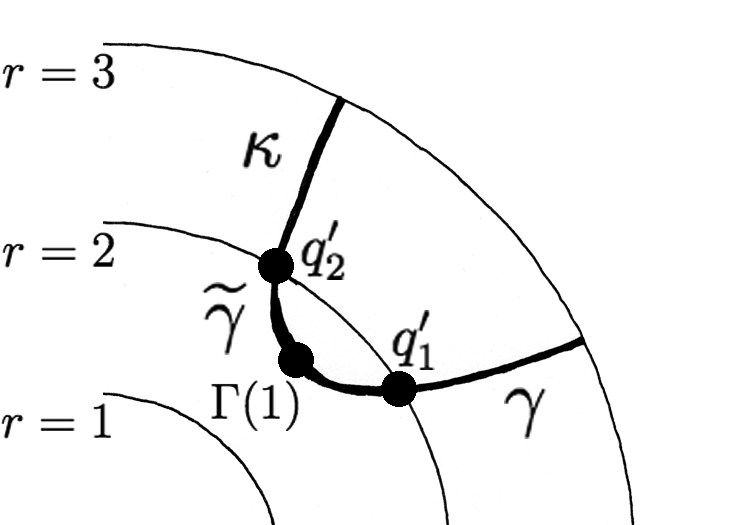}
\caption{\label{fig:Gammapositive}The curve $\G=\g\cup\wtg\cup
  \kappa$ in $\partial_h^-\mW$}}
\end{figure}

  Let $\G:[0,2]\to L_1^-\subset \partial_h^-\mW$ be a curve such that:
\begin{itemize}
\item $r(\G(0))=r(\G(2))=3$;
\item $r(\G(t))>r(\G(1))$ for all $t\neq 1$ and $r(\G(1))<2$;
\item $\G$ is transverse to the circle $\{r=2\}$.
\end{itemize}
We will call such a curve a \emph{traversing curve}.
Let $t_1<t_2$ be such that $r(\G(t_1))=r(\G(t_2))=2$. We can divide $\G$ in three overlapping parts,  as illustrated in Figure~\ref{fig:Gammapositive}:

\begin{itemize}
\item $\g$ the curve $\G([0,t_1])$;
\item $\wtg$ the curve $\G([t_1,t_2])$;
\item $\kappa$ the curve $\G([t_2,2])$.
\end{itemize}

Label the interior endpoints of these curves as $q_1' = \G(t_1) \in \g \cap \wtg$ and $q_2' = \G(t_2) \in \wtg \cap \kappa$.
Note that the forward orbits of the   points $\{q_1', q_2'\} \subset \partial_h^-\mW$  by the flow $\Psi_t$  spiral upward to the periodic orbit $\cO_1$, and so are trapped in the Wilson plug $\mW$.

Recall from Section~\ref{sec-radiuslevel} that we say points $x' \in \partial_h^- \mW$ and $y' \in \partial_h^+ \mW$ are \emph{facing}, and we  write $x' \equiv y'$,  
if      $x' = (r, \theta, -2)$ and $y' = (r, \theta, 2)$ for some $r$ and $\theta$.  Let $\ovG \subset \partial_h^+ \mW$ be the facing curve to $\G$, and denote by $\ovg$, $\ovk$ and $\ovtg$ the corresponding segments of $\ovG$ facing to the curves $\g$, $\kappa$ and $\wtg$. Then also define $\oq_1' = \ovG(t_1) \in \ovg \cap \ovtg$ and $\oq_2' = \ovG(t_2) \in \ovtg \cap \ovk$. The antisymmetry assumption on the vector field $\cW$ implies that the forward $\Psi_t$ flow of a point $\G(t)\in \partial_h^- \mW$ terminates in the facing point  $\ovG(t)\in \partial_h^+ \mW$, except for the cases when $t = t_1 ,  t_2$.

Now consider the surface $P_\G\subset \mW$ generated by the $\Psi_t$-flow   of a traversing curve $\G$. The flows of the points 
$q_1' , q_2' \in \G$ are trapped in $\mW$, hence  the surface $P_\G$ is non-compact. 
On the other hand, by the conditions (W1) to (W6) on the function $f$ in the definition of the Wilson flow $\cW$ in Section~\ref{subsec-wilson}, the  vector field $\cW$ is transverse to $\bRt$ away from the boundary cylinders $r=1$ and $r=3$, and the annulus $\{z=0\} = \cA$. It follows that  each connected component of $P_{\G} \cap \bRt$ is a closed   embedded curve in $\bRt$. We next analyze the properties of these curves.

 Observe that the $\Psi_t$-flows of the curve segments $\g$ and $\kappa$ generate infinite propellers $P_\g$ and $P_\kappa$ as in Definition~\ref{def-infpropeller}, and illustrated in Figure~\ref{fig:propeller},   whose trace on $\bRt$ is pictured in Figure~\ref{fig:arcspropeller}(A). We denote the curves in $P_\g\cap\bRt$ by $\g_0(\ell)$, for $\ell> 0$ and unbounded. Recall that these are simple arcs whose endpoints lie on the vertical line $\{r=2\}\cap \bRt$. The lower endpoint $p_1^-(1)$ of $\g_0(1)$ is the first intersection of the forward $\cW$-orbit of $q_1'$ with $\bRt$. Then for $\ell > 1$, the lower endpoint $p_1^-(\ell)$ of $\g_0(\ell)$ is $\whPsi^{\ell -1}(p_1^-(1))$. 
 Likewise, the upper endpoint $p_1^+(1)$ of $\g_0(1)$ is the first
 intersection of the backward $\cW$-orbit of $\oq_1'$ with $\bRt$, and for $\ell > 1$ the upper endpoint of $\g_0(\ell)$ is $\whPsi^{1-\ell}(p_1^+(1))$.

 Analogously, we denote  the curves in $P_\kappa\cap\bRt$ by $\kappa_0(\ell)$, for $\ell>0$ and unbounded. Again,  these are simple arcs whose endpoints lie on the vertical line $\{r=2\}\cap \bRt$. 
The lower endpoint $p_2^-(1)$ of $\kappa_0(1)$ is the first intersection of the forward $\cW$-orbit of $q_2'$ with $\bRt$. Then for $\ell > 1$, the lower endpoint $p_2^-(\ell)$ of $\kappa_0(\ell)$ is $\whPsi^{\ell -1}(p_2^-(1))$. 
 Likewise, the upper endpoint $p_2^+(1)$ of $\kappa_0(1)$ is the first
 intersection of the backward $\cW$-orbit of $\oq_2'$ with $\bRt$, and for $\ell > 1$ the upper endpoint of $\kappa_0(\ell)$ is $\whPsi^{1-\ell}(p_2^+(1))$.

Note  that these two infinite families of curves $\{\g_0(\ell) \mid \ell \geq 1\}$ and $\{\kappa_0(\ell) \mid \ell \geq 1\}$ are interlaced: between two $\g_0$-curves there is a $\kappa_0$-curve, and vice-versa.  To be more precise, the geometry of the Wilson flow maps the     lower endpoint $q_2'$ of the ``upper'' curve $\kappa$ in Figure~\ref{fig:Gammapositive} to the point 
$p_2^-(1) \in \bRt$ , while he     lower endpoint $q_1'$ of the ``lower'' curve $\g$ is mapped to $p_1^-(1) \in \bRt$ which lies \emph{above} $p_2^-(1)$. That is, we have $z(p_1^-(1)) > z(p_2^-(1))$. The return map of the Wilson flow  then preserves this local order, so that we have
\begin{equation}\label{eq-interlaced}
z(p_2^-(1)) <z(p_1^-(1)) < z(p_2^-(2)) <z(p_1^-(2)) < \cdots < z(p_2^-(\ell)) <z(p_1^-(\ell)) < \cdots < -1 \ .
\end{equation}

In order to complete the description of the curves in $P_{\G} \cap
\bRt$ we must consider the $\Psi_t$ flow of the third curve segment $\wtg$, and how the endpoints  of the curves in its intersections with $\bRt$ are attached to the endpoints of the curves $\g_0(\ell)$ and $\kappa_0(\ell)$ for $\ell \geq 1$.
 
The curve $\wtg$ coincides with $\G\cap \{r\leq2\}$ and its endpoints are
$q_1' , q_2' \in \G$. These points   are trapped in $\mW$ for the forward
$\Psi_t$ flow,  hence  the forward flow of $\wtg$  under $\Psi_t$ defines a
non-compact surface $P_{\wtg}$ in the region $\{r\leq2\}\subset \mW$. 

We denote the curves in $P_{\wtg}\cap \bRt$ by
$\wtg_0(\ell)$, for $\ell> 0$ and unbounded. Consider the point
$\G(1)\in \wtg$, whose radius coordinate is less than 2. It follows that the
$\cW$-orbit of $\G(1)$ is finite, thus intersects $\bRt$ in a finite
number of points $n$.  The integer $n$ depends on the choice of the insertion map $\sigma_1^\e$, but we omit this dependence from the notation, as it simplifies the presentation and does not impact the results below.

Moreover,   without loss of
generality, we may assume that the intersection of the $\cW$-orbit of $\G(1)$ with
the annulus $\cA$ is not contained in $\bRt$. The symmetry of the Wilson flow with respect to the annulus $\cA$ implies that the trace of the $\cW$-orbit of $\G(1)$ on $\bRt$ forms a symmetric pattern, where points of intersection are paired, each point    below the center line $\{z=0\} \cap \bRt$ paired with a symmetric copy above this line.    Hence, $n$ is an even number. 

\begin{figure}[!htbp]
\centering
{\includegraphics[width=45mm]{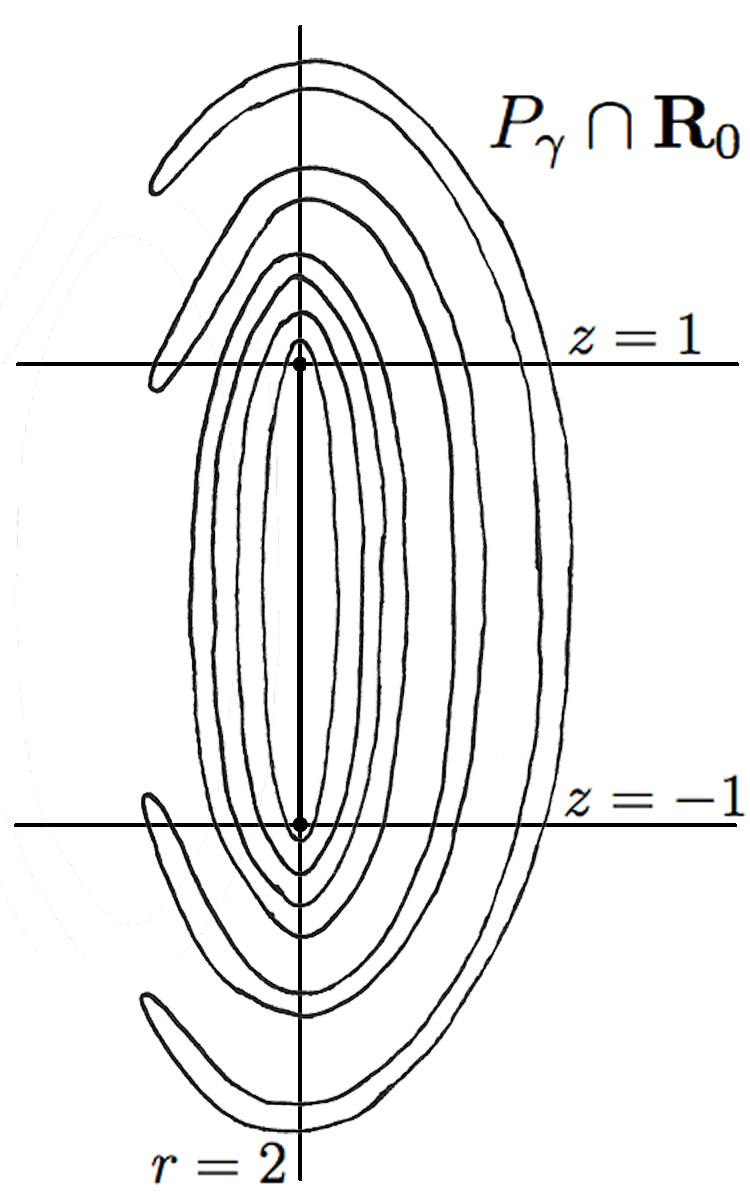}
\caption{\label{fig:tracePG}Trace of the surface $P_\g$ on $\bRt$}}
 \vspace{-6pt}
\end{figure}

The intersection  $P_{\wtg}\cap \bRt$  consists of an infinite
collection of arcs, with endpoints in $\{r=2\} \cap \bRt$. For $1\leq
\ell\leq n/2$ label the   curve containing with endpoints $\{p_1^-(\ell) , p_2^-(\ell) \}$ by $\wtg_0(\ell)$.  Thus, $\wtg_0(\ell)$
 is a ``parabolic curve'' in the region $\{z < 0\} \cap \bRt$, as illustrated in Figure~\ref{fig:tracePG}. 
 
By the anti-symmetry of the flow $\Psi_t$,  for each $1\leq \ell\leq n/2$ there is a corresponding ``parabolic curve'' $\ovtg_0(\ell)$ in the region $\{z > 0\} \cap \bRt$ with endpoints $\{p_1^+(\ell) , p_2^+(\ell) \}$. Let $\G_0(\ell) \subset \bRt$ denote the closed curve 
 obtained by joining the endpoints of the curve $\wtg_0(\ell)$ with the endpoints of the curve $\ovtg_0(\ell)$ via the curves 
 $\{\g_0(\ell),  \kappa_0(\ell)\}$.
The curves $\G_0(\ell)$ for $1 \leq \ell \leq n/2$ are illustrated in Figure~\ref{fig:tracePG} as the closed curves that \emph{do not contain} the vertical arc $\cR \cap \bRt$ in their interiors.

The trace of $P_{\wtg}$ on $\bRt$  contains also an infinite
number of connected arcs that do not contain points of the orbit of $\G(1)$,
and each   of these arcs is again symmetric with respect to $\{z=0\}$. 
For each point $p_1^-(\ell)$ with $\ell > n/2$ we obtain an arc in $\bRt$ with  $p_1^-(\ell)$ as the lower endpoint, and  $p_1^+(\ell)$ as the upper endpoint. When joined  with the arc $\g_0(\ell)$ having the same endpoints, we obtain a closed curve 
$\G_0^{\g}(\ell)$ for $\ell > n/2$, 
that  contains the vertical arc $\cR \cap \bRt$ in its interior. 

Similarly,   for each point $p_2^-(\ell)$ with $\ell > n/2$ we obtain an arc in $\bRt$ with  $p_2^-(\ell)$ as the lower endpoint, and  $p_2^+(\ell)$ as the upper endpoint. When joined  with the arc $\kappa_0(\ell)$ having the same endpoints, we obtain a closed curve $\G_0^{\kappa}(\ell)$ for $\ell > n/2$,  that  also contains the vertical arc $\cR \cap \bRt$ in its interior. 
Moreover, the closed curves $\G_0^{\g}(\ell)$  and $\G_0^{\kappa}(\ell)$  are nested for $\ell > n/2$.
These curves are illustrated in Figure~\ref{fig:tracePG} as the infinite sequence of nested closed curves that   contain the vertical arc $\cR \cap \bRt$ in their interiors.
 
 The curves $\G_0(\ell) \subset \bRt$ for $1 \leq \ell \leq n/2$, and
 $\G_0^{\g}(\ell), \G_0^{\kappa}(\ell) \subset   \bRt$ for $\ell >
 n/2$, can be visualized in terms of a modification of the
 illustration of a standard propeller in
 Figure~\ref{fig:propeller}. Since the curve $\Gamma$ crosses the
 circle $\{r=2\} \cap \partial_h^-\mW$, the flow of the curves
 $\g\cup\tilde{\g}$ and $\tilde{\g}\cup\kappa$ each develops a
 singularity along this circle, which results in an infinite cylinder
 attached to the surface pictured in Figure~\ref{fig:propeller}. This
 infinite cylinder wraps around and is asymptotic to the Reeb cylinder
 $\cR$ in $\mW$. Its intersections with $\bRt$ yields the nested
 family  $\G_0^{\g}(\ell)$ for $\ell > n/2$. A similar statement holds
 for the flow of the curve $\kappa$, yielding the nested curves
 $\G_0^{\kappa}(\ell)$ for $\ell > n/2$. The curves $\G_0(\ell)$ for
 $1 \leq \ell \leq n/2$, are obtained from the flow of the point $\G(1)$, which now lies in the region $r < 2$.

\subsection{The transverse hypotheses}   \label{subsec-transverse}
 
We formulate the conditions on a Kuperberg plug  $\mK_\e$ for $\e > 0$ which will be assumed in the following  sections of the text, where horseshoe dynamics are exhibited  for  the action of the    $\psg$   $\whcGe^*$  on $\bRt$. 

Assume that  the flow $\Psi_t$ on the Wilson Plug satisfies Hypothesis~\ref{hyp-genericW}, and that the construction of the flow $\Phi_t^\e$ on $\mK_\e$  satisfies the conditions (K1) to (K8) in Section~\ref{subsec-kuperberg}, and Hypotheses~\ref{hyp-monotone} and  \ref{hyp-SRI}.

\begin{figure}[!htbp]
\centering
{\includegraphics[width=60mm]{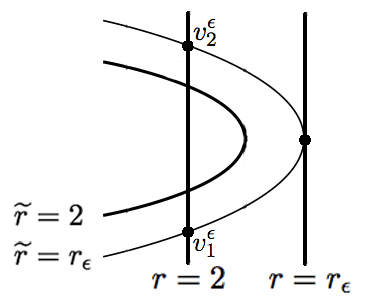}
\caption{\label{fig:Gprime}The points $v_1^\e$ and $v_2^\e$ in $\bRt$}}
 \vspace{-6pt}
\end{figure}

Observe
that for a plug $\mK_\e$ defined by the embeddings $\sigma_i^\e$
satisfying the Parametrized Radius Inequality (K8), the curve
$\theta' \mapsto \sigma_1^\e(2,\theta',-2)$ in $\cL_1^{\e-}$
contains two points with $r=2$. In other words,  the image under
$\sigma_i^\e$ of the circle of radius $r =  2$ intersects twice the
cylinder $\cC\subset \mW$ of points with radius  $r= 2$. Thus the
curve $\wtr=2$ is $\bRt$ intersects the vertical line $r=2$ in
two points, as in
Figure~\ref{fig:Gprime}. 

   Hypothesis~\ref{hyp-monotone} and \ref{hyp-SRI} imply that   the parabolic
   curve of points in $\bRt$ whose forward orbit hits $E^{\e}_1$ in  points
   of radius $r_\e$ is tangent to the vertical line of radius equal to
   $r_\e$. That is, the curve $\wtr=r_\e$ is tangent to the
   vertical line $r=r_\e$. Moreover, by the Parametrized Radius Inequality,  the parabolic  curve of points $\wtr=2$ is tangent to the vertical line of radius equal to $2+\e$, and by Lemma~\ref{lem-r0}, we have that $r_\e>2+\e$, as illustrated in Figure~\ref{fig:Gprime}.

Next, assume that $\e > 0$ is sufficiently small so that  $2+\e < r_\e < 2 + \e_0/2$. This implies that the parabolic curves $\wtr=2$ and $\wtr=r_\e$ in  Figure~\ref{fig:Gprime} intersect the vertical line $r=2$
transversally.  Recall that the orbits of the points $\xi\in \bRt$ such that $\wtr(\xi)=r_\e$
hit $E^{\e}_1$ in points of radius $r_\e$.

 Label the points of intersection of the curve $\wtr=r_\e$ with the
 vertical line $r=2$ by
\begin{equation}\label{eq-boundarypoints}
\{r=2\} \cap \{\wtr=r_\e\} \cap \bRt = \{ v_1^\e, v_2^\e\} \ ,
\end{equation}
 where $z(v_1^\e) < z(v_2^\e) < 0$ by the Hypothesis (K4). We will assume in addition that  $z(v_1^\e) > -1$, so that the curve $\wtr=r_\e$ has a vertical offset.   Note that for $\e=0$, the Radius Inequality implies that $r_\e = 2$ and $v_1^\e = v_2^\e$ is the unique point $\omega_1 \in \bRt$ defined by \eqref{eq-omegas}, so that $z(v_1^\e)  +1 =0$.

 Observe that Hypothesis~\ref{hyp-monotone} implies that if $r(\xi)\geq r_\e$ then $\wtr(\xi)\geq r_\e$,   while condition (K8) implies that if  $r(\xi)<r_\e$ then $\wtr(\xi)\geq r(\xi)-\e$.
 Moreover,  Hypothesis~\ref{hyp-SRI} implies that the curves $\{\wtr=cst\}$ are parabolic, that is the $r$-coordinate depends quadratically on the $z$-coordinate.

Finally, the domain $U_{\phi_1^+}$ of the map $\phi_1^+$ contains a
neighborhood of the set $\{\wtr=2\}$, and thus  we can assume that  for $\e > 0 $   sufficiently small  we have that 
\begin{equation}\label{eq-funddomain}
V_0 = \{x\in \bRt \,|\, r(x)\geq 2, \, \wtr(x)\leq r_\e\} \subset  U_{\phi_1^+} ~.
\end{equation}
That is, the set bounded by the vertical line of radius 2 and the parabolic curve of points whose forward orbit hit $E^{\e}_1$ in  points of radius $r_\e$, is contained in $U_{\phi_1^+}$.

\subsection{Invariant Cantor sets}   \label{subsec-horseshoemap}

 We   assume that    the conditions of Section~\ref{subsec-transverse} are satisfied. 
 
 Let $\G'\subset \cL_1^{\e-}$ be the intersection of the cylinder $\{r=2\}\in \mW$ and $\cL_1^{\e-}$, and set $\G=(\sigma_1^\e)^{-1}(\G')\subset L_1^-$. Observe that $\G$ is a traversing curve with parabolic shape, with its endpoints  in the circle $\{r=3\}$, and admits a parametrization as   in Section~\ref{subsec-props}. In particular, we can divide it in three parts,  $\G=\g\cup \wtg\cup \kappa$ with $\G\cap \{r\geq 2\}=\g\cup\kappa$, as in Figure~\ref{fig:Gammapositive}. 
 
\begin{figure}[!htbp]
\centering
{\includegraphics[width=66mm]{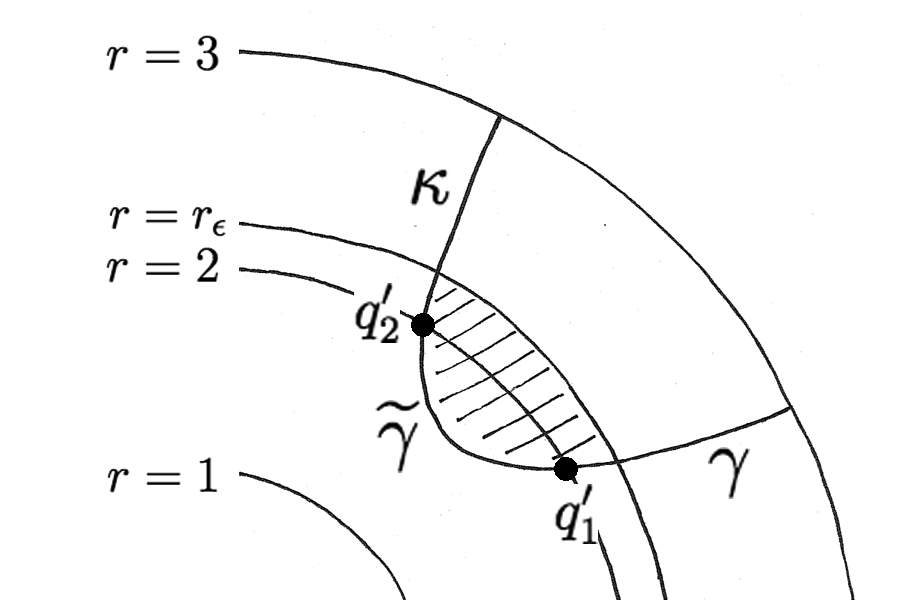}
\caption{\label{fig:UV} The shaded region $V$ in $L_1^-$}}
 \vspace{-6pt}
\end{figure}

  Let $V \subset L_1^-$ be the compact subset bounded by $\G$ and the circle $\{r=r_\e\}$, as illustrated in   Figure~\ref{fig:UV}.
Let $V' = \tau(V) \subset E^{\e}_1 \subset \mK$ be its image in the entry region $E^{\e}_1$. 
Then the set $V_0$ defined in \eqref{eq-funddomain} is identified with
the set of points in $\bRt$ for which the first transition point of
their $\Phi_t^\e$-flow lies in $V'\subset E^{\e}_1$.

Recall that the \emph{trace} of the Reeb cylinder in $\bRt$  is the   line segment 
 $\cR \cap \bRt = \{(r,z) \mid r=2, -1\leq z\leq 1\}$.
 
 Let  $P_\G \subset \mW$ be  the propeller surface generated by the
 $\Psi_t$ flow   of the traversing curve $\G$.
  Let $n > 0$ denote  the number of intersections of the flow of the point $\G(1)$ with the surface $\bRt$, which we assume is  an even integer as discussed in   Section~\ref{subsec-props}. 
Then the trace  $P_\G \cap \bRt$ consists of an infinite collection of closed curves, what we call the $\G_0$-curves. Recall that these have two types:    the closed curves  that \emph{do not contain} the trace of the Reeb cylinder in their interior, and are denoted by     $\G_0(\ell) \subset \bRt$ for for $0 < \ell \leq  n/2$. 
We  also have the two infinite families of closed curves  
 $\{ \G_0^{\g}(\ell) \mid \ell > n/2\}$ and $\{ \G_0^{\kappa}(\ell) \mid \ell > n/2\}$     which contain the trace of the Reeb cylinder in their interior, as  illustrated in  Figure~\ref{fig:tracePG}.  
 For $\ell>n/2$, set $\G_0(\ell)=\G_0^\g(\ell)\cup \G_0^\kappa(\ell)$.

By the symmetry of the Wilson flow, each $\G_0$ curve is symmetric with respect to $\{z=0\}$, so that the observations for the  forward $\Psi_t$ flow of curves in $\partial_h^- \mW$ also apply to the reverse flow of curves in    $\partial_h^+ \mW$.

The traversing curve  $\G=\g\cup \wtg\cup \kappa$  intersects the circle $\{r=2\} \cap L_1^-$ twice, in the points $\{q_1' , q_2'\}$, 
where $q_1'$ is the inner endpoint of $\g$, and $q_2'$ is the
inner endpoint of $\kappa$. The $\cW$-orbits of these points
intersect the vertical segment $\{r=2, \, z<-1\}$ in an increasing
sequence of   interlaced points, as in \eqref{eq-interlaced}, which
limit to $p_1^-$ as $\ell \to \infty$. For each $\ell>0$, the lower intersection point  $p_2^-(\ell)$ is the lower endpoint of the curve $\kappa_0(\ell)$, and the upper intersection point $p_1^-(\ell)$ is the lower endpoint of the curve $\g_0(\ell)$. 

By Hypothesis~\ref{hyp-SRI}, the curves $\g_0(\ell)$ and
$\kappa_0(\ell)$ have parabolic shape near their intersection with the
vertical line $r=2$; that is, the $z$-coordinate depends quadratically on the $r$-coordinate. 
As  the points in the intersection $\G_0(\ell)\cap\{r=2\}$ tend to $p_1^-$, the curves $\g_0(\ell)$ and $\kappa_0(\ell)$ accumulate on the trace of the Reeb cylinder,    $\cR \cap \bRt$,  which is a vertical line segment, 
and these curves become increasingly vertical as they approach $\cR \cap \bRt$. 
 
We next require a technical result, that states that   the    curves $\g_0(\ell)$ and $\kappa_0(\ell)$     are in ``general position'' with respect to the curves $\wtr = 2$ and $\wtr = r_\e$.
Recall that we      assume      the conditions of Section~\ref{subsec-transverse} are satisfied by the given map $\sigma_1^\e$, and in particular, the vertical offset $z(v_1^\e)>-1$ for the points defined by \eqref{eq-boundarypoints}.

\begin{lemma} \label{lem-perturbation}
For $\e > 0$ sufficiently small, there exists $\ell(\e) > 0$ such that for   $\ell \geq \ell(\e)$,  the curves $\g_0(\ell)$ and $\kappa_0(\ell)$
intersect the curves  $\wtr = 2$ and $\wtr = r_\e$ in four points, where the intersections are transverse.  
\end{lemma}
\proof
The assumption that the vertical coordinate $z(v_1^\e) > -1$ implies that for $\ell$ sufficiently large, the curves $\g_0(\ell)$ and $\kappa_0(\ell)$  intersect the curve in $\bRt$ defined by $\wtr=2$ transversely. Let $\ell(\e)$  be the first index for which this holds. Then by the nested properties  of the curves   $\g_0(\ell)$ and $\kappa_0(\ell)$, the transversality property     holds for all $\ell \geq \ell(\e)$ as well.
\endproof
 
\begin{figure}[!htbp]
\centering
{\includegraphics[width=80mm]{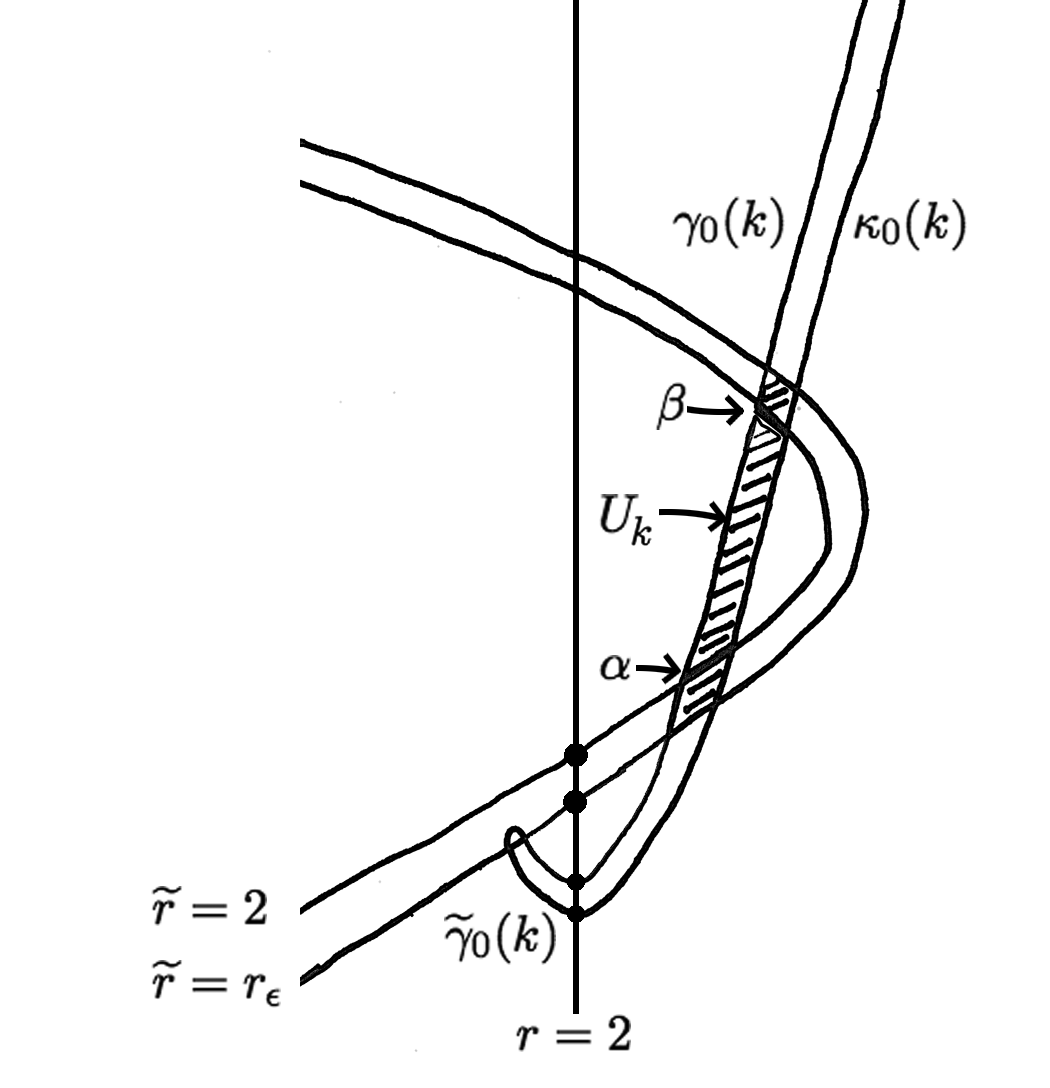}
\caption{\label{fig:positive1}The region $U_k$ and the curves $\alpha$ and $\beta$ in $\bRt$}}
 \vspace{-6pt}
\end{figure}

 The conclusion of Lemma~\ref{lem-perturbation} is   illustrated in Figure~\ref{fig:positive1}. Note that the constant $\ell(\e)$ depends also on the choice of the embedding $\sigma_1^\e$ but for simplicity of notation we omit this dependence from the notation.  
 
\begin{prop}\label{prop-iteration0}
Let $\e > 0$, and let  $V_0$ and  $\ell(\e)>0$  be defined as above.  Fix $k \geq \ell(\e)$, and define $\vp_k=\whPsi^k\circ\phi_1^+ \in \whcGe^*$.  Then for   $V_0$ as defined in \eqref{eq-funddomain}, the following results hold :
\begin{enumerate}
\item  $U_k=\vp_k(V_0)\cap V_0\neq \emptyset$;
\item $U_k$ intersects the curve $\{\wtr=2\}$ along two disjoint arcs $\alpha$ and $\beta$.
\end{enumerate}
\end{prop}

\proof The image $\phi_1^+(V_0)\subset \bRt$ is   the set of first intercept points in $\bRt$   for the $\Psi_t$-flow of the region $V  \subset L_1^-$, where 
 $V$ is illustrated  in Figure~\ref{fig:UV}. 
Then the trace on $\bRt\subset \mW$ of the $\cW$-orbits of
$V$ is the union of the sets  $V_{\ell}=\whPsi^{\ell} \circ \phi_1^+(V_0)\subset \bRt$ for $\ell>0$.
The region $V_\ell$ is bounded by a segment in the vertical line     $\{r= r_\e\} \cap \bRt$, and by the curve $\G_0(\ell)$, where we identify the rectangles $\bRt$ in $\mW$ and $\mK_\e$.

For  $k \geq \ell(\e)$, $\G_0(k)$ intersects the curves $\{\wtr=2\}$ and $\{\wtr=r_\e\}$ 
transversely,  and each intersection consists of four points. Then $V_k$ intersects the region bounded by the curve $\{\wtr=r_\e\}$, implying that $U_k\neq \emptyset$ and that $U_k\cap \{\wtr=2\}$ has two connected components, which are labeled  $\alpha$ for the lower one, and   $\beta$ for the upper one, as in Figure~\ref{fig:positive1}. Then $U_k$ is bounded by $\g_0(k)$, $\kappa_0(k)$, the vertical line $\{r=2\}$ and the curve $\{\wtr=r_\e\}$.

Observe that $U_k$ is bounded by parts of the curves $\g_0(k)$ and $\kappa_0(k)$. In Figure~\ref{fig:positive1} we are assuming that $k\leq n/2$, and the case $k>n/2$ is analogous since the shape of $\g_0(k)$ and $\kappa_0(k)$ near $V_0$ is analogous for the two cases of $k$, as depicted in Figure~\ref{fig:tracePG}.
 \endproof

\begin{remark}\label{rmk-curves}
A priori, the boundary of $U_k$ is contained in the union of
$\g_0(k)$, $\kappa_0(k)$, the curve $\{\wtr=r_\e\}$ and the vertical
line $\{r=2\}$. The assumption that $z(v_1^\e) > -1$ implies that this region is disjoint from the vertical line $\{r=2\}$. 
\end{remark}

\subsection{The shift map}   \label{subsec-labelhorseshoemap}

Let $\e > 0$, and let  $V_0$ and  $\ell(\e)>0$  be defined as in Section~\ref{subsec-horseshoemap}.  Fix $k \geq \ell(\e)$   and define $\vp_k=\whPsi^k\circ\phi_1^+$ and  $U_k=\vp_k(V_0)\cap V_0$, as in Proposition~\ref{prop-iteration0}.   
  We   describe in detail the set $U_k \cap \vp_k(U_k)$ and show that it is composed of two connected components, and  give some of the details of the description of the set $U_k\cap \vp_k(U_k\cap\vp_k(U_k))$. The recursive construction will then give us a collection  of $2^n$  disjoint compact regions in the image of $\vp_k^n$, and their infinite intersection defines a   Cantor set invariant under the action of $\vp_k$, for which the restricted action is conjugate to the full shift.

\begin{figure}[!htbp]
\centering
{\includegraphics[width=62mm]{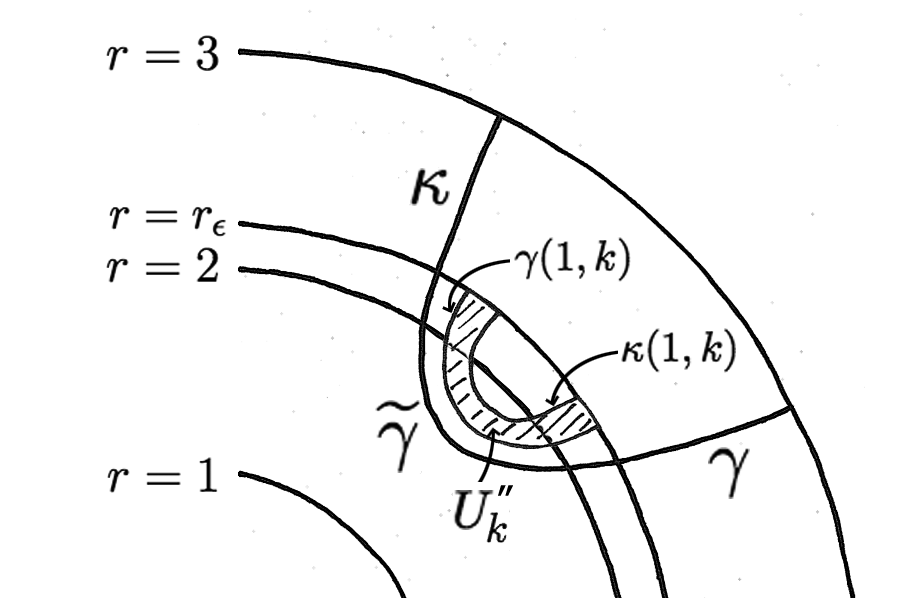}
\caption{\label{fig:U}The region $U_k''$ and parts of the curves
  $\kappa(1,k)$ and $\g(1,k)$ in $L_1^-$}}
 \vspace{-10pt}
\end{figure}

Observe that $U_k\subset V_0\subset U_{\phi_1^+}$. The forward $\Phi_t^\e$-flow 
of $U_k$ intersects the secondary entry region $E^{\e}_1$ in the subset $U_k'\subset
E^{\e}_1$. Then $U_k'$    is bounded by the curve $\{r=r_\e\}$, 
and the curves $\g'(1,k)$ and $\kappa'(1,k)$, obtained by flowing
forward to $E^{\e}_1$ the curves $\g_0(k)$ and $\kappa_0(k)$,
respectively.  Set $U_k''=\tau^{-1}(U_k')\subset L_1^-$, bounded by the
circle $\{r=r_\e\}$ and the curves $\g(1,k)=\tau^{-1}(\g'(1,k))$ and
$\kappa(1,k)=\tau^{-1}(\kappa'(1,k))$. By construction, $\g(1,k)$
intersects twice the circle $\{r=2\}$ and  $\kappa(1,k)$ also
intersects twice $\{r=2\}$, thus both $\g(1,k)$ and $\kappa(1,k)$ are
traversing curves, as in Figure~\ref{fig:U}.

The surfaces $P_{\g(1,k)}$ and $P_{\kappa(1,k)}$ in $\mW$, generated
by the curves $\g(1,k)$ and $\kappa(1,k)$, intersect the rectangle
$\bRt$ in a similar manner as $P_\G$, as in
Figure~\ref{fig:tracePG}. Let $\g_0(1,k;\ell)$ and
$\kappa_0(1,k;\ell)$ for $\ell\geq 1$ and unbounded be the curves in
the trace of the surfaces $P_{\g(1,k)}$ and $P_{\kappa(1,k)}$ on
$\bRt$. The shape of these curves is analogous to the one of
  the curves $\G_0(\ell)$ in the trace of $P_\G$. Observe that for each $\ell>0$, the curves $\g_0(1,k;\ell)$ and $\kappa_0(1,k;\ell)$ are contained in the region bounded by $\G_0(\ell)$ and as $\ell\to \infty$ these curves accumulate on $\cR\cap \bRt$. 

 Since $U_k'' \subset L_1^-$  is bounded by $\g(1,k)$, $\kappa(1,k)$ and the circle $\{r=r_\e\}$, the region $\whPsi^{\ell}\circ\phi_1^+(U_k)$ is bounded by the curves $\g_0(1,k;\ell)$, $\kappa_0(1,k;\ell)$ and the vertical line $\{r=r_\e\}$.

\begin{figure}[!htbp]
\centering
{\includegraphics[width=76mm]{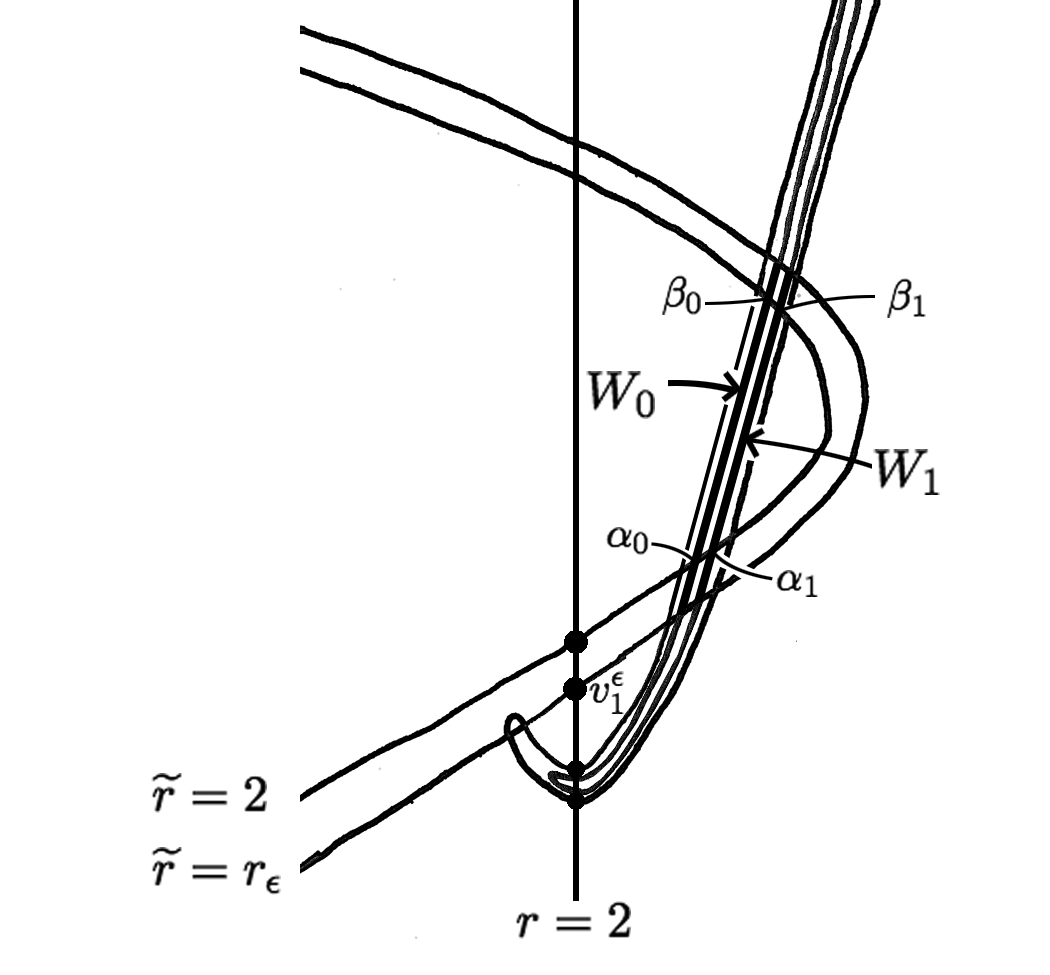}
\caption{\label{fig:positive2}The regions $W_0$ and $W_1$ in $\bRt$}}
 \vspace{-10pt}
\end{figure}

Consider now the case $\ell=k\geq \ell(\e)$.  Then
$\vp_k(U_k)=\whPsi^k\circ\phi_1^+(U_k)$   intersects the set $\{r<
2\}$ along one connected component that is U-shaped and bounded by
$\g_0(1,k;k)\cap \{r<2\}$ and $\kappa_0(1,k;k)\cap \{r<2\}$. Also,
$\vp_k(U_k)$  has two connected components in  $\{r\geq 2\}$, each
corresponding to one of the connected components in $\g_0(1,k;k)\cap
\{r\geq 2\}$. Hence $U_k \cap \vp_k(U_k)$ has two connected
components. Let $W_0$ be the component on the left and $W_1$ the
component on the right as in Figure~\ref{fig:positive2}. The sets
$W_0$ and $W_1$ depend on $k$, but we omit this dependence in the notation.

\begin{lemma}\label{lem-iteration1}
For $i=0,1$, the sets $W_i$ are non-empty and intersect the curve $\{\wtr=2\}$ along two arcs $\alpha_i\subset \alpha$ and $\beta_i\subset \beta$. 
\end{lemma}

\proof
The discussion above implies that the sets $W_i$ are non-empty. Also, 
by construction  $W_i\cap\alpha$ and $W_i\cap\beta$ are non-empty for $i=0,1$, set $\alpha_i=W_i\cap \alpha$ and $\beta_i=W_i\cap \beta$, for $i=0,1$.
\endproof

We can now iterate the construction. For $i_1=0,1$, consider the set $\vp_k(W_{i_1})=\whPsi^k\circ\phi_1^+(W_{i_1})$. Since the curves bounding $W_{i_1}$ cross the curve $\{\wtr=2\}$, the set  $\vp_k(W_{i_1})\subset \vp_k(U_k)$ intersects the side $\{r<2\}$ in a U-shaped set and $\vp_k(W_{i_1})\cap \{r\geq 2\}$ has two connected components one inside $W_0$ and the other inside $W_1$. Thus we can define the four sets
$$W_{i_1,i_2}\,=\,\vp_k(W_{i_1})\cap W_{i_2}\,\subset\,W_{i_2} \qquad \text{for}\qquad i_1,i_2=0,1.$$
In complete analogy with Lemma~\ref{lem-iteration1}, we have that each
$W_{i_1,i_2}$ is non-empty and intersects $\{\wtr=2\}$ along two arcs $\alpha_{i_1,i_2}\subset \alpha_{i_2}$ and $\beta_{i_1,i_2}\subset \beta_{i_2}$. Observe that $\vp_k^{-1}(W_{i_1,i_2})\subset W_{i_1}$.

Iterating this construction, we conclude that for any $n>0$ and for
any sequence $I=\{i_1,i_2, \ldots, i_n\}\in \{0,1\}^n$ the set $W_I=W_{i_1,i_2, \ldots, i_n}$ satisfies that:
\begin{itemize}
\item $W_I$ is non-empty;
\item $W_I$ intersects the curve $\{\wtr=2\}$ along two arcs $\alpha_I$ and $\beta_I$.
\end{itemize}

We make two observations. For a sequence $I\in \{0,1\}^n$, $I=\{i_1,i_2,\ldots, i_n\}$ and a point $\xi\in W_I$, we have 
\begin{eqnarray}
\vp_k^{-1}(\xi) & \in & W_{i_1,i_2,\ldots, i_{n-1}} \label{eq-shift1}\\
\vp_k^{-2}(\xi)=\vp_k^{-1}\circ \vp_k^{-1}(\xi) & \in & W_{i_1,i_2,\ldots, i_{n-2}},  \label{eq-shift2}
\end{eqnarray}
and so forth, so that by induction we have $\vp_k^{-(n-1)}(\xi)\in W_{i_1}$. Also by construction,
$$W_I\subset W_{i_2,i_3,\ldots,i_n}\subset W_{i_3,i_4,\ldots,i_n}\subset \cdots \subset W_{i_n}.$$

We can now define 
  \begin{equation}\label{eq-cantorset}
W(\vp_k) = \bigcap_{n\geq 1} \, \bigcap_{I=\{i_1,i_2,\ldots, i_n\}} ~ W_I
\end{equation}
which is a Cantor set that is invariant under the map $\vp_k$. 

Observe that each $\xi \in W(\vp_k)$ is uniquely defined by an
infinite string  $\ds I_{\xi} = \{ i_1,i_2,\ldots, i_n, \ldots\}\in \{0,1\}^\mN$, which we call the shift coordinates on $W(\vp_k)$.  The identities \eqref{eq-shift1}, \eqref{eq-shift2} and their generalization imply that under this identification, the map $\vp_k$ acts on the points of $W(\vp_k)$ via the right shift on the shift coordinates. We then have the standard observation about such maps:
\begin{lemma}\label{lem-periodic}
The periodic orbits for the restricted action $\vp_k \colon W(\vp_k) \to W(\vp_k)$ are dense in $W(\vp_k)$.
\end{lemma}

Finally, observe that   for $k' > k \geq \ell(\e)$,  the regions bounded by the curves $\G_0(k')$ and $\G_0(k)$ are disjoint, hence we have that   $W(\vp_{k'}) \cap W(\vp_k) = \emptyset$.

 Combining the above results and observations, we have shown:

\begin{thm}\label{thm-horseshoe}
For $0<\e<\e_0$ sufficiently small,   let  $V_0$ and  $\ell(\e)>0$  be defined as above.  Then for each $k \geq \ell(\e)$,   define $\vp_k=\whPsi^k\circ\phi_1^+ \in \whcGe^*$, and let $W(\vp_k)$ be defined by \eqref{eq-cantorset}.   
Then $W(\vp_k)$ is a Cantor set which is invariant under the action of $\vp_k$, and the action admits a dynamical coding as a full-shift space.  Thus, the action of the $\psg$ $\whcGe^*$ on $\bRt$ contains an infinite number of disjoint horseshoe dynamical systems, each with a dense set of periodic orbits. 
\end{thm}

\section{Topological entropy for $\e > 0$}\label{sec-entropy}

The pseudogroup $ \whcGe$ acting on $\bRt$ was introduced in Definition~\ref{def-pseudogroup1}, with generators obtained from the return maps $\whPsi$ and $\whPhie$.
Then in Section~\ref{sec-horseshoe},       
for $\e > 0$   sufficiently small and  after  imposing a sequence of assumptions on the geometry of the construction of the flow $\Phi_t^\e$, it was shown that the action of $\whcGe$ contains disjoint families of invariant Cantor sets. In this section, we impose one further condition on the construction of the flow $\Phi_t^\e$ which suffices to imply that each of these ``horseshoe dynamical systems'' is realized by the  flow $\Phi_t^\e$. It follows immediately from this that the flow  $\Phi_t^\e$ has positive topological entropy and has infinite families of periodic orbits.    
   
Recall that the map   $\vp_k=\whPsi^k\circ\phi_1^+$ was introduced in Proposition~\ref{prop-iteration0}, and realizes a part of the dynamical properties of the action of the $\psg$ $\whcGe^*$ on $\bRt$. However, the map $\vp_k$ 
 need not be realized by the return map of the flow $\Phi_t^\e$. Recall that the   map $\phi_1^ + = \whPhie | U_{\phi_1^+}$, where $\phi_1^ +$ is defined at  $\xi \in U_{\phi_1^+}$ with $\eta = \phi_1^+(\xi)$, if  there  is a $\cK_\e$-arc $[\xi, \eta]_{\cK_\e}$
  which contains a single transition point $x \in E^{\e}_1$. On the other hand, we have that $\zeta = \whPsi^k(\eta)$ if there exists a $\cW$-arc 
$[\eta, \zeta]_{\cW}$ in $\mW$, which is independent of the return map $\whPhie$.  

The strategy in this section is to apply Corollary~\ref{cor-positivemain} to conclude that there is also  a  $\cK_\e$-arc $[\eta, \zeta]_{\cK_\e}$ in $\mK_\e$ between $\eta$ and  $\zeta$. When applied to points $\xi \in W(\vp_k)$, this will imply that the map  
$\vp_k \colon W(\vp_k) \to W(\vp_k)$   represents a subsystem of the return map $\whPhie$, and hence gives information about the dynamics of the flow $\Phi_t^\e$.     The difficulty is that Corollary~\ref{cor-positivemain} assumes that $\eta \in U_{r_\e}$ where   $U_{r_\e}$ is defined by   \eqref{eq-domainre}. 
We first obtain    conditions on the   map $\sigma_1^\e$ which will imply that this requirement is satisfied for $\xi \in W(\vp_k)$, so that we can then use Theorem~\ref{thm-horseshoe} to obtain a proof of Theorem~\ref{thm-main2}.

 \subsection{Wilson dynamics}\label{subsec-wilsonsteps}
The orbits of the map $\whPsi$ on $\bRt$ are simple to describe, in that for any point $\xi \in \bRt$ with $r(\xi) \ne 2$, the orbit is finite in both forward and backward directions. However,   as the value of $r(\xi)$ tends to $r=2$, the lengths of these finite orbits increase, as   condition \eqref{eq-generic1} implies that the vertical distance   between the iterations $\whPsi^{\ell}(\xi)$ becomes arbitrarily small for $r(\whPsi^{\ell}(\xi))$ near to $2$ and its vertical coordinate $z(\whPsi^{\ell}(\xi))$   near to either $z= \pm 1$.  We give an approximation of the  distance between points $\whPsi^{\ell}(\xi)$ and $\whPsi^{\ell+1}(\xi)$, following the same approach as in  Chapter 17 of \cite{HR2016}. These estimates are then used to impose restrictions on the insertion maps $\sigma_1^\e$ so that Corollary~\ref{cor-positivemain} can be applied.

   Recall the functions $f$ and $g$ were chosen in Section~\ref{subsec-wilson}, which are constant in the coordinate $\theta$, with
\begin{equation}\label{eq-wilson2}
\cW =g(r, \theta, z)  \frac{\partial}{\partial  z} + f(r, \theta, z)  \frac{\partial}{\partial  \theta} ~ .
\end{equation}
  Hypothesis~\ref{hyp-genericW}  and condition \eqref{eq-generic1} imply there exists constants $A_g,B_g,C_g$ such that  the quadratic form  $Q_g(u,v) = A_g \cdot u^2 + 2B_g \cdot uv + C_g \cdot v^2$ defined by the Hessian of $g$ at $\omega_1$ is positive definite. Set 
  $$Q_0(r,z) = \left( d_{\bRt}((r -2 ) , (z+1))\right)^2 = (r-2)^2 + (z+1)^2$$
 then  it follows that  there exists   $D_g > 0$ such that 
  \begin{equation}\label{eq-quadraticest1}
| g(r, \theta, z)  - Q_g(r-2, z+1) | ~ \leq ~   D_g \cdot (|r-2|^3 + |z+1|^3) \quad {\rm for} ~   Q_0(r,z)  \leq \e_0^2
\end{equation}
where $\e_0$ is the constant defined in \eqref{eq-generic1}. The condition \eqref{eq-quadraticest1} implies that for $(r,z)$ sufficiently close to $(2,-1)$,  the error term on the right-hand-side can be made arbitrarily small relative to the distance squared   from the special point $\omega_1 = (2,-1)$. We also observe that    \eqref{eq-quadraticest1}   implies there exists constants $0 < \lambda_1 \leq \lambda_2$ such that
  \begin{equation}\label{eq-quadraticest2}
\lambda_1 \cdot Q_0(r,z)  \leq  g(r, \theta, z) \leq \lambda_2 \cdot Q_0(r,z)   \quad {\rm for} ~   Q_0(r,z)  \leq \e_0^2 \ . 
\end{equation}
Next, consider the action of the maps $\whPsi^{\ell}$ for $\ell > 0$.
   Let $\xi  \in \bRt$ with $2 \leq r(\xi) \leq 2 + \e_0$ and $-7/4 \leq z(\xi) \leq -1/4$, such that $\whPsi(\xi)$ is defined and   $z(\whPsi(\xi))< 0$.
   Let $T(\xi) > 0$ be defined by $\whPsi(\xi) = \Psi_{T(\xi) }(\xi)$. Then the $z$-coordinate of $\whPsi(\xi)$ is given by 
   \begin{equation}\label{eq-coordinates}
  z(\whPsi(\xi)) -  z(\xi) ~ = ~ \int_0^{T(\xi) } ~ g(\Psi_s(\xi)) ~ ds  ~ \geq ~ 0  \ .
\end{equation}
If $\xi \ne \omega_1$ then $g(\Psi_s(\xi))$ is positive along the orbit segment for $0 \leq s \leq T(\xi)$, hence  $ z(\whPsi(\xi)) -  z(\xi) > 0$. Moreover, combining the estimates \eqref{eq-quadraticest2} and \eqref{eq-coordinates} and the estimate $T(\xi) \geq 4\pi$ on the return time for the flow $\Psi_t$ to $\bRt$ for $r(\xi) \geq 2$, we obtain:
\begin{lemma}\label{lem-stepesitimates}
For $\delta > 0$, suppose that $\xi = (r,z) \in \bRt$ satisfies $\delta < d_{\bRt}((r -2 ) , (z+1)) <\e_0$. Then 
\begin{equation}\label{eq-stepesitimates}
 z(\whPsi(\xi)) -  z(\xi) ~ \geq  ~ 4 \pi \cdot \lambda_1 \delta^2 \ .
\end{equation}
\end{lemma}

\subsection{Admissible deformations}\label{subsec-deformations}

We   assume that  for $0 < \e < \e_0$  the conditions   of Section~\ref{sec-horseshoe} are satisfied, so that  the hypotheses of Proposition~\ref{prop-iteration0} and Theorem~\ref{thm-horseshoe} are satisfied. 

Recall that for  the insertion map $\sigma_1^\e$      the points
$\{v_1^\e , v_2^\e \}\subset \bRt$  are defined by  \eqref{eq-boundarypoints} and satisfy 
\begin{equation}\label{eq-offset0}
r(v_1^\e) = r(v_2^\e) = 2 \quad , \quad -1 < z(v_1^\e) < z(v_2^\e) < 0 \ .
\end{equation}
 Define 
\begin{equation}\label{eq-offset1}
\delta(\sigma_1^\e) = z(\whPsi^{-1}(v_1^\e)) + 1 > 0\ .
\end{equation}

\begin{hyp}\label{hyp-offset}
Assume that $\e > 0$ is sufficiently small so that 
\begin{equation}\label{eq-offset2}
 z(v_2^\e) -  z(v_1^\e) < 4 \pi \cdot \lambda_1 \cdot \delta(\sigma_1^\e)^2 \ .
\end{equation}
\end{hyp}
Hypothesis~\ref{hyp-SRI} implies that  the boundary curve 
$\{\wtr=r_\e\}\subset \bRt$ of the region $V_0$ is parabolic, as illustrated in
Figure~\ref{fig:Gprime}, which implies that  for a fixed   $v_1^\e \in \bRt$ then      
   for $\e >0$ sufficiently small, the condition \eqref{eq-offset2} will be satisfied.

Let  $V_0$ and  $\ell(\e)>0$  be as defined as in Theorem~\ref{thm-horseshoe}.
   Fix $k \geq \ell(\e)$, and define $\vp_k=\whPsi^k\circ\phi_1^+$.   
Set  $U_k=\vp_k(V_0)\cap V_0$ as defined in  Proposition~\ref{prop-iteration0}.
Then let $W(\vp_k)$ be defined by \eqref{eq-cantorset}. By Theorem~\ref{thm-horseshoe} we have that 
$\vp_k \colon W(\vp_k) \to W(\vp_k)$.

\begin{prop}\label{prop-realized}
Assume that Hypothesis~\ref{hyp-offset} holds for  the insertion map $\sigma_1^\e$. Then for each $\xi \in W(\vp_k)$ with $\eta = \vp_k(\xi)$, there exists a $\cK_\e$-arc $[\xi, \eta]_{\cK_\e}$ in $\mK_\e$. Thus, there exists $\ell_{\xi} \geq 1$ such that 
$\eta = \whPhie^{\ell_{\xi}}(\xi)$. 
\end{prop}
\proof
Let $(r,z) \in \bRt$ satisfy $Q_0(r,z) \leq \e_0^2$ and $z \geq z(\whPsi^{-1}(v_1^\e))$. 
Then  $Q_0(r,z) \geq Q_0(2,z) \geq \delta(\sigma_1^\e)^2$, so that by \eqref{eq-stepesitimates} we have 
\begin{equation}
z(\whPsi(r,z)) -z ~ \geq  ~ 4 \pi \cdot \lambda_1 \cdot \delta(\sigma_1^\e)^2 >  z(v_2^\e) -  z(v_1^\e) \ .
\end{equation}
Now let $\eta \in V_0$, for $V_0$ defined by \eqref{eq-funddomain}. Then $z(v_1^\e) \leq z(\eta) \leq z(v_2^\e)$ by the convexity of the region $V_0$. It then follows from \eqref{eq-stepesitimates} applied to $(r,z) = \eta$ that 
\begin{equation}
z(\whPsi^{-1}(\eta)) = z + \left(z(\whPsi^{-1}(\eta)) -z \right)   \leq z - 4 \pi \cdot \lambda_1 \cdot \delta(\sigma_1^\e)^2   <  z(v_2^\e) - \left( z(v_2^\e) -  z(v_1^\e) \right)  = z(v_1^\e) \ .
\end{equation}
That is, the image $\whPsi^{-1}(V_0)$ lies below the line $\{z=
z(v_1^\e)\}$, and hence lies outside the region bounded by the curve $\wtr = r_\e$ whose lower edge lies on this line by the choice of $v_1^\e$ in \eqref{eq-boundarypoints} and the convexity of the region $V_0$.  That is, for all $\eta \in V_0$ we have $\wtr(\whPsi^{-1}(\eta))> r_\e$.

Let  $\xi \in W(\vp_k)$ and set  $\eta = \vp_k(\xi) = \whPsi^k\circ\phi_1^+(\xi)   \in W(\vp_k)  \subset V_0$.
 Then for  the point $\zeta = \whPsi^{-1}(\eta)$ we have $\wtr(\zeta) > r_\e$ by the above estimates.

 Following the $\cK_\e$-orbit of $\zeta$ backwards and applying  Proposition~\ref{prop-positivemain} and Corollary~\ref{cor-positivemain} we obtain that   $\xi$ and $\zeta$ are in the same $\cK_\e$-orbit.

 Now consider the $\cW$-orbit segment  $[\zeta, \eta]_{\cW}$. 
 If $[\zeta, \eta]_{\cW}$ does not contain any transition 
 points, then $\eta$ is in the $\cK_\e$-orbit of $\zeta$ and thus in the $\cK_\e$-orbit 
 of $\xi$. 
 If $[\zeta, \eta]_{\cW}$  contains transition points, let $x_1$ be the first 
 transition point. Then $x_1$ is a secondary entry point in $E^{\e}_1$ with $r(x_1)> r_\e$ and 
 Corollary~\ref{cor-positivemain} implies that its $\cK_\e$-orbit contains the facing point 
 $y_1\in S^{\e}_1$.  The $\cK_\e$-orbit of    $y_1$ continues to the first intersection with 
 $\bRt$ which is the point $\eta$. Hence $\eta$ is in the $\cK_\e$-orbit of $\xi$.
\endproof

We next deduce an important consequence of Proposition~\ref{prop-realized}. Observe that for every $\xi \in W(\vp_k)$ we have $z(\xi) < 0$ and thus the flows $\Psi_t$ and $\Phi_t^\e$ intersect the set $W(\vp_k) \subset \bRt$ transversally. Let $\ell_{\xi} \geq 1$ be the integer defined in Proposition~\ref{prop-realized} so that $\vp_k(\xi) = \whPhie^{\ell_{\xi}}(\xi)$. Then both maps are continuous at $\xi$ and there is an open neighborhood $U_{\xi} \subset V_0$ for which the identity $\vp_k(\xi') = \whPhie^{\ell_{\xi}}(\xi')$ holds for all $\xi' \in U_{\xi}$. The set $W(\vp_k)$ is compact, so there exists a finite collection of such open sets which cover $W(\vp_k)$. Hence, there exists $N_k$ so that for all $\xi \in W(\vp_k)$ there exists $1 \leq \ell_{\xi} \leq N_k$ for which 
$\vp_k(\xi) = \whPhie^{\ell_{\xi}}(\xi)$.  As a consequence, we have:

\begin{cor}\label{cor-realized}
Assume that Hypothesis~\ref{hyp-offset} holds for  the insertion map $\sigma_1^\e$. Then   there exists $T_k > 0$ such that 
for each $\xi \in W(\vp_k)$ there exists $0 < t_{\xi} \leq T_k$ such that $\vp_k(\xi) = \Phi_{t_{\xi}}^\e(\xi)$. 
\end{cor}

\subsection{Topological entropy}\label{subsec-entropy}
  
We define the entropy of a flow $\vp_t$  on a compact metric space $(X, d_X)$ using a variation of the Bowen formulation of   topological entropy  \cite{Bowen1971,Walters1982}.   The definition we adopt    is symmetric in the role of the time variable $t$.   For $T > 0$, define a metric on $X$ by 
\begin{equation}
d_X^T(x,y) = \max \ \left\{ d_X(\vp_t(x),\vp_t(y)) \mid -T \leq t \leq T \right\} \ , \quad {\rm for} ~ x,y \in X \ . 
\end{equation}
 Two points $x,y\in X$ are said to be  \emph{$(\vp_t , T, \delta)$-separated} if  $d_X^T(x,y)>\delta$.
A set $E \subset X$ is \emph{$(\vp_t , T, \delta)$-separated} if  all pairs of distinct 
points in $E$ are $(\vp_t , T, \delta)$-separated. Let $s(\vp_t , T, \delta)$ be the maximal
cardinality of a $(\vp_t , T, \delta)$-separated set in $X$. Then the topological entropy
is defined by
\begin{equation}\label{eq-separated}
h_{top}(\vp_t)= \frac{1}{2} \cdot \lim_{\delta\to  0} \left\{ \limsup_{T\to\infty}\frac{1}{T}\log(s(\vp_t , T, \delta)) \right\} .
\end{equation}
It is a standard fact that for a compact space $X$, the entropy $h_{top}(\vp_t)$ is independent of the choice of the metric $d_X$ on $X$.

Given a $\vp_t$-invariant subset $K \subset X$, we can define the restricted topological entropy $h_{top}(\vp_t, K)$   by the same formula \eqref{eq-separated}, where we now require that  the $(\vp_t , T, \delta)$-separated sets in the definition must be subsets of $K$. It follows immediately that we have the estimate
\begin{equation} \label{eq-Kseparated}
h_{top}(\vp_t) \geq h_{top}(\vp_t, K) \ .
\end{equation}

Let $\e > 0$ be given, and assume that $\Phi_t^\e$ is a Kuperberg flow as constructed above which satisfies the \emph{generic hypotheses} of Section~\ref{sec-generic},   the \emph{geometric hypotheses} of Section~\ref{sec-horseshoe}, and the   Hypothesis~\ref{hyp-offset} above. Then we have:

\begin{thm}\label{thm-entropypositive}
The topological entropy $h_{top}(\Phi_t^\e) > 0$. 
\end{thm}
\proof
Let  $V_0$ and  $\ell(\e)>0$  be defined as in Theorem~\ref{thm-horseshoe}.  Choose $k \geq \ell(\e)$, and  define $\vp_k=\whPsi^k\circ\phi_1^+ \in \whcGe^*$. Then let $W(\vp_k) \subset \bRt$ be the $\vp_k$-invariant Cantor set defined by \eqref{eq-cantorset}. By Proposition~\ref{prop-realized}, the set $W(\vp_k) $ is invariant under the return map $\whPhie$ for the flow $\Phi_t^\e$. Let $\whW(\vp_k)$ denote the flow saturation of $W(\vp_k)$ which is then a compact invariant set for $\Phi_t^\e$. By \eqref{eq-Kseparated} it will suffice to show that $h_{top}(\Phi_t^\e, \whW(\vp_k))  > 0$.

Note that $W(\vp_k) \subset \whW(\vp_k)$ is a transverse  section for the flow $\Phi_t^\e$ restricted to $\whW(\vp_k)$, and by Corollary~\ref{cor-realized} the flow has bounded return times to the section $W(\vp_k)$. 

\begin{lemma}\label{lem-returnmap}
The map $\vp_k \colon W(\vp_k) \to W(\vp_k)$ is the first return map for the flow $\Phi_t^\e$ restricted to $\whW(\vp_k)$.
\end{lemma}
\proof
Note that $\whPhie$ is the first return map of the flow $\Phi_t^\e$ to the section $\bRt$. 
Given $\xi \in W(\vp_k)$, Proposition~\ref{prop-realized} states that
there is some least $\ell_{\xi}$ such that $\vp_k(\xi) =
\whPhie^{\ell_{\xi}}(\xi) \in W(\vp_k)$. Thus, it suffices to show the
$\cK_\e$-orbit segment $[\xi,\vp_k(\xi)]_{\cK_\e}$ is disjoint from
$\whW(\vp_k)\subset U_k$ except at the endpoints. In fact, we prove   that
for any $\xi\in U_k$ such that $\vp_k(\xi)\in U_k$, the $\cK_\e$-orbit
segment $[\xi,\vp_k(\xi)]_{\cK_\e}$ does not intersect $U_k$.

Since $\xi \in W(\vp_k)\subset U_k\subset U_{\phi_1^+}$, we have that
$\whPhie(\xi)=\phi_1^+(\xi)=\eta$ and thus $r(\eta) \leq r_\e$ since
$U_k\subset V_0$ and the points in $V_0$ have $\wtr$-coordinate
smaller or equal to $r_\e$.  
Consider the $\cK_\e$-orbit segment
$[\eta,\whPsi(\eta)]_{\cK_\e}$. 
If $[\eta,\whPsi(\eta)]_{\cK_\e}$ does not contains transition points,
then it does not intersect $\bRt$ in its interior. If not, let $x_1$ be the first transition point that is
a secondary entry point in $E^{\e}_1$. Since $z(\eta)<z(\whPhi^{k-1}(\eta))$, then
by the arguments in the proof of Proposition~\ref{prop-realized}, we
have 
$$
z(\eta)<z(\whPhi^{k-1}(\eta))<z(v_1^\e),
$$
thus $x_1$ is a secondary entry point outside the region $U'\subset
E^{\e}_1$. In particular, $r(x_1)>r_\e$. Let $y_1\in S^{\e}_1$ be the secondary
exit point such that $x_1\equiv y_1$, then by construction $y_1$ is
the last transition point in $[\eta,\whPsi(\eta)]_{\cK_\e}$. By
Proposition~\ref{prop-positivemain} and Hypothesis~\ref{hyp-monotone} all the $\cW$-arcs in the $\cK_\e$-orbit segment
$[\eta,\whPsi(\eta)]_{\cK_\e}$  have radius greater than $r_\e$. Hence any point in the intersection of the interior of the orbit segment
$[\eta,\whPsi(\eta)]_{\cK_\e}$ with $\bRt$ has radius greater than
$r_\e$ and thus does not belongs to $U_k$, in particular does not
belongs to $\whW(\vp_k)$.

Repeating the argument, consider the   $\cK_\e$-orbit segment
$[\whPsi^\ell(\eta),\whPsi^{\ell+1}(\eta)]_{\cK_\e}$ for any $1\leq \ell
\leq k-1$.
If $[\whPsi^{\ell}(\eta),\whPsi^{\ell+1}(\eta)]_{\cK_\e}$ does not contains transition points,
then it does not intersects $\bRt$ in its interior. If not, let $x_1$ be the first transition point that is
a secondary entry point in $E^{\e}_1$. Since $z(\whPsi^{\ell}(\eta))\leq z(\whPhi^{k-1}(\eta))$, then
we have 
$$
z(\whPsi^{\ell}(\eta))\leq z(\whPhi^{k-1}(\eta))<z(v_1^\e),
$$
thus $x_1$ is a secondary entry point outside the region $U'\subset
E^{\e}_1$. In particular, $r(x_1)>r_\e$. Let $y_1\in S^{\e}_1$ be the secondary
exit point such that $x_1\equiv y_1$, then by construction $y_1$ is
the last transition point in $[\whPsi^{\ell}(\eta),\whPsi^{\ell+1}(\eta)]_{\cK_\e}$. Then all the $\cW$-arcs in the $\cK_\e$-orbit segment
$[\whPsi^{\ell}(\eta),\whPsi^{\ell+1}(\eta)]_{\cK_\e}$  have radius
greater than $r_\e$. Hence any point in the intersection of the interior of the orbit segment
$[\whPsi^{\ell}(\eta),\whPsi^{\ell+1}(\eta)]_{\cK_\e}$ with $\bRt$ has radius greater than
$r_\e$ and thus does not belongs to $U_k$, in particular does not
belongs to $\whW(\vp_k)$.
\endproof

We claim that the restricted map 
$\vp_k \colon W(\vp_k) \to W(\vp_k) $ has positive entropy, and it then follows by standard techniques that $h_{top}(\Phi_t^\e, \whW(\vp_k))  > 0$.

In the following, we use Proposition~\ref{prop-iteration0} and the labeling system in Section~\ref{subsec-labelhorseshoemap}. 
In  particular, let 
   $W_0, W_1 \subset \bRt$ be the disjoint sets introduced in    Lemma~\ref{lem-iteration1}.
Let  $\delta>0$ be such that the $d_{\bRt}$-distance   between the sets $W_0$ and $W_1$ is greater than $\delta$. 
 We obtain a lower bound estimate on the maximum number $s(\vp_k , n, \delta)$ of points in an $(\vp_k , n, \delta)$-separated subset of $W(\vp_k)$   as a function of $n$.

 For a sequence $I\in \{0,1\}^n$, recall that $W_I   =  W_{i_1,i_2,\ldots,i_n}$ is a closed and open cylinder set in the Cantor set $W(\vp_k)$. 
 For any pair of distinct sequences $I,J\in \{0,1\}^n$, choose  points 
 \begin{equation}
\xi\in W_I \cap W(\vp_k)  \quad , \quad \eta \in W_J \cap W(\vp_k)   \ .
\end{equation}
Let $1 \leq m\leq n$ be the largest index such that $i_m\neq j_m$.
 Then 
\begin{eqnarray*}
 \vp_k^{-(n-m)}(\xi) & \in & W_{i_1,i_2,\ldots, i_m} \, \subset\, W_{i_m}\\
  \vp_k^{-(n-m)}(\eta) & \in & W_{j_1,j_2,\ldots, j_m} \, \subset\, W_{j_m}.
  \end{eqnarray*}

 Fix  $n \geq 1$ and let $E_n \subset U_k \subset \bRt$ be a collection of $2^n$ points obtained by choosing one point from each of the sets  $W_I \cap W(\vp_k) =  W_{i_1,i_2,\ldots,i_n} \cap W(\vp_k)$ where $I \in \{0,1\}^n$. Then $E_n$ is a $(\vp_k,  n, \delta)$-separated set, and so $s(\vp_k ,   n , \delta) \geq 2^n$.  
 It follows that  $s(\vp_k , n, \delta)  \geq 2^n$, and thus $h_{top}(\vp_k, W(\vp_k) ) \geq \ln(2)> 0$. 
\endproof
  
  We point out another consequence of Theorem~\ref{thm-horseshoe} and Proposition~\ref{prop-realized}. 
  
  \begin{cor}\label{cor-periodic}
 The restriction of  the flow $\Phi_t^\e$   to $\whW(\vp_k)$ has a countably dense set of periodic orbits.
  \end{cor} 
 \proof
 The  restriction of the map $\vp_k$  to $\whW(\vp_k)$ has horseshoe dynamics by Theorem~\ref{thm-horseshoe}, and hence has a dense set of periodic orbits. Each of these orbits for $\vp_k$  is also a periodic orbit for the return map $\whPhie$ by Proposition~\ref{prop-realized}.
 \endproof

  Note that the conclusions of Theorem~\ref{thm-horseshoe} apply to the map $\vp_k$ for any choice of $k \geq \ell(\e)$, and for each such flow the set $W(\vp_k) \subset U_k$ as illustrated in Figure~\ref{fig:positive1}. See also Remark~\ref{rmk-curves}. In particular, the set 
  $W(\vp_k)$ is contained in the region between the   curves $\g_0(k),  \kappa_0(k) \subset  \bRt$. 
  As remarked in Section~\ref{subsec-horseshoemap}, these curves are asymptotic to the trace $\cR \cap \bRt$ of the Reeb cylinder as $k$ increases, so the Cantor sets $W(\vp_k)$ also limit to subsets of $\cR$.

Now assume that the flow $\Phi_t^\e$ also satisfies the offset Hypothesis~\ref{hyp-offset}, 
  and let $\cP_\e$ denote the union of all periodic orbits of  $\Phi_t^\e$ and let $\ocP_{\e}$ denote the closure of $\cP_{\e}$ in $\mK_\e$. It then follows from Corollary~\ref{cor-periodic} that the intersection  $\ocP_\e \cap\cR$ is non-empty.    It would be interesting to be able to describe the closure $\ocP_{\e}$, or at least the trace $\ocP_\e \cap\cR$ of this closure on $\cR$, as some sort of analog of the \emph{bad set} for compact flows as discussed in \cite{EMS1977,EpsteinVogt1978}, though that   seems like a rather difficult task.
  
\subsection{Smooth admissible deformations}\label{subsec-admissible}

Suppose that we are given a standard Kuperberg flow $\Phi_t=\Phi_t^0$ on $\mK$ constructed starting with a Wilson flow    
 $\Psi_t$ on the Wilson Plug $\mW$ which satisfies Hypothesis~\ref{hyp-genericW}, and that the construction of the flow $\Phi_t$ on $\mK$  satisfies the conditions (K1) to (K8) in Section~\ref{subsec-kuperberg}, and the \emph{Strong Radius Inequality} in Hypotheses~\ref{hyp-SRI} for $\e = 0$. 
 We  construct a   family of flows $\Phi_t^\e$  for $0 < \e < \e_0$ sufficiently small, which give a smooth deformation of $\Phi_t$ such that for each $\e > 0$, the flow $\Phi_t^\e$ satisfies the hypotheses of Theorem~\ref{thm-entropypositive}. 
 
 The construction of the flows  $\Phi_t^\e$ uses a two-step process. In the first step, we introduce a small vertical offset to the insertion map  $\sigma_i^\e \colon D_1 \to \mW$. After this, we modify the insertion by increasing the radius of the insertion of the vertex of the parabola as in Figure~\ref{fig:modifiedradius}(C).
 
Recall that $(2, \overline{\theta_1}, -1) \in \cR$ is the special point on the Reeb cylinder defined in Hypotheses~\ref{hyp-SRI}, for the case $\e=0$.
Let $0< \e' < \e_0$ and choose a point $v_1(\e') \in \cR$ which satisfies $z(v_1(\e')) = -1 + \e'$ and 
$\theta(v_1(\e')) = \overline{\theta_1}$. That is, $v_1(\e')$ is a perturbation of the special point $(2, \overline{\theta_1}, -1) $  in the vertical direction along $\cR$. 
 Let $\widetilde{\sigma_1^{\e'}} \colon   D_1 \to \mW$ denote the vertical translate of $\sigma_1$ so that the vertex of the image curve 
 $\ds \widetilde{\sigma_1^{\e'}}\left( \{r=2\} \cap L_1^-\right)$ is the point $v_1(\e')$.

Next, slide the map $\ds \widetilde{\sigma_1^{\e'}}$ along the parabolic curve  $\ds \widetilde{\sigma_1^{\e'}}\left( \{r=2\} \cap L_1^-\right)$ so that the vertex now lies on the cylinder $\{r=2+\e\}$ where $\e$ is sufficiently small so that the offset Hypothesis~\ref{hyp-offset} is satisfied for the given Wilson flow $\Psi_t$ on $\mW$. This yields the embedding $\ds \sigma_1^{\e}  \colon   D_1 \to \mW$ for which $v_1(\e')$ is $\cW$-flow of the lower intersection point defined in \eqref{eq-boundarypoints}.  
Intuitively, the value of $\e$ is proportional to $(\e')^2$ as the graph of $\ds \widetilde{\sigma_1^{\e'}}\left( \{r=2\} \cap L_1^-\right)$ as illustrated in Figure~\ref{fig:modifiedradius}(C) is quadratic.   This yields our family of insertion maps $\sigma_1^\e$.

The smooth family of second insertion maps $\ds \sigma_2^\e \colon   D_2 \to \mW$ can then be defined similarly. However, we now point out that in the proof of Theorem~\ref{thm-entropypositive}, and its preparatory results, the role of the second insertion only arises in one manner. We require that the entry/exit condition for orbits of $\cK_\e$ entering a face of the second insertion be satisfied for points $\xi \in E^{\e}_2$ when $r(\xi) > r_\e$. In particular, this condition will be satisfied if the insertion $\sigma_2^\e = \sigma_2$ is the given map for the construction of the flow $\Phi_t^0$, which satisfies the original radius inequality,  and thus as noted previously, has corresponding value of $r_\e = 2$. Thus, in order to obtain  the claim of Theorem~\ref{thm-entropypositive}, it is not even necessary to modify the upper insertion.  

Thus,  we obtain a smooth family of deformations of $\Phi_t$  where each map $\Phi_t^\e$ satisfies the hypotheses of Theorem~\ref{thm-entropypositive}, and hence has positive entropy with infinitely many periodic orbits, as asserted in Theorem~\ref{thm-main2}.


\begin{thebibliography}{10}

 
    
\bibitem{Bowen1971}
{R.~Bowen},
\newblock {\it Entropy for group endomorphisms and homogeneous spaces},
\newblock {\bf Trans. Amer. Math. Soc.}, 153:401--414, 1971.
       
 
\bibitem{EMS1977}
{R.~Edwards, K.~Millett, and D.~Sullivan}, 
\newblock {\it Foliations with all leaves compact}, 
\newblock {\bf Topology}, 16:13--32, 1977.


\bibitem{EpsteinVogt1978}
{D.B.A.~Epstein and E.~Vogt},
\newblock {\it A counterexample to the periodic orbit conjecture in codimension {$3$}},
\newblock {\bf Ann. of Math. (2)}, 108:539--552, 1978.

 
\bibitem{Ghys1995}
{{\'E.}~Ghys},
\newblock {\it Construction de champs de vecteurs sans orbite p\'eriodique (d'apr\`es {K}rystyna {K}uperberg)},
\newblock {S{\'e}minaire Bourbaki, Vol. 1993/94, Exp.\ No.\ 785},
\newblock {\bf Ast\'erisque}, Vol. 227: 283--307, 1995.

\bibitem{Harrison1988}
{J.C.~Harrison},
\newblock {\it {$C^2$} counterexamples to the {S}eifert conjecture},
\newblock {\bf Topology}, 27:249--278, 1988.

 


\bibitem{HR2016}
{S.~Hurder and A.~Rechtman},
\newblock {\it The dynamics of  generic Kuperberg flows},
\newblock  {\bf Ast\'erisque}, Vol. 377:1--250, 2016.

\bibitem{HKR2016}
{S.~Hurder and A.~Rechtman},
\newblock {\it Perspectives on    Kuperberg flows},
\newblock  {\bf submitted}, 2016;  {arXiv:1607.00731}.

 
\bibitem{Katok1980}
{A.~Katok},
\newblock {\it Lyapunov exponents, entropy and periodic orbits for diffeomorphisms},
\newblock {\bf Inst. Hautes \'Etudes Sci. Publ. Math.}, 51:137--173, 1980.

 
\bibitem{Kuperberg1994}
{K.~Kuperberg},
\newblock {\it A smooth counterexample to the {S}eifert conjecture},
\newblock {\bf Ann. of Math. (2)}, 140:723--732, 1994.

\bibitem{Kuperbergs1996}
{G.~Kuperberg and K.~Kuperberg},
\newblock {\it Generalized counterexamples to the {S}eifert conjecture},
\newblock {\bf Ann. of Math. (2)}, 144:239--268, 1996.


      
\bibitem{Matsumoto1995}
{S.~Matsumoto},
\newblock {\it K.{M}. {K}uperberg's {$C^\infty$} counterexample to the {S}eifert conjecture},
\newblock {\bf S\=ugaku}, Mathematical Society of Japan, 47:38--45, 1995.
\newblock {Translation:} {\bf Sugaku Expositions}, 11:39--49, 1998, 
\newblock {Amer. Math. Soc.} 

 \bibitem{Matsumoto2010}
{S.~Matsumoto}, 
\newblock {\it The unique ergodicity of equicontinuous laminations}, 
\newblock  {\bf Hokkaido Math. J.},   39:389--403, 2010.


\bibitem{PercelWilson1977}
{P.B.~Percell and F.W.~Wilson, Jr.},
\newblock {\it Plugging flows},
\newblock {\bf Trans. Amer. Math. Soc.}, 233:93--103, 1977.

 

 \bibitem{Rechtman2009}
{A.~Rechtman},
\newblock {\it Pi\`{e}ges dans la th\'{e}orie des feuilletages:  exemples et contre-exemples},
\newblock  {\bf Th\`{e}se}, l'Universit\'{e} de Lyon -- \'{E}cole Normale Sup\'{e}rieure de Lyon, 2009.

\bibitem{Schweitzer1974}
{P.A.~Schweitzer},
\newblock {\it Counterexamples to the {S}eifert conjecture and opening closed leaves of foliations},
\newblock  {\bf Ann. of Math. (2)}, 100:386--400, 1974.

\bibitem{Seifert1950}
{H.~Seifert},
\newblock {\it Closed integral curves in {$3$}-space and isotopic two-dimensional deformations},
\newblock {\bf Proc. Amer. Math. Soc.}, 1:287--302, 1950.

     
\bibitem{Sullivan1976}
{D.~Sullivan},
\newblock {\it A counterexample to the periodic orbit conjecture},
\newblock {\bf Inst. Hautes \'Etudes Sci. Publ. Math.}, 46:5--14, 1976.
     

 
\bibitem{Walters1982}
{P.~Walters},
\newblock {\bf An {I}ntroduction to {E}rgodic {T}heory},
\newblock {Graduate Texts in Math. Vol. 79}, Springer-Verlag, New York, 1982. 
    
\bibitem{Wilson1966}
{F.W.~Wilson, Jr.},
\newblock {\it On the minimal sets of non-singular vector fields},
\newblock {\bf Ann. of Math. (2)}, 84:529--536, 1966.


\end{thebibliography}
\end{document}